\documentclass{article}
\usepackage{amssymb}
\usepackage{graphicx}
\usepackage{amsmath}

\newtheorem{theorem}{Theorem}

\newtheorem{corollary}[theorem]{Corollary}

\newtheorem{example}[theorem]{Example}

\newtheorem{lemma}[theorem]{Lemma}

\newtheorem{proposition}[theorem]{Proposition}

\newtheorem{summary}[theorem]{Summary}

\renewcommand{\summary}[1]{\textbf{Summary}\ }
\begin{document}

\title{A Generalized Ito Formula}
\author{ Kenneth L. Kuttler \ and Li Ji \\
klkuttle@math.byu.edu\\
{\small Department of Mathematics}\\
{\small Brigham Young University}\\
{\small Provo, UT 84602} \\
liji84852418@gmail.com \\
{\small Institute for mathematics and its applications,}\\
{\small Minneapolis, MN, 55455}}
\maketitle

\begin{abstract}
An Ito formula is developed in a context consistent with the development of
abstract existence and uniqueness theorems for nonlinear stochastic partial
differential equations, which are singular or degenerate. This is a
generalization of an earlier Ito formula for Gelfand triples. After this, an existence theorem
is presented for some singular and degenerate stochastic equations followed by a few examples.
\end{abstract}

\vskip8pt


Classification numbers: 60G44,60H05,60H15,60H20,35R60,35M13

\section{Introduction}

The Ito formula describes $F\left( X\right) $ where
\begin{equation*}
X\left( t\right) =X_{0}+\int_{0}^{t}\phi \left( s\right) ds+\int_{0}^{t}\Phi
dW
\end{equation*}%
in which the last term is an Ito integral and $W$ is a Wiener process. The
above integral equation is the precise meaning for the stochastic
differential equation
\begin{equation*}
dX=\phi dt+\Phi dW,\ X\left( 0\right) =X_{0}.
\end{equation*}%
There are various forms for the Ito formula depending on where $X$ takes its
values. When $X$ has values in a separable Hilbert space and $F$ is
sufficiently smooth, the Ito formula takes the form
\begin{equation*}
F\left( t,X\left( t\right) \right) =F\left( 0,X_{0}\right)
+\int_{0}^{t}F_{X}\left( \cdot ,X\left( \cdot \right) \right) \Phi dW+
\end{equation*}%
\begin{equation*}
\int_{0}^{t}F_{t}\left( s,X\left( s\right) \right) +F_{X}\left( s,X\left(
s\right) \right) \phi \left( s\right) ds+\frac{1}{2}\int_{0}^{t}\left(
F_{XX}\left( s,X\left( s\right) \right) \Phi ,\Phi \right) _{\mathcal{L}%
_{2}\left( Q^{1/2}U,H\right) }ds
\end{equation*}%
In this formula, $\Phi $ is a stochastically square integrable function
having values in the Hilbert space of Hilbert Schmidt operators $\mathcal{L}%
_{2}\left( Q^{1/2}U,H\right) \ $where $Q$ is a nonnegative self adjoint
operator defined on a Hilbert space $U$.

In addition, there is a version of the Ito formula in the context of a
Gelfand triple of spaces
\begin{equation*}
V\subseteq H=H^{\prime }\subseteq V^{\prime }
\end{equation*}%
in which
\begin{equation}
X\left( t\right) =X_{0}+\int_{0}^{t}Y\left( s\right) ds+\int_{0}^{t}Z\left(
s\right) ds,  \label{6june1h}
\end{equation}%
the equation holding in $V^{\prime }$ for $t\in \left[ 0,T\right] $ almost
everywhere. In this case it is known that if for some $p>1$
\begin{equation*}
X\in L^{p}\left( \left[ 0,T\right] \times \Omega ,V\right) \cap L^{2}\left( %
\left[ 0,T\right] \times \Omega ,H\right) ,\ Y\in L^{p^{\prime }}\left( %
\left[ 0,T\right] \times \Omega ,V^{\prime }\right)
\end{equation*}%
\begin{equation*}
Z\in L^{2}\left( \left[ 0,T\right] \times \Omega ,\mathcal{L}_{2}\left(
Q^{1/2}U,H\right) \right)
\end{equation*}%
Then
\begin{equation*}
\left\vert X\left( t\right) \right\vert _{H}^{2}=\left\vert X_{0}\right\vert
^{2}+2\int_{0}^{t}\left\langle Y\left( s\right) ,\bar{X}\left( s\right)
\right\rangle ds+\int_{0}^{t}\left\Vert Z\left( s\right) \right\Vert _{%
\mathcal{L}_{2}\left( Q^{1/2}U,H\right) }^{2}ds+\mathcal{M}\left( t\right)
\end{equation*}%
where $\mathcal{M}\left( t\right) $ is a local martingale defined as a
stochastic integral, $\mathcal{M}\left( 0\right) =0$. Thus, one can obtain
the important estimate
\begin{eqnarray*}
E\left( \left\vert X\left( t\right) \right\vert _{H}^{2}\right) &=&E\left(
\left\vert X_{0}\right\vert ^{2}\right) +2E\left( \int_{0}^{t}\left\langle
Y\left( s\right) ,\bar{X}\left( s\right) \right\rangle ds\right) \\
&&+E\left( \int_{0}^{t}\left\Vert Z\left( s\right) \right\Vert _{\mathcal{L}%
_{2}\left( Q^{1/2}U,H\right) }^{2}ds\right)
\end{eqnarray*}%
A discussion of this formula and its applications is found in \cite{pre07}
it appears to be due to Krylov and is in Russian \cite{kr79}. This is a much
more difficult result. It is shown in this reference that this Ito formula
is the fundamental idea in developing general existence and uniqueness
theorems for nonlinear stochastic partial differential equations in the
context of variational formulations involving Gelfand triples. The formula
itself, without the stochastic terms, is fairly familiar to those who
formulate partial differential equations in this way, but it is much more
profound and difficult than the standard results for deterministic problems
because the presence of the stochastic integral causes a loss of weak time
derivatives. As is well known, the Wiener process is nowhere differentiable.
There are other major technical difficulties related to the minimal
assumption that $Z\in L^{2}\left( \left[ 0,T\right] \times \Omega ,\mathcal{L%
}_{2}\left( Q^{1/2}U,H\right) \right) $. These considerations require the
use of the Burkholder Davis Gundy inequality.

For deterministic evolution equations, an interesting generalization was the
step from evolution equations
\begin{equation*}
y^{\prime }+Ay=f
\end{equation*}%
to implicit or degenerate evolution equations%
\begin{equation*}
\left( By\right) ^{\prime }+Ay=f
\end{equation*}%
in which $B$ is an operator which may vanish. Since $B$ may fail to be one
to one, it may be impossible to consider such an equation as an evolution
equation. Instead it is called an implicit evolution equation or sometimes a
degenerate evolution equation. It could also happen that $B$ comes from some
sort of differential operator and may even be a Riesz map or as a special
case, the identity map on a Hilbert space in the context of a Gelfand triple.

In the case of deterministic equations, this was a very natural
generalization studied by many authors including Lions \cite{Lio69}, Brezis
\cite{bre68}, Showalter \cite{car76}, Bardos \cite{bar69}, \cite{Ken86}. and
many others. It led to interesting theorems including abstract existence and
uniqueness results for partial differential equations of mixed type, simple
ways to include systems of equations which involved coupling an elliptic
equation with a parabolic equation, and more transparent treatments of
equations like the porous media equation. If a good theory of implicit
stochastic equations can be obtained, many of the same interesting
applications will also have an extension to stochastic problems. The Ito
formula discussed above is a way to do integration by parts arguments for
stochastic evolution equations, and the version in this paper will provide
similar justification of integration by parts procedures for degenerate or
implicit stochastic equations. Thus many of the interesting deterministic
examples of the last forty years which are in terms of degenerate or partial
differential equations of mixed type will have generalizations to stochastic
versions.

In this paper, there will be a reflexive separable Banach space $V$ and a
separable Hilbert space $W,$ such that $V$ is dense in $W.$ Thus it is
possible to consider the following generalization of a Gelfand triple.
\begin{equation*}
V\subseteq W,\ \ W^{\prime }\subseteq V^{\prime },
\end{equation*}%
The usual pivot space $H$ is replaced with the pair $W,W^{\prime }$. It is
also assumed
\begin{equation}
BX\left( t\right) =BX_{0}+\int_{0}^{t}Y\left( s\right) ds+B\int_{0}^{t}ZdW,
\label{6octe1}
\end{equation}%
where it is known that
\begin{equation*}
X\in L^{p}\left( \left[ 0,T\right] \times \Omega ,V\right) ,\ BX\in
L^{2}\left( \left[ 0,T\right] \times \Omega ,W^{\prime }\right) ,\ Y\in
L^{p^{\prime }}\left( \left[ 0,T\right] \times \Omega ,V^{\prime }\right)
\end{equation*}%
\begin{equation*}
Z\in L^{2}\left( \left[ 0,T\right] \times \Omega ,\mathcal{L}_{2}\left(
Q^{1/2}U,W\right) \right)
\end{equation*}%
In terms of stochastic differential equations it is formally written as
\begin{equation*}
d\left( BX\right) =Ydt+BZdW,\ BX\left( 0\right) =BX_{0}.
\end{equation*}%
It will be assumed $B$ is a bounded nonnegative self adjoint operator which
maps $W$ to $W^{\prime }$. The case that $B$ is not one to one is included.
Then the Ito formula gives the justification for integration by parts
manipulations commonly used in the study of evolution equations.

It is necessary to have the stochastic part of \ref{6octe1} to vanish in
case $B=0,$ since otherwise, you might obtain an Ito integral equal to a
deterministic integral. However, the Ito integral will likely be nowhere
differentiable, due to this property which is possessed by the Wiener
process, \cite{str94}, \cite{str011} but the deterministic integral will
have a derivative a.e. Thus the above formula for $BX\left( t\right) $ is a
reasonable generalization of the case of evolution equations \ref{6june1h}.

\ When the formula for this more general situation is obtained, the more
standard result like one obtained in \cite{pre07} the context of a Gelfand
triple is recovered by letting $W=H$ and $B=I$.

To begin with, the paper considers some preliminary results and then the
proof of the Ito formula is presented. The techniques generalize those used
in \cite{pre07} to the situation where $V\subseteq W,W^{\prime }\subseteq
V^{\prime }$ instead of the more usual Gelfand triple. All spaces will be
assumed real and separable in the paper. Furthermore, there is the usual
filtration determined from increments of the Wiener process with respect to
which all martingale considerations are defined. This filtration is denoted
by $\mathcal{F}_{t}$ and it is assumed to be a normal filtration \cite{pre07}
so that each $\mathcal{F}_{t}$ is complete and $\mathcal{F}_{t+}=\cap _{s>t}%
\mathcal{F}_{s}=\mathcal{F}_{t}.$

In Section \ref{preliminary} we give a brief discussion of background
results. In Section \ref{integralequation} we give a fundamental equation
which will serve as the basis for the proof of the Ito formula. In Section %
\ref{mainestimate} a remarkable estimate is obtained along with some other
assertions. Section \ref{technicalsimplification} is devoted to obtaining a
technical simplification. It is this which allows us to consider the most
general initial conditions. Section \ref{itoformula} has the main result of
the paper.

\section{Preliminary results \label{preliminary}}

The entire presentation is based on the following fundamental lemma \cite%
{kar91}.

\begin{lemma}
\label{z23decl1g}Let $\Phi :\left[ 0,T\right] \times \Omega \rightarrow E,$
be $\mathcal{B}\left( \left[ 0,T\right] \right) \times \mathcal{F}$
measurable and suppose
\begin{equation*}
\Phi \in K\equiv L^{p}\left( \left[ 0,T\right] \times \Omega ;E\right) ,\
p\geq 1
\end{equation*}%
Then there exists a sequence of nested partitions, $\mathcal{P}_{k}\subseteq
\mathcal{P}_{k+1},$
\begin{equation*}
\mathcal{P}_{k}\equiv \left\{ t_{0}^{k},\cdots ,t_{m_{k}}^{k}\right\}
\end{equation*}%
such that the step functions given by
\begin{eqnarray*}
\Phi _{k}^{r}\left( t\right) &\equiv &\sum_{j=1}^{m_{k}}\Phi \left(
t_{j}^{k}\right) \mathcal{X}_{[t_{j-1}^{k},t_{j}^{k})}\left( t\right) \\
\Phi _{k}^{l}\left( t\right) &\equiv &\sum_{j=1}^{m_{k}}\Phi \left(
t_{j-1}^{k}\right) \mathcal{X}_{[t_{j-1}^{k},t_{j}^{k})}\left( t\right)
\end{eqnarray*}%
both converge to $\Phi $ in $K$ as $k\rightarrow \infty $ and
\begin{equation*}
\lim_{k\rightarrow \infty }\max \left\{ \left\vert
t_{j}^{k}-t_{j+1}^{k}\right\vert :j\in \left\{ 0,\cdots ,m_{k}\right\}
\right\} =0.
\end{equation*}%
Also, each $\Phi \left( t_{j}^{k}\right) ,\Phi \left( t_{j-1}^{k}\right) $
is in $L^{p}\left( \Omega ;E\right) $. One can also assume that $\Phi \left(
0\right) =0$. The mesh points $\left\{ t_{j}^{k}\right\} _{j=0}^{m_{k}}$ can
be chosen to miss a given set of measure zero.
\end{lemma}

There is also a known result on quadratic variation which we use later. \cite%
{chu83}

\begin{theorem}
\label{31mayt1h} Let $H$ be a Hilbert space and suppose $\left( M,\mathcal{F}%
_{t}\right) ,t\in \left[ 0,T\right] $ is a uniformly bounded continuous
martingale with values in $H$. Also let $\left\{ t_{k}^{n}\right\}
_{k=1}^{m_{n}}$ be a sequence of partitions satisfying
\begin{equation*}
\lim_{n\rightarrow \infty }\max \left\{ \left\vert
t_{i}^{n}-t_{i+1}^{n}\right\vert ,i=0,\cdots ,m_{n}\right\} =0,\ \left\{
t_{k}^{n}\right\} _{k=1}^{m_{n}}\subseteq \left\{ t_{k}^{n+1}\right\}
_{k=1}^{m_{n+1}}.
\end{equation*}%
Then
\begin{equation*}
\left[ M\right] \left( t\right) =\lim_{n\rightarrow \infty
}\sum_{i=0}^{m_{n}-1}\left\vert M\left( t\wedge t_{k+1}^{n}\right) -M\left(
t\wedge t_{k}^{n}\right) \right\vert _{H}^{2}
\end{equation*}%
the limit taking place in $L^{2}\left( \Omega \right) $. In case $M$ is just
a continuous local martingale, the above limit happens in probability.
\end{theorem}

In order to deal with the possibly degenerate operator $B,$ we have the
following interesting generalization of standard material involving inner
products.

\begin{lemma}
\label{7mayl1h}Suppose $V,W$ are separable Banach spaces, $W$ also a Hilbert
space such that $V$ is dense in $W$ and $B\in \mathcal{L}\left( W,W^{\prime
}\right) $ satisfies%
\begin{equation*}
\left\langle Bx,x\right\rangle \geq 0,\ \left\langle Bx,y\right\rangle
=\left\langle By,x\right\rangle ,B\neq 0.
\end{equation*}%
Then there exists a countable set $\left\{ e_{i}\right\} $ of vectors in $V$
such that
\begin{equation*}
\left\langle Be_{i},e_{j}\right\rangle =\delta _{ij}
\end{equation*}%
and for each $x\in W,$%
\begin{equation*}
\left\langle Bx,x\right\rangle =\sum_{i=1}^{\infty }\left\vert \left\langle
Bx,e_{i}\right\rangle \right\vert ^{2},
\end{equation*}%
and also
\begin{equation*}
Bx=\sum_{i=1}^{\infty }\left\langle Bx,e_{i}\right\rangle Be_{i},
\end{equation*}%
the series converging in $W^{\prime }$.
\end{lemma}

\textbf{Proof:\ } Let $\left\{ g_{k}\right\} _{k=1}^{\infty }$ be linearly
independent vectors of $V$ whose span is dense in $V$. This is possible
because $V$ is separable. Let $n_{1}$ be the first index such that $%
\left\langle Bg_{n_{1}},g_{n_{1}}\right\rangle \neq 0.$

\textbf{Claim:\ }If there is no such index, then $B=0.$

\textbf{Proof of claim:\ }First note that if $\left\langle Bg,g\right\rangle
=0,$ then%
\begin{equation*}
\left\vert \left\langle Bg,x\right\rangle \right\vert \leq \left\vert
\left\langle Bg,g\right\rangle \right\vert ^{1/2}\left\vert \left\langle
Bx,x\right\rangle \right\vert ^{1/2}=0
\end{equation*}%
and so $Bg=0$. Therefore, if $x$ is given, you could take $x_{k}$ in the
span of $\left\{ g_{1},\cdots ,g_{k}\right\} $ such that $\left\Vert
x_{k}-x\right\Vert _{W}\rightarrow 0$. Then
\begin{equation*}
\left\vert \left\langle Bx,y\right\rangle \right\vert =\lim_{k\rightarrow
\infty }\left\vert \left\langle Bx_{k},y\right\rangle \right\vert \leq
\lim_{k\rightarrow \infty }\left\langle Bx_{k},x_{k}\right\rangle
^{1/2}\left\langle By,y\right\rangle ^{1/2}=0
\end{equation*}%
because $Bx_{k}$ is zero, being the sum of scalars times $Bg_{i}$ for
finitely many $i$. Since $y$ is arbitrary, this shows $Bx=0$.

Thus assume there is such a first index. Let
\begin{equation*}
e_{1}\equiv \frac{g_{n_{1}}}{\left\langle Bg_{n_{1}},g_{n_{1}}\right\rangle
^{1/2}}
\end{equation*}%
Then $\left\langle Be_{1},e_{1}\right\rangle =1.$ Now if you have
constructed $e_{j}$ for $j\leq k,$
\begin{equation*}
e_{j}\in \mbox{span}\left( g_{n_{1}},\cdots ,g_{n_{k}}\right) ,\
\left\langle Be_{i},e_{j}\right\rangle =\delta _{ij},
\end{equation*}%
$g_{n_{j+1}}$ being the first for which%
\begin{equation*}
\left\langle Bg_{n_{j+1}}-\sum_{i=1}^{j}\left\langle
Bg_{n_{j+1}},e_{i}\right\rangle
Be_{i},g_{n_{j+1}}-\sum_{i=1}^{j}\left\langle Bg_{nj},e_{i}\right\rangle
e_{i}\right\rangle \neq 0,
\end{equation*}%
and%
\begin{equation*}
\mbox{span}\left( g_{n_{1}},\cdots ,g_{n_{k}}\right) =\mbox{span}\left(
e_{1},\cdots ,e_{k}\right) ,
\end{equation*}%
let $g_{n_{k+1}}$ be such that $g_{n_{k+1}}$ is the first in the list $%
\left\{ g_{n_{j}}\right\} _{j=1}^{\infty }$ such that
\begin{equation*}
\left\langle Bg_{n_{k+1}}-\sum_{i=1}^{k}\left\langle
Bg_{n_{k+1}},e_{i}\right\rangle
Be_{i},g_{n_{k+1}}-\sum_{i=1}^{k}\left\langle
Bg_{n_{k+1}},e_{i}\right\rangle e_{i}\right\rangle \neq 0
\end{equation*}

\textbf{Claim:\ }If there is no such first $g_{n_{k+1}},$ then $B\left( %
\mbox{span}\left( e_{i},\cdots ,e_{k}\right) \right) =BW$ so in this case, $%
\left\{ Be_{i}\right\} _{i=1}^{k}$ is actually a basis for $BW$.

\textbf{Proof:\ }Let $x\in W$. Let $x_{r}\in \mbox{span}\left( g_{1},\cdots
,g_{r}\right) ,r>n_{k}$ such that $\lim_{r\rightarrow \infty }x_{r}=x$ in $W$%
. Then
\begin{equation}
x_{r}=\sum_{i=1}^{k}c_{i}^{r}e_{i}+\sum_{i\notin \left\{ n_{1},\cdots
,n_{k}\right\} }^{r}d_{i}^{r}g_{i}\equiv y_{r}+z_{r}  \label{13mare7g}
\end{equation}%
If $l\notin \left\{ n_{1},\cdots ,n_{k}\right\} ,$ then by the construction
and the above assumption, for some $j\leq k$
\begin{equation*}
\left\langle Bg_{l}-\sum_{i=1}^{j}\left\langle Bg_{l},e_{i}\right\rangle
Be_{i},g_{l}-\sum_{i=1}^{j}\left\langle Bg_{l},e_{i}\right\rangle
e_{i}\right\rangle =0
\end{equation*}%
The reasoning is as follows. If $l\leq k$ and if the above is nonzero for
all $j\leq k,$ then $l$ would have been chosen but it wasn't. Thus in this
case that $l\leq k,$ there exists $j$ such that
\begin{equation*}
Bg_{l}=\sum_{i=1}^{j}\left\langle Bg_{l},e_{i}\right\rangle Be_{i}
\end{equation*}%
If $l>n_{k},$ then by assumption, the above is never nonzero for $j=k$.
Thus, in any case, it follows that for each $l\notin \left\{ n_{1},\cdots
,n_{k}\right\} ,$%
\begin{equation*}
Bg_{l}\in B\left( \mbox{span}\left( e_{i},\cdots ,e_{k}\right) \right) .
\end{equation*}

Now it follows from \ref{13mare7g} that
\begin{eqnarray*}
Bx_{r} &=&\sum_{i=1}^{k}c_{i}^{r}Be_{i}+\sum_{i\notin \left\{ n_{1},\cdots
,n_{k}\right\} }^{r}d_{i}^{r}Bg_{i} \\
&=&\sum_{i=1}^{k}c_{i}^{r}Be_{i}+\sum_{i\notin \left\{ n_{1},\cdots
,n_{k}\right\} }^{r}d_{i}^{r}\sum_{j=1}^{k}c_{j}^{i}Be_{j}
\end{eqnarray*}%
and so $Bx_{r}\in B\left( \mbox{span}\left( e_{i},\cdots ,e_{k}\right)
\right) .$ Then $Bx=\lim_{r\rightarrow \infty }Bx_{r}=\lim_{r\rightarrow
\infty }By_{r}$ where $y_{r}\in \mbox{span}\left( e_{i},\cdots ,e_{k}\right)
$. Say
\begin{equation*}
Bx_{r}=\sum_{i=1}^{k}a_{i}^{r}Be_{i}
\end{equation*}%
It follows easily that $\left\langle Bx_{r},e_{j}\right\rangle =a_{j}^{r}.$
(Act on $e_{j}$ by both sides and use $\left\langle
Be_{i},e_{j}\right\rangle =\delta _{ij}.$) Now since $x_{r}$ is bounded, it
follows that these $a_{j}^{r}$ are also bounded. Hence, defining $%
y_{r}\equiv \sum_{i=1}^{k}a_{i}^{r}e_{i},$ it follows that $y_{r}$ is
bounded in $\mbox{span}\left( e_{i},\cdots ,e_{k}\right) $ and so, there
exists a subsequence, still denoted by $r$ such that $y_{r}\rightarrow y\in %
\mbox{span}\left( e_{i},\cdots ,e_{k}\right) $. Therefore, $%
Bx=\lim_{r\rightarrow \infty }By_{r}=By$. In other words, $BW=B\left( %
\mbox{span}\left( e_{i},\cdots ,e_{k}\right) \right) $ as claimed. This
proves the claim.

If this happens, the process being described stops. You have found what is
desired which has only finitely many vectors involved.

As long as the process does not stop, let
\begin{equation*}
e_{k+1}\equiv \frac{g_{n_{k+1}}-\sum_{i=1}^{k}\left\langle
Bg_{n_{k+1}},e_{i}\right\rangle e_{i}}{\left\langle B\left(
g_{n_{k+1}}-\sum_{i=1}^{k}\left\langle Bg_{n_{k+1}},e_{i}\right\rangle
e_{i}\right) ,g_{n_{k+1}}-\sum_{i=1}^{k}\left\langle
Bg_{n_{k+1}},e_{i}\right\rangle e_{i}\right\rangle ^{1/2}}
\end{equation*}%
Thus, as in the usual argument for the Gram Schmidt process, $\left\langle
Be_{i},e_{j}\right\rangle =\delta _{ij}$ for $i,j\leq k$.

Consider
\begin{equation}
\left\langle Bg_{p}-B\left( \sum_{i=1}^{k}\left\langle
Bg_{p},e_{i}\right\rangle e_{i}\right) ,g_{p}-\sum_{i=1}^{k}\left\langle
Bg_{p},e_{i}\right\rangle e_{i}\right\rangle  \label{17may1f}
\end{equation}%
If $p$ is never one of the $n_{k},$ then there exists $k$ such that $p\in
\left( n_{k},n_{k+1}\right) $ so \ref{17may1f} equals 0. If $p=n_{k}$ for
some $k,$ then from the construction, $g_{n_{k}}=g_{p}\in \mbox{span}\left(
e_{1},\cdots ,e_{k}\right) $ and therefore,%
\begin{equation*}
g_{p}=\sum_{j=1}^{k}a_{j}e_{j}
\end{equation*}%
which requires easily that
\begin{equation*}
Bg_{p}=\sum_{i=1}^{k}\left\langle Bg_{p},e_{i}\right\rangle Be_{i},
\end{equation*}%
and \ref{17may1f} equals 0, the above holding for all $k$ large enough. It
follows that for any $x\in \mbox{span}\left( \left\{ g_{k}\right\}
_{k=1}^{\infty }\right) ,$ (finite linear combination of vectors in $\left\{
g_{k}\right\} _{k=1}^{\infty }$).
\begin{equation}
Bx=\sum_{i=1}^{\infty }\left\langle Bx,e_{i}\right\rangle Be_{i}
\label{13mare4g}
\end{equation}%
because for all $k$ large enough,
\begin{equation*}
Bx=\sum_{i=1}^{k}\left\langle Bx,e_{i}\right\rangle Be_{i}
\end{equation*}%
Also note that for such $x\in \mbox{span}\left( \left\{ g_{k}\right\}
_{k=1}^{\infty }\right) ,$%
\begin{eqnarray*}
\left\langle Bx,x\right\rangle &=&\left\langle \sum_{i=1}^{k}\left\langle
Bx,e_{i}\right\rangle Be_{i},x\right\rangle =\sum_{i=1}^{k}\left\langle
Bx,e_{i}\right\rangle \left\langle Bx,e_{i}\right\rangle \\
&=&\sum_{i=1}^{k}\left\vert \left\langle Bx,e_{i}\right\rangle \right\vert
^{2}=\sum_{i=1}^{\infty }\left\vert \left\langle Bx,e_{i}\right\rangle
\right\vert ^{2}
\end{eqnarray*}%
Now for $x$ arbitrary, let $x_{k}\rightarrow x$ in $W$ where $x_{k}\in %
\mbox{span}\left( \left\{ g_{k}\right\} _{k=1}^{\infty }\right) .$ Then by
Fatou's lemma,

\begin{eqnarray}
\sum_{i=1}^{\infty }\left\vert \left\langle Bx,e_{i}\right\rangle
\right\vert ^{2} &\leq &\lim \inf_{k\rightarrow \infty }\sum_{i=1}^{\infty
}\left\vert \left\langle Bx_{k},e_{i}\right\rangle \right\vert ^{2}  \notag
\\
&=&\lim \inf_{k\rightarrow \infty }\left\langle Bx_{k},x_{k}\right\rangle
=\left\langle Bx,x\right\rangle  \label{13mare6g} \\
&\leq &\left\Vert Bx\right\Vert _{W^{\prime }}\left\Vert x\right\Vert
_{W}^{2}\leq \left\Vert B\right\Vert \left\Vert x\right\Vert _{W}^{2}  \notag
\end{eqnarray}%
Thus the series on the left converges. Then also, from the above inequality,
\begin{equation*}
\left\vert \left\langle \sum_{i=p}^{q}\left\langle Bx,e_{i}\right\rangle
Be_{i},y\right\rangle \right\vert \leq \sum_{i=p}^{q}\left\vert \left\langle
Bx,e_{i}\right\rangle \right\vert \left\vert \left\langle
Be_{i},y\right\rangle \right\vert
\end{equation*}%
\begin{eqnarray*}
&\leq &\left( \sum_{i=p}^{q}\left\vert \left\langle Bx,e_{i}\right\rangle
\right\vert ^{2}\right) ^{1/2}\left( \sum_{i=p}^{q}\left\vert \left\langle
By,e_{i}\right\rangle \right\vert ^{2}\right) ^{1/2} \\
&\leq &\left( \sum_{i=p}^{q}\left\vert \left\langle Bx,e_{i}\right\rangle
\right\vert ^{2}\right) ^{1/2}\left( \sum_{i=1}^{\infty }\left\vert
\left\langle By,e_{i}\right\rangle \right\vert ^{2}\right) ^{1/2}
\end{eqnarray*}%
\begin{equation*}
\leq \left( \sum_{i=p}^{q}\left\vert \left\langle Bx,e_{i}\right\rangle
\right\vert ^{2}\right) ^{1/2}\left( \left\Vert B\right\Vert \left\Vert
y\right\Vert _{W}^{2}\right) ^{1/2}\leq \left( \sum_{i=p}^{q}\left\vert
\left\langle Bx,e_{i}\right\rangle \right\vert ^{2}\right) ^{1/2}\left\Vert
B\right\Vert ^{1/2}\left\Vert y\right\Vert
\end{equation*}%
It follows that
\begin{equation}
\sum_{i=1}^{\infty }\left\langle Bx,e_{i}\right\rangle Be_{i}
\label{13mare2g}
\end{equation}%
converges in $W^{\prime }$ because it was just shown that
\begin{equation*}
\left\Vert \sum_{i=p}^{q}\left\langle Bx,e_{i}\right\rangle
Be_{i}\right\Vert _{W^{\prime }}\leq \left( \sum_{i=p}^{q}\left\vert
\left\langle Bx,e_{i}\right\rangle \right\vert ^{2}\right) ^{1/2}\left\Vert
B\right\Vert ^{1/2}
\end{equation*}%
and so the partial sums of the series \ref{13mare2g} constitute a Cauchy
sequence in $W^{\prime }$. Also, the above estimate shows that
\begin{equation}
\left\Vert \sum_{i=1}^{\infty }\left\langle Bx,e_{i}\right\rangle
Be_{i}\right\Vert _{W^{\prime }}\leq \left( \sum_{i=1}^{\infty }\left\vert
\left\langle Bx,e_{i}\right\rangle \right\vert ^{2}\right) ^{1/2}\left\Vert
B\right\Vert ^{1/2}  \label{13mare5g}
\end{equation}

Now for $x$ arbitrary, let $x_{k}\in \mbox{span}\left( \left\{ g_{j}\right\}
_{j=1}^{\infty }\right) $ and $x_{k}\rightarrow x$ in $W.$ Then for a fixed $%
k$ large enough,%
\begin{equation*}
\left\Vert Bx-\sum_{i=1}^{\infty }\left\langle Bx,e_{i}\right\rangle
Be_{i}\right\Vert \leq \left\Vert Bx-Bx_{k}\right\Vert
\end{equation*}%
\begin{equation*}
+\left\Vert Bx_{k}-\sum_{i=1}^{\infty }\left\langle
Bx_{k},e_{i}\right\rangle Be_{i}\right\Vert +\left\Vert \sum_{i=1}^{\infty
}\left\langle Bx_{k},e_{i}\right\rangle Be_{i}-\sum_{i=1}^{\infty
}\left\langle Bx,e_{i}\right\rangle Be_{i}\right\Vert
\end{equation*}%
\begin{equation*}
\leq \varepsilon +\left\Vert \sum_{i=1}^{\infty }\left\langle B\left(
x_{k}-x\right) ,e_{i}\right\rangle Be_{i}\right\Vert ,
\end{equation*}%
the middle term equaling 0 by \ref{13mare4g}. From \ref{13mare5g} and \ref%
{13mare6g},
\begin{eqnarray*}
&\leq &\varepsilon +\left\Vert B\right\Vert ^{1/2}\left( \sum_{i=1}^{\infty
}\left\vert \left\langle B\left( x_{k}-x\right) ,e_{i}\right\rangle
\right\vert ^{2}\right) ^{1/2} \\
&\leq &\varepsilon +\left\Vert B\right\Vert ^{1/2}\left\langle B\left(
x_{k}-x\right) ,x_{k}-x\right\rangle ^{1/2}<2\varepsilon
\end{eqnarray*}%
whenever $k$ is large enough. Therefore,
\begin{equation*}
Bx=\sum_{i=1}^{\infty }\left\langle Bx,e_{i}\right\rangle Be_{i}
\end{equation*}%
in $W^{\prime }$. It follows that
\begin{equation*}
\left\langle Bx,x\right\rangle =\lim_{k\rightarrow \infty }\left\langle
\sum_{i=1}^{k}\left\langle Bx,e_{i}\right\rangle Be_{i},x\right\rangle
=\lim_{k\rightarrow \infty }\sum_{i=1}^{k}\left\vert \left\langle
Bx,e_{i}\right\rangle \right\vert ^{2}\equiv \sum_{i=1}^{\infty }\left\vert
\left\langle Bx,e_{i}\right\rangle \right\vert ^{2}\ \blacksquare
\end{equation*}

The details of the definition of the stochastic integral are in \cite{pre07},%
\cite{dap92}. For completeness, here is a short summary. Consider the
following diagram in which $J$ is a one to one Hilbert Schmidt operator and $%
Q$ is a nonnegative and self adjoint operator defined on the Hilbert space $%
U $.
\begin{equation*}
\begin{array}{cccl}
&  &  & U \\
&  &  &
\begin{array}{cc}
\downarrow & Q^{1/2}%
\end{array}
\\
U_{1} & \supseteq JQ^{1/2}U & \underset{1-1}{\overset{J}{\leftarrow }} &
Q^{1/2}U \\
&
\begin{array}{cc}
&  \\
\Phi _{n} & \searrow%
\end{array}
&  &
\begin{array}{cc}
\downarrow & \Phi%
\end{array}
\\
&  &  & W%
\end{array}%
\end{equation*}%
The idea is to define $\int_{0}^{t}\Phi dW$ where $\Phi \in L^{2}\left( %
\left[ 0,T\right] \times \Omega ;\mathcal{L}_{2}\left( Q^{1/2}U,W\right)
\right) $, $\mathcal{L}_{2}\left( Q^{1/2}U,W\right) $ being the Hilbert
Schmidt operators mapping $Q^{1/2}U$ to $W$ and $J$ a Hilbert Schmidt
operator. Here $W\left( t\right) $ is the process
\begin{equation*}
W\left( t\right) =\sum_{i=1}^{\infty }\psi _{i}\left( t\right) Jg_{i}\text{
in }U_{1},
\end{equation*}%
where the $\psi _{i}\left( t\right) $ are real, independent Wiener
processes. It is a $Q_{1}$ Wiener process on $U_{1}$ for $Q_{1}=JJ^{\ast }$.
To get $\int_{0}^{t}\Phi dW,$ $\Phi \circ J^{-1}$ was approximated by a
sequence of elementary functions, $\left\{ \Phi _{n}\right\} ,$ adapted step
functions having finitely many values in $\mathcal{L}\left( U_{1},W\right) .$
Then the stochastic integral was defined in the usual way. For%
\begin{equation*}
\Phi _{n}\left( t\right) =\sum_{i=0}^{m-1}\phi _{i}\mathcal{X}%
_{[t_{i},t_{i+1})}\left( t\right) ,\ \phi _{i}\text{ being }\mathcal{F}%
_{t_{i}}\text{ measurable,}
\end{equation*}%
\begin{equation*}
\int_{0}^{t}\Phi _{n}dW\equiv \sum_{i=0}^{m-1}\phi _{i}\left( W\left(
t\wedge t_{i+1}\right) -W\left( t\wedge t_{i}\right) \right) .
\end{equation*}%
Then it is shown that this sequence of processes converges in $L^{2}\left(
\Omega ,W\right) $ and
\begin{equation*}
\int_{0}^{t}\Phi dW\equiv \lim_{n\rightarrow \infty }\int_{0}^{t}\Phi _{n}dW
\end{equation*}%
It can be shown that this integral $\int_{0}^{t}\Phi dW$ satisfies the Ito
isometry and is independent of the choice of $U_{1}$ and $J$.

In all that follows, $Q$ will be a nonnegative self adjoint operator defined
on a separable Hilbert space $U.$ Also $Z$ will be progressively measurable
and in $L^{2}\left( \left[ 0,T\right] \times \Omega ,\mathcal{L}_{2}\left(
Q^{1/2}U,W\right) \right) $ while $J:Q^{1/2}U\rightarrow U_{1}$ will be a
one to one Hilbert Schmidt operator.

Now here is a technical result which will be needed.$\ $This is a technical
application of the above description of the stochastic integral.

\begin{theorem}
\label{z19mayt1h}Let $Z$ be progressively measurable and in
\begin{equation*}
L^{2}\left( \left[ 0,T\right] \times \Omega ,\mathcal{L}_{2}\left(
Q^{1/2}U,W\right) \right) .
\end{equation*}%
Also suppose $P$ is progressively measurable and in $L^{2}\left( \left[ 0,T%
\right] \times \Omega ,W^{\prime }\right) $. Let $\left\{ t_{j}^{n}\right\}
_{j=0}^{m_{n}}$ be a sequence of partitions of the sort in Lemma \ref%
{z23decl1g} such that if
\begin{equation*}
P_{n}\left( t\right) \equiv \sum_{j=0}^{m_{n}-1}P\left( t_{j}^{n}\right)
\mathcal{X}_{[t_{j}^{n},t_{j+1}^{n})}\left( t\right) \equiv P_{n}^{l}\left(
t\right)
\end{equation*}%
then $P_{n}\rightarrow P$ in $L^{2}\left( \left[ 0,T\right] \times \Omega
,W\right) .$ Then the expression
\begin{equation}
\sum_{j=0}^{m_{n}-1}\left\langle P\left( t_{j}^{n}\right)
,\int_{t_{j}^{n}\wedge t}^{t_{j+1}^{n}\wedge t}ZdW\right\rangle
\label{z2june6h}
\end{equation}%
is a local martingale which can be written as a stochastic integral in the
form
\begin{equation*}
\int_{0}^{t}\left( Z\circ J^{-1}\right) ^{\ast }P_{n}^{l}\circ JdW
\end{equation*}
\end{theorem}

\textbf{Proof: }Note that $P_{n}$ is right continuous and progressively
measurable. Thus one can define the stopping time
\begin{equation*}
\tau _{n}^{p}\equiv \inf \left\{ t:\left\Vert P_{n}\left( t\right)
\right\Vert _{W}>p\right\} ,
\end{equation*}%
the first hitting time of an open set. We need the formula in \ref{z2june6h}
as a stochastic integral. First note that $W$ has values in $U_{1}$.

Consider one of the terms of the sum more simply as
\begin{equation*}
\left\langle P\left( a\right) ,\int_{a}^{b}ZdW\right\rangle ,\
a=t_{k}^{n}\wedge t,\ b=t_{k+1}^{n}\wedge t.
\end{equation*}%
Then from the definition of the integral, let $Z_{n}$ be a sequence of
elementary functions converging to $Z\circ J^{-1}$ in $L^{2}\left( \left[ a,b%
\right] \times \Omega ,\mathcal{L}_{2}\left( JQ^{1/2}U,W\right) \right) $
and
\begin{equation*}
\left\Vert \int_{a}^{t}ZdW-\int_{a}^{t}Z_{n}dW\right\Vert _{L^{2}\left(
\Omega ,W\right) }\rightarrow 0
\end{equation*}%
Using a maximal inequality and the fact that the two integrals are
martingales along with the Borel Cantelli lemma, there exists a set of
measure 0 $N$ such that for $\omega \notin N,$ the convergence of a suitable
subsequence of these integrals, still denoted by $n$, is uniform for $t\in $
$\left[ a,b\right] $. It follows that for such $\omega ,$%
\begin{equation}
\left\langle P\left( a\right) ,\int_{a}^{t}ZdW\right\rangle
=\lim_{n\rightarrow \infty }\left\langle P\left( a\right)
,\int_{a}^{t}Z_{n}dW\right\rangle .  \label{5june1h}
\end{equation}%
Say
\begin{equation*}
Z_{n}\left( u\right) =\sum_{k=0}^{m_{n}-1}Z_{k}^{n}\mathcal{X}%
_{[t_{k}^{n},t_{k+1}^{n})}\left( u\right)
\end{equation*}%
where $Z_{k}^{n}$ has finitely many values in $\mathcal{L}\left(
U_{1},W\right) _{0},$ the restrictions of maps in $\mathcal{L}\left(
U_{1},W\right) $ to $JQ^{1/2}U,$ and the $t_{k}^{n}$ refer to a partition of
$\left[ a,b\right] $. Then the product on the right in \-\ref{5june1h} is of
the form%
\begin{equation*}
\sum_{k=0}^{m_{n}-1}\left\langle P\left( a\right) ,Z_{k}^{n}\left( W\left(
t\wedge t_{k+1}^{n}\right) -W\left( t\wedge t_{k}^{n}\right) \right)
\right\rangle _{W^{\prime },W}
\end{equation*}%
Note that it makes sense because $Z_{k}^{n}$ is the restriction to $J\left(
Q^{1/2}U\right) $ of a map from $U_{1}$ to $W.$ Thus the above equals
\begin{equation*}
=\sum_{k=0}^{m_{n}-1}\left\langle P\left( a\right) ,Z_{k}^{n}\left( W\left(
t\wedge t_{k+1}^{n}\right) -W\left( t\wedge t_{k}^{n}\right) \right)
\right\rangle _{W^{\prime },W}
\end{equation*}%
\begin{eqnarray*}
&=&\sum_{k=0}^{m_{n}-1}\left\langle \left( Z_{k}^{n}\right) ^{\ast }P\left(
a\right) ,\left( W\left( t\wedge t_{k+1}^{n}\right) -W\left( t\wedge
t_{k}^{n}\right) \right) \right\rangle _{U_{1}^{\prime },U_{1}} \\
&=&\sum_{k=0}^{m_{n}-1}\left( Z_{k}^{n}\right) ^{\ast }P\left( a\right)
\left( W\left( t\wedge t_{k+1}^{n}\right) -W\left( t\wedge t_{k}^{n}\right)
\right) \\
&=&\int_{a}^{t}Z_{n}^{\ast }P\left( a\right) dW
\end{eqnarray*}%
Note that the restriction of $\left( Z_{n}\right) ^{\ast }P\left( a\right) $
is in
\begin{equation*}
\mathcal{L}\left( U_{1},\mathbb{R}\right) _{0}\subseteq \mathcal{L}%
_{2}\left( JQ^{1/2}U,\mathbb{R}\right) .
\end{equation*}%
Recall also that the space on the left is dense in the one on the right. Now
let $\left\{ g_{i}\right\} $ be an orthonormal basis for $Q^{1/2}U,$ so that
$\left\{ Jg_{i}\right\} $ is an orthonormal basis for $JQ^{1/2}U.$ Then
\begin{equation*}
\sum_{i=1}^{\infty }\left\vert \left( \left( Z_{n}\right) ^{\ast }P\left(
a\right) -\left( Z\circ J^{-1}\right) ^{\ast }P\left( a\right) \right)
\left( Jg_{i}\right) \right\vert ^{2}
\end{equation*}%
\begin{equation*}
=\sum_{i=1}^{\infty }\left\vert \left\langle P\left( a\right) ,\left(
Z_{n}-Z\circ J^{-1}\right) \left( Jg_{i}\right) \right\rangle \right\vert
^{2}
\end{equation*}%
\begin{eqnarray*}
&\leq &\sum_{i=1}^{\infty }\left\Vert P\left( a\right) \right\Vert
^{2}\left\Vert \left( Z_{n}-Z\circ J^{-1}\right) \left( Jg_{i}\right)
\right\Vert ^{2} \\
&=&\left\Vert P\left( a\right) \right\Vert ^{2}\left\Vert Z_{n}-Z\circ
J^{-1}\right\Vert _{\mathcal{L}_{2}\left( JQ^{1/2}U,W\right) }^{2}
\end{eqnarray*}%
When integrated over $\left[ a,b\right] \times \Omega ,$ it is given that
this converges to 0 as $n\rightarrow \infty $, assuming that $\left\Vert
P\left( a\right) \right\Vert \in L^{\infty }\left( \Omega \right) ,$ which
is assumed for now. It follows that
\begin{equation*}
Z_{n}^{\ast }P\left( a\right) \rightarrow \left( Z\circ J^{-1}\right) ^{\ast
}P\left( a\right)
\end{equation*}%
in $L^{2}\left( \left[ a,b\right] \times \Omega ,\mathcal{L}_{2}\left(
JQ^{1/2}U,\mathbb{R}\right) \right) .$ Writing this differently, it says
\begin{equation*}
Z_{n}^{\ast }P\left( a\right) \rightarrow \left( \left( Z\circ J^{-1}\right)
^{\ast }P\left( a\right) \circ J\right) \circ J^{-1}\text{ in }L^{2}\left( %
\left[ a,b\right] \times \Omega ,\mathcal{L}_{2}\left( JQ^{1/2}U,\mathbb{R}%
\right) \right)
\end{equation*}%
It follows from the definition of the integral that the Ito integrals
converge. Therefore,
\begin{equation*}
\left\langle P\left( a\right) ,\int_{a}^{t}ZdW\right\rangle
=\int_{a}^{t}\left( Z\circ J^{-1}\right) ^{\ast }P\left( a\right) \circ JdW
\end{equation*}%
The term on the right is a martingale.

Next it is necessary to drop the assumption that $\left\Vert P\left(
a\right) \right\Vert \in L^{\infty }\left( \Omega \right) $. This involves
the above stopping time. From localization,
\begin{eqnarray*}
\left\langle P\left( a\right) ,\int_{a\wedge \tau _{p}^{n}}^{t\wedge \tau
_{p}^{n}}ZdW\right\rangle &=&\left\langle P\left( a\right) ,\int_{a}^{t}%
\mathcal{X}_{\left[ 0,\tau _{p}^{n}\right] }ZdW\right\rangle \\
&=&\int_{a}^{t}\left( \mathcal{X}_{\left[ 0,\tau _{p}^{n}\right] }Z\circ
J^{-1}\right) ^{\ast }P\left( a\right) \circ JdW \\
&=&\int_{a\wedge \tau _{p}^{n}}^{t\wedge \tau _{p}^{n}}\left( Z\circ
J^{-1}\right) ^{\ast }P\left( a\right) \circ JdW
\end{eqnarray*}%
Then it follows that, using the stopping time,%
\begin{equation*}
\sum_{j=0}^{m_{n}-1}\left\langle P\left( t_{j}^{n}\right)
,\int_{t_{j}^{n}\wedge t\wedge \tau _{p}^{n}}^{t_{j+1}^{n}\wedge t\wedge
\tau _{p}^{n}}ZdW\right\rangle =\int_{0}^{t\wedge \tau _{p}^{n}}\left(
Z\circ J^{-1}\right) ^{\ast }P_{n}^{l}\circ JdW
\end{equation*}%
where $P_{n}^{l}$ is the step function
\begin{equation*}
P_{n}^{l}\left( t\right) =\sum_{k=0}^{m_{n}-1}P\left( t_{k}^{n}\right)
\mathcal{X}_{[t_{k}^{n},t_{k+1}^{n})}\left( t\right) .
\end{equation*}%
Thus the given sum equals the local martingale
\begin{equation*}
\int_{0}^{t}\left( Z\circ J^{-1}\right) ^{\ast }P_{n}^{l}\circ JdW.\ \
\blacksquare
\end{equation*}

The original formula does not depend on $J$ and so the same is true of this
last expression although it does not look like it. The unaesthetic
appearance of the above integral can be improved, but such an effort is of
no significance in what follows.

The next question is whether the above stochastic integral converges as $%
n\rightarrow \infty $ in some sense to an integral
\begin{equation}
\int_{0}^{t}\left( Z\circ J^{-1}\right) ^{\ast }P\circ JdW.  \label{z2june7h}
\end{equation}

The problem is that the integrand is not known to be in $L^{2}\left( \left[
0,T\right] \times \Omega ;\mathcal{L}_{2}\left( Q^{1/2}U,\mathbb{R}\right)
\right) $. It would be useful to define a stopping time
\begin{equation}
\tau _{n}\equiv \inf \left\{ t:\left\Vert P\left( t\right) \right\Vert
_{W^{\prime }}>n\right\}  \label{z2june9h}
\end{equation}%
because then, you could localize and define the integral in \ref{z2june7h}
as a local martingale. However, to do this would require the stopping time
to make sense. It is not known that $P$ is continuous or even right
continuous. Therefore, we need other assumptions.

\begin{lemma}
Suppose $t\rightarrow P\left( t\right) $ is weakly continuous into $%
W^{\prime }$ for a.e. $\omega $ and that $P$ is adapted. Then $\tau _{n}$
described in \ref{z2june9h} is well defined. It also satisfies $%
\lim_{n\rightarrow \infty }\tau _{n}=\infty $.
\end{lemma}

\textbf{Proof:\ }Let $O\equiv \left\{ y\in W:\left\Vert y\right\Vert
_{W^{\prime }}>n\right\} .$ Then the complement of $O$ is a closed convex
set. It follows that $O^{C}$ is also weakly closed. Hence $O$ must be weakly
open. Now $t\rightarrow P\left( t\right) $ is adapted as a function mapping
into the topological space consisting of $W^{\prime }$ with the weak
topology. Hence $\tau _{n}$ is the first hitting time of an open set by a
continuous process, so $\tau _{n}$ is a stopping time. Also, by the
assumption that $t\rightarrow P\left( t\right) $ is weakly continuous, it
follows from the uniform boundedness theorem that $\left\Vert P\left(
t\right) \right\Vert $ is bounded on $\left[ 0,T\right] .$ Hence for a.e. $%
\omega ,$ $\tau _{n}=\infty $ for all $n$ large enough. $\blacksquare $

It follows that it is possible to define the stochastic integral of \ref%
{z2june7h} as a local martingale when $t\rightarrow P\left( t\right) $ is
weakly continuous.\ In the derivation which follows, the computations will
pertain to such a weakly continuous process.

It remains to consider the convergence of a suitable subsequence of
\begin{equation*}
\int_{0}^{t}\left( Z\circ J^{-1}\right) ^{\ast }P_{n}^{l}\circ JdW
\end{equation*}%
to the integral of \ref{z2june7h}. The desired result follows. The proof is
similar to that given in \cite{pre07} for a similar situation in the context
of a Gelfand triple.

\begin{lemma}
\label{z24janl9gaa}In the above context, let $P\left( s\right)
-P_{k}^{l}\left( s\right) \equiv \Delta _{k}\left( s\right) .$ Let
\begin{equation*}
Z\in L^{2}\left( \left[ a,b\right] \times \Omega ,\mathcal{L}_{2}\left(
JQ^{1/2}U,W\right) \right)
\end{equation*}%
and let $P\in L^{2}\left( \left[ 0,T\right] \times \Omega ,W^{\prime
}\right) $ with both $P$ and $Z$ progressively measurable. Also suppose $%
t\rightarrow P\left( t\right) $ is weakly continuous. Then the integral
\begin{equation*}
\int_{0}^{t}\left( Z\circ J^{-1}\right) ^{\ast }P\circ JdW
\end{equation*}%
exists as a local martingale and the following limit is valid for a suitable
subsequence, still denoted by $k$%
\begin{equation*}
\lim_{k\rightarrow \infty }P\left( \left[ \sup_{t\in \left[ 0,T\right]
}\left\vert \int_{0}^{t}\left( Z\circ J^{-1}\right) ^{\ast }\Delta _{k}\circ
JdW\right\vert \geq \varepsilon \right] \right) =0.
\end{equation*}%
That is,
\begin{equation*}
\sup_{t\in \left[ 0,T\right] }\left\vert \int_{0}^{t}\left( Z\circ
J^{-1}\right) ^{\ast }\Delta _{k}\circ JdW\right\vert
\end{equation*}%
converges to 0 in probability.
\end{lemma}

\textbf{Proof:\ \ }The existence of the integral was dealt with earlier. Let
$k$ denote a subsequence for which for a.e. $\omega ,$
\begin{equation*}
P_{k}^{l}\left( \cdot ,\omega \right) \rightarrow P\left( \cdot ,\omega
\right)
\end{equation*}%
in $L^{p}\left( \left[ 0,T\right] ,W^{\prime }\right) $ and also $%
P_{k}^{l}\left( t,\omega \right) \rightarrow P\left( t,\omega \right) $ for
a.e. $t$. This is done as follows.
\begin{equation*}
P\left( \left\Vert P_{k}^{l}-P\right\Vert _{L^{p}\left( 0,T,W^{\prime
}\right) }>\lambda \right) \leq \frac{1}{\lambda }\int_{\Omega }\left\Vert
P_{k}^{l}-P\right\Vert _{L^{p}\left( 0,T,W^{\prime }\right) }dP
\end{equation*}%
and the integral on the right is small provided $k$ is large. Therefore,
there exists a subsequence still called $k$ such that
\begin{equation*}
P\left( \left\Vert P_{k}^{l}-P\right\Vert _{L^{p}\left( 0,T,W^{\prime
}\right) }>2^{-k}\right) <2^{-k}
\end{equation*}%
Then this satisfies the desired conditions.

From the assumption of weak continuity, there exists for a.e. $\omega $ a
constant, $C\left( \omega \right) $ such that
\begin{equation*}
\sup_{t\in \left[ 0,T\right] }\left\Vert P\left( t\right) \right\Vert \leq
C\left( \omega \right) .
\end{equation*}%
For the first part of the argument, assume that $C\left( \omega \right) $ is
independent of $\omega $ off a set of measure zero. Let $\left\{
e_{i}\right\} $ be an orthonormal basis of vectors in $W$. Thus $R\left(
e_{i}\right) $ is an orthonormal basis of vectors in $W^{\prime }$ where $R$
is the Riesz map. Hence
\begin{equation*}
P=\sum_{i=1}^{\infty }\left( P,R\left( e_{i}\right) \right) _{W^{\prime
}}R\left( e_{i}\right) =\sum_{i=1}^{\infty }\left( R^{-1}P,e_{i}\right)
_{W}R\left( e_{i}\right) =\sum_{i=1}^{\infty }\left\langle
P,e_{i}\right\rangle R\left( e_{i}\right)
\end{equation*}%
It follows that
\begin{equation*}
Px=\sum_{i=1}^{\infty }\left\langle P,e_{i}\right\rangle \left\langle
R\left( e_{i}\right) ,x\right\rangle
\end{equation*}

Let
\begin{equation*}
\pi _{n}P\equiv \sum_{i=1}^{n}\left\langle P,e_{i}\right\rangle R\left(
e_{i}\right)
\end{equation*}%
\begin{equation*}
P\left( \left[ \sup_{t\in \left[ 0,T\right] }\left\vert \int_{0}^{t}\left(
Z\circ J^{-1}\right) ^{\ast }\Delta _{k}\circ JdW\right\vert \geq
\varepsilon \right] \right)
\end{equation*}%
\begin{equation*}
\leq P\left( \left[ \sup_{t\in \left[ 0,T\right] }\left\vert
\int_{0}^{t}\left( Z\circ J^{-1}\right) ^{\ast }\pi _{n}\Delta _{k}\circ
JdW\right\vert \geq \varepsilon /3\right] \right) +
\end{equation*}%
\begin{equation}
P\left( \left[ \sup_{t\in \left[ 0,T\right] }\left\vert \int_{0}^{t}\left(
Z\circ J^{-1}\right) ^{\ast }\left( I-\pi _{n}\right) P\circ JdW\right\vert
\geq \varepsilon /3\right] \right) +  \label{5june2h}
\end{equation}%
\begin{equation}
P\left( \left[ \sup_{t\in \left[ 0,T\right] }\left\vert \int_{0}^{t}\left(
Z\circ J^{-1}\right) ^{\ast }\left( I-\pi _{n}\right) P_{k}^{l}\circ
JdW\right\vert \geq \varepsilon /3\right] \right)  \label{5june3h}
\end{equation}%
Using the Burkholder Davis Gundy inequality on \ref{5june2h} along with the
description of the quadratic variation given above,
\begin{equation*}
P\left( \left[ \sup_{t\in \left[ 0,T\right] }\left\vert \int_{0}^{t}\left(
Z\circ J^{-1}\right) ^{\ast }\left( I-\pi _{n}\right) P\circ JdW\right\vert
\geq \varepsilon /3\right] \right)
\end{equation*}%
\begin{equation*}
\leq \frac{3}{\varepsilon }\int_{\Omega }\sup_{t\in \left[ 0,T\right]
}\left\vert \int_{0}^{t}\left( Z\circ J^{-1}\right) ^{\ast }\left( I-\pi
_{n}\right) P\circ JdW\right\vert dP
\end{equation*}%
\begin{equation*}
\leq \frac{3C}{\varepsilon }\int_{\Omega }\left( \int_{0}^{T}\left\Vert
Z\right\Vert ^{2}\left\Vert \left( I-\pi _{n}\right) P\right\Vert
^{2}dt\right) ^{1/2}dP
\end{equation*}

\begin{equation*}
\leq \frac{3C}{\varepsilon }\left( \int_{\Omega }\int_{0}^{T}\left\Vert
Z\right\Vert ^{2}\left\Vert \left( I-\pi _{n}\right) P\right\Vert
^{2}dtdP\right) ^{1/2}
\end{equation*}%
This integral converges to 0 as $n\rightarrow \infty $ by the assumption
that $P$ is bounded along with the dominated convergence theorem applied to
the finite measure $\left\Vert Z\right\Vert ^{2}dtdP$. Letting $\eta >0$ be
given, choose $n$ large enough that the above term is less than $\eta .$
From now on, use this $n$. Thus \ref{5june2h} $\leq \eta .$

Next consider \ref{5june3h}. By the Burkholder Davis Gundy inequality again,%
\begin{equation*}
P\left( \left[ \sup_{t\in \left[ 0,T\right] }\left\vert \int_{0}^{t}\left(
Z\circ J^{-1}\right) ^{\ast }\left( I-\pi _{n}\right) P_{k}^{l}\circ
JdW\right\vert \geq \varepsilon /3\right] \right)
\end{equation*}%
\begin{equation*}
\leq \frac{3}{\varepsilon }\int_{\Omega }\sup_{t\in \left[ 0,T\right]
}\left\vert \int_{0}^{t}\left( Z\circ J^{-1}\right) ^{\ast }\left( I-\pi
_{n}\right) P_{k}^{l}\circ JdW\right\vert dP
\end{equation*}%
\begin{equation*}
\leq \frac{3C}{\varepsilon }\int_{\Omega }\left( \int_{0}^{T}\left\Vert
Z\right\Vert ^{2}\left\Vert \left( I-\pi _{n}\right) P_{k}^{l}\right\Vert
^{2}dt\right) ^{1/2}dP
\end{equation*}%
\begin{equation*}
\leq \frac{3C}{\varepsilon }\left( \int_{\Omega }\int_{0}^{T}\left\Vert
Z\right\Vert ^{2}\left\Vert \left( I-\pi _{n}\right) P_{k}^{l}\right\Vert
^{2}dtdP\right) ^{1/2}
\end{equation*}%
Next,
\begin{equation*}
\int_{\Omega }\int_{0}^{T}\left\Vert Z\right\Vert ^{2}\left\Vert \left(
I-\pi _{n}\right) P_{k}^{l}\right\Vert ^{2}dtdP
\end{equation*}%
\begin{equation}
\leq \int_{\Omega }\int_{0}^{T}\left\Vert Z\right\Vert ^{2}\left\Vert
P_{k}^{l}-P\right\Vert ^{2}dtdP+\int_{\Omega }\int_{0}^{T}\left\Vert
Z\right\Vert ^{2}\left\Vert \left( I-\pi _{n}\right) P\right\Vert ^{2}dtdP
\label{1may1f}
\end{equation}%
Now
\begin{equation*}
\frac{3C}{\varepsilon }\left( \int_{\Omega }\int_{0}^{T}\left\Vert
Z\right\Vert ^{2}\left\Vert \left( I-\pi _{n}\right) P\right\Vert
^{2}dtdP\right) ^{1/2}\leq \eta
\end{equation*}%
and so \ref{1may1f} is dominated by
\begin{equation*}
\leq \left( \frac{\eta \varepsilon }{3C}\right) ^{2}+\int_{\Omega
}\int_{0}^{T}\left\Vert Z\right\Vert ^{2}\left\Vert P_{k}^{l}-P\right\Vert
^{2}dtdP
\end{equation*}%
Therefore, \ref{5june3h} can be estimated as follows.
\begin{equation*}
P\left( \left[ \sup_{t\in \left[ 0,T\right] }\left\vert \int_{0}^{t}\left(
Z\circ J^{-1}\right) ^{\ast }\left( I-\pi _{n}\right) P_{k}^{l}\circ
JdW\right\vert \geq \varepsilon /3\right] \right) \leq
\end{equation*}%
\begin{equation*}
\leq \frac{3C}{\varepsilon }\left( \int_{\Omega }\int_{0}^{T}\left\Vert
Z\right\Vert ^{2}\left\Vert \left( I-\pi _{n}\right) P_{k}^{l}\right\Vert
^{2}dtdP\right) ^{1/2}
\end{equation*}%
\begin{eqnarray*}
&\leq &\frac{3C}{\varepsilon }\left( \left( \frac{\eta \varepsilon }{3C}%
\right) ^{2}+\int_{\Omega }\int_{0}^{T}\left\Vert Z\right\Vert
^{2}\left\Vert P_{k}^{l}-P\right\Vert ^{2}dtdP\right) ^{1/2} \\
&\leq &\eta +\frac{3C}{\varepsilon }\left( \int_{\Omega
}\int_{0}^{T}\left\Vert Z\right\Vert ^{2}\left\Vert P_{k}^{l}-P\right\Vert
^{2}dtdP\right) ^{1/2}
\end{eqnarray*}

\vspace{1pt}It follows that
\begin{equation*}
P\left( \left[ \sup_{t\in \left[ 0,T\right] }\left\vert \int_{0}^{t}\left(
Z\circ J^{-1}\right) ^{\ast }\Delta _{k}\circ JdW\right\vert \geq
\varepsilon \right] \right)
\end{equation*}%
\begin{eqnarray*}
&\leq &2\eta +\frac{3C}{\varepsilon }\left( \int_{\Omega
}\int_{0}^{T}\left\Vert Z\right\Vert ^{2}\left\Vert P_{k}^{l}-P\right\Vert
^{2}dtdP\right) ^{1/2}+ \\
&&+P\left( \left[ \sup_{t\in \left[ 0,T\right] }\left\vert
\int_{0}^{t}\left( Z\circ J^{-1}\right) ^{\ast }\pi _{n}\Delta _{k}\circ
JdW\right\vert \geq \varepsilon /3\right] \right)
\end{eqnarray*}

\vspace{1pt}By weak continuity,
\begin{equation*}
\lim_{k\rightarrow \infty }\pi _{n}\left( P\left( s\right) -P_{k}^{l}\left(
s\right) \right) =\lim_{k\rightarrow \infty }\sum_{i=1}^{n}\left\langle
P\left( s\right) -P_{k}^{l}\left( s\right) ,e_{i}\right\rangle R\left(
e_{i}\right) =0\ \ \text{ }a.e\ \omega
\end{equation*}%
It follows that
\begin{equation*}
\lim_{k\rightarrow \infty }\int_{\Omega }\int_{0}^{T}\left\Vert \pi
_{n}\left( P\left( s\right) -P_{k}^{l}\left( s\right) \right) \right\Vert
_{W^{\prime }}^{2}\left\Vert Z\left( s\right) \right\Vert _{\mathcal{L}%
_{2}}^{2}dsdP=0
\end{equation*}%
because you can apply the dominated convergence theorem with respect to the
measure $\left\Vert Z\left( s\right) \right\Vert ^{2}dsdP$ along with the
assumption that $\left\Vert P\left( t\right) \right\Vert $ is bounded
independent of $\omega $.

Therefore, it follows that
\begin{equation*}
\lim_{k\rightarrow \infty }P\left( \left[ \sup_{t\in \left[ 0,T\right]
}\left\vert \int_{0}^{t}\left( Z\circ J^{-1}\right) ^{\ast }\pi _{n}\Delta
_{k}\circ JdW\right\vert \geq \varepsilon /3\right] \right) =0.
\end{equation*}%
Here is why. By the Burkholder Davis Gundy theorem, and the description of
the quadratic variation of a stochastic integral,%
\begin{equation*}
\int_{\Omega }\sup_{t\in \left[ 0,T\right] }\left\vert \int_{0}^{t}\left(
Z\circ J^{-1}\right) ^{\ast }\pi _{n}\Delta _{k}\circ JdW\right\vert dP
\end{equation*}%
\begin{eqnarray*}
&\leq &C\int_{\Omega }\left( \int_{0}^{T}\left\Vert \pi _{n}\left( P\left(
s\right) -P_{k}^{l}\left( s\right) \right) \right\Vert _{W^{\prime
}}^{2}\left\Vert Z\left( s\right) \right\Vert _{\mathcal{L}%
_{2}}^{2}ds\right) ^{1/2}dP \\
&\leq &C\left( \int_{\Omega }\int_{0}^{T}\left\Vert \pi _{n}\left( P\left(
s\right) -P_{k}^{l}\left( s\right) \right) \right\Vert _{W^{\prime
}}^{2}\left\Vert Z\left( s\right) \right\Vert _{\mathcal{L}%
_{2}}^{2}dsdP\right) ^{1/2}
\end{eqnarray*}%
which converges to 0 as $k\rightarrow \infty $.

\ If $k$ is large enough, this implies%
\begin{equation*}
P\left( \left[ \sup_{t\in \left[ 0,T\right] }\left\vert \int_{0}^{t}\left(
Z\circ J^{-1}\right) ^{\ast }\Delta _{k}\circ JdW\right\vert \geq
\varepsilon \right] \right)
\end{equation*}%
\begin{equation*}
\leq 3\eta +\frac{3C}{\varepsilon }\left( \int_{\Omega
}\int_{0}^{T}\left\Vert Z\right\Vert ^{2}\left\Vert \left(
P_{k}^{l}-P\right) \right\Vert ^{2}dtdP\right) ^{1/2}
\end{equation*}%
The last integral also converges to 0 because of the assumption that $P$,
hence $P_{k}^{l}$ is bounded, and the dominated convergence theorem. Thus
\begin{equation*}
\lim_{k\rightarrow \infty }\sup_{t\in \left[ 0,T\right] }\left\vert
\int_{0}^{t}\left( Z\circ J^{-1}\right) ^{\ast }\Delta _{k}\circ
JdW\right\vert =0\text{ in probability.}
\end{equation*}

Now to finish the argument, define the stopping time
\begin{equation*}
\tau _{m}\equiv \inf \left\{ t>0:\left\Vert P\left( t\right) \right\Vert
>m\right\} .
\end{equation*}%
As discussed earlier, this is a valid stopping time by the weak continuity
of $t\rightarrow P\left( t\right) $. Then
\begin{equation*}
X^{\tau _{m}}-\left( X_{k}^{l}\right) ^{\tau _{m}}
\end{equation*}%
still converges pointwise to 0 as $k\rightarrow \infty $. Let $\Delta
_{k}^{\tau _{m}}=\left( P\left( s\right) -P_{k}^{l}\left( s\right) \right)
^{\tau _{m}}$

Now consider
\begin{equation*}
A_{k\varepsilon }\equiv \left[ \sup_{t\in \left[ 0,T\right] }\left\vert
\int_{0}^{t}\left( Z\circ J^{-1}\right) ^{\ast }\Delta _{k}\circ
JdW\right\vert \geq \varepsilon \right]
\end{equation*}%
Then
\begin{equation*}
P\left( A_{k\varepsilon }\cap \left[ \tau _{m}=\infty \right] \right) \leq
P\left( \left[ \sup_{t\in \left[ 0,T\right] }\left\vert \int_{0}^{t}\left(
Z\circ J^{-1}\right) ^{\ast }\Delta _{k}^{\tau _{m}}\circ JdW\right\vert
\geq \varepsilon \right] \right)
\end{equation*}%
which converges to 0 as $k\rightarrow \infty $ by the first part of the
argument. This is because $\left\Vert \left( P^{\tau _{m}}\right)
_{k}^{l}\right\Vert $ and $\left\Vert P^{\tau _{m}}\right\Vert $ are both
bounded by $m$ uniformly off a set of measure zero. Now $A_{k\varepsilon }$
can be partitioned in the following way.
\begin{equation*}
A_{k\varepsilon }=\cup _{m=1}^{\infty }A_{k\varepsilon }\cap \left( \left[
\tau _{m}=\infty \right] \setminus \left[ \tau _{m-1}<\infty \right] \right)
\end{equation*}%
Thus
\begin{equation*}
P\left( A_{k\varepsilon }\right) =\sum_{m=1}^{\infty }P\left(
A_{k\varepsilon }\cap \left( \left[ \tau _{m}=\infty \right] \setminus \left[
\tau _{m-1}<\infty \right] \right) \right)
\end{equation*}%
and $P\left( A_{k\varepsilon }\cap \left( \left[ \tau _{m}=\infty \right]
\setminus \left[ \tau _{m-1}<\infty \right] \right) \right) \leq P\left(
\left( \left[ \tau _{m}=\infty \right] \setminus \left[ \tau _{m-1}<\infty %
\right] \right) \right) $ which is summable, since the sets are disjoint.
Hence one can apply the dominated convergence theorem and conclude that
\begin{equation*}
\lim_{k\rightarrow \infty }P\left( A_{k\varepsilon }\right)
=\sum_{m=1}^{\infty }\lim_{k\rightarrow \infty }P\left( A_{k\varepsilon
}\cap \left( \left[ \tau _{m}=\infty \right] \setminus \left[ \tau
_{m-1}<\infty \right] \right) \right) =0.\ \blacksquare
\end{equation*}

The notation used in the above is inelegant. In fact it is more attractive
to write this in the form
\begin{equation*}
\int_{0}^{t}\left( Z\circ J^{-1}\right) ^{\ast }P\circ JdW=\int_{0}^{t}PdM
\end{equation*}%
where $M\left( t\right) =\int_{0}^{t}ZdW.$ However, this is of no importance
in what follows, and written in the above inelegant form, assertions about
the quadratic variation are possibly more obvious.

\section{The Integral Equation\label{integralequation}}

For a set of measure zero $N_{0}$ and $\omega \notin N_{0},$
\begin{equation}
P\left( t\right) \equiv BX_{0}+\int_{0}^{t}Y\left( s\right)
ds+B\int_{0}^{t}Z\left( s\right) dW\left( s\right) ,  \label{z24jane1g}
\end{equation}%
for all $t$. Let $X$ be progressively measurable into $V$ such that for a
set of measure zero $S\subseteq \left[ 0,T\right] $ and $t\notin S,$%
\begin{equation*}
BX\left( t\right) =P\left( t\right) \text{ for }a.e.\omega .
\end{equation*}%
The exceptional set in the above may depend on $t$. In short, for all $%
t\notin S,$
\begin{equation}
BX\left( t\right) =BX_{0}+\int_{0}^{t}Y\left( s\right)
ds+B\int_{0}^{t}Z\left( s\right) dW\left( s\right) \text{ a.e. }\omega ,
\label{25oct1h}
\end{equation}%
the exceptional set possibly depending on $t$.

Also assume that $X_{0}\in L^{2}\left( \Omega ;W\right) $ and is $\mathcal{F}%
_{0}$ measurable, where $Z$ is $\mathcal{L}_{2}\left( Q^{1/2}U,W\right) $
progressively measurable and
\begin{equation*}
\left\Vert Z\right\Vert _{L^{2}\left( \left[ 0,T\right] \times \Omega ,%
\mathcal{L}_{2}\left( Q^{1/2}U,W\right) \right) }<\infty .
\end{equation*}%
This is what is needed to define the stochastic integral in the above
formula.

$X,Y$ satisfy
\begin{equation*}
X\in K\equiv L^{p}\left( \left[ 0,T\right] \times \Omega ;V\right) \cap
L^{2}\left( \left[ 0,T\right] \times \Omega ,W\right) ,Y\in K^{\prime
}=L^{p^{\prime }}\left( \left[ 0,T\right] \times \Omega ;V^{\prime }\right)
\end{equation*}%
where $1/p^{\prime }+1/p=1,\ p>1,$ and $X,Y$ are progressively measurable
into $V$ and $V^{\prime }$ respectively.

Also, by enlarging $N_{0}$ if necessary, one can assume that off $N_{0},$
the stochastic integral in \ref{z24jane1g} is continuous into $W^{\prime }$
and $Y\left( \cdot ,\omega \right) \in L^{p^{\prime }}\left( 0,T,V^{\prime
}\right) $ so that the deterministic integral in this equation is also
continuous as a function with values in $V^{\prime }$, also that $%
t\rightarrow X\left( t,\omega \right) \in L^{p}\left( 0,T,V\right) $ for $%
\omega \notin N_{0}$. From now on, let $N_{0}$ be so enlarged.

The goal is to prove the following Ito formula for $P\left( t\right) $
defined as the right side of \ref{z24jane1g}.
\begin{equation*}
\left\langle P\left( t\right) ,X\left( t\right) \right\rangle =\left\langle
BX_{0},X_{0}\right\rangle +\int_{0}^{t}\left( 2\left\langle Y\left( s\right)
,X\left( s\right) \right\rangle +\left\langle BZ,Z\right\rangle _{\mathcal{L}%
_{2}}\right) ds
\end{equation*}%
\begin{equation}
+\int_{0}^{t}\left( Z\circ J^{-1}\right) ^{\ast }P\circ JdW
\label{z27jane2g}
\end{equation}%
The most significant feature of the last term is that it is a local
martingale. The term $\left\langle BZ,Z\right\rangle _{\mathcal{L}_{2}}$
will be discussed later. In the stochastic integral, $\left( Z\circ
J^{-1}\right) ^{\ast }P\circ J$ has values in $\mathcal{L}_{2}\left(
Q^{1/2}U,\mathbb{R}\right) $ and so it makes sense to consider this
stochastic integral.

The argument for the Ito formula will be based on a formula which follows in
the next lemma.

\begin{lemma}
\label{z24janl3g}In the situation of the above integral equation, the
following formula holds for a.e. $\omega $ for $0<$ $s<t,$ where $s,t\notin
S $ where $M\left( t\right) \equiv \int_{0}^{t}Z\left( u\right) dW\left(
u\right) $ which has values in $W$. In the following, $\left\langle \cdot
,\cdot \right\rangle $ denotes the duality pairing between $V,V^{\prime }.$
\begin{equation*}
\left\langle BX\left( t\right) ,X\left( t\right) \right\rangle =\left\langle
BX\left( s\right) ,X\left( s\right) \right\rangle +
\end{equation*}%
\begin{equation*}
+2\int_{s}^{t}\left\langle Y\left( u\right) ,X\left( t\right) \right\rangle
du+\left\langle B\left( M\left( t\right) -M\left( s\right) \right) ,M\left(
t\right) -M\left( s\right) \right\rangle
\end{equation*}%
\begin{equation*}
-\left\langle B\left( X\left( t\right) -X\left( s\right) -\left( M\left(
t\right) -M\left( s\right) \right) \right) ,X\left( t\right) -X\left(
s\right) -\left( M\left( t\right) -M\left( s\right) \right) \right\rangle
\end{equation*}%
\begin{equation}
+2\left\langle BX\left( s\right) ,M\left( t\right) -M\left( s\right)
\right\rangle  \label{z24jane2g}
\end{equation}%
Also for $t>0$
\begin{equation*}
\left\langle BX\left( t\right) ,X\left( t\right) \right\rangle =\left\langle
BX_{0},X_{0}\right\rangle +2\int_{0}^{t}\left\langle Y\left( u\right)
,X\left( t\right) \right\rangle du+2\left\langle BX_{0},M\left( t\right)
\right\rangle +
\end{equation*}%
\begin{equation}
\left\langle BM\left( t\right) ,M\left( t\right) \right\rangle -\left\langle
B\left( X\left( t\right) -X_{0}-M\left( t\right) \right) ,X\left( t\right)
-X_{0}-M\left( t\right) \right\rangle  \label{z26jane1g}
\end{equation}
\end{lemma}

\textbf{Proof:} From the formula which is assumed to hold, \textbf{\ }
\begin{eqnarray*}
BX\left( t\right) &=&BX_{0}+\int_{0}^{t}Y\left( u\right) du+BM\left( t\right)
\\
BX\left( s\right) &=&BX_{0}+\int_{0}^{s}Y\left( u\right) du+BM\left( s\right)
\end{eqnarray*}%
Then
\begin{equation*}
BM\left( t\right) -BM\left( s\right) +\int_{s}^{t}Y\left( u\right)
du=BX\left( t\right) -BX\left( s\right)
\end{equation*}%
It follows that%
\begin{equation*}
\left\langle B\left( M\left( t\right) -M\left( s\right) \right) ,M\left(
t\right) -M\left( s\right) \right\rangle -
\end{equation*}%
\begin{equation*}
\left\langle B\left( X\left( t\right) -X\left( s\right) -\left( M\left(
t\right) -M\left( s\right) \right) \right) ,X\left( t\right) -X\left(
s\right) -\left( M\left( t\right) -M\left( s\right) \right) \right\rangle
\end{equation*}%
\begin{equation*}
+2\left\langle BX\left( s\right) ,M\left( t\right) -M\left( s\right)
\right\rangle
\end{equation*}%
\begin{eqnarray*}
&=&\left\langle B\left( M\left( t\right) -M\left( s\right) \right) ,M\left(
t\right) -M\left( s\right) \right\rangle -\left\langle B\left( X\left(
t\right) -X\left( s\right) \right) ,X\left( t\right) -X\left( s\right)
\right\rangle \\
&&+2\left\langle B\left( X\left( t\right) -X\left( s\right) \right) ,M\left(
t\right) -M\left( s\right) \right\rangle
\end{eqnarray*}%
\begin{equation*}
-\left\langle B\left( M\left( t\right) -M\left( s\right) \right) ,M\left(
t\right) -M\left( s\right) \right\rangle +2\left\langle BX\left( s\right)
,M\left( t\right) -M\left( s\right) \right\rangle
\end{equation*}%
Some terms cancel and this yields
\begin{equation*}
=-\left\langle B\left( X\left( t\right) -X\left( s\right) \right) ,X\left(
t\right) -X\left( s\right) \right\rangle +2\left\langle BX\left( t\right)
,M\left( t\right) -M\left( s\right) \right\rangle
\end{equation*}%
\begin{equation*}
=-\left\langle B\left( X\left( t\right) -X\left( s\right) \right) ,X\left(
t\right) -X\left( s\right) \right\rangle +2\left\langle B\left( M\left(
t\right) -M\left( s\right) \right) ,X\left( t\right) \right\rangle
\end{equation*}%
\begin{eqnarray*}
&=&-\left\langle B\left( X\left( t\right) -X\left( s\right) \right) ,X\left(
t\right) -X\left( s\right) \right\rangle \\
&&+2\left\langle B\left( X\left( t\right) -X\left( s\right) \right)
-\int_{s}^{t}Y\left( u\right) du,X\left( t\right) \right\rangle
\end{eqnarray*}%
\begin{eqnarray*}
&=&-\left\langle BX\left( t\right) ,X\left( t\right) \right\rangle
-\left\langle BX\left( s\right) ,X\left( s\right) \right\rangle \\
&&+2\left\langle BX\left( t\right) ,X\left( s\right) \right\rangle
+2\left\langle BX\left( t\right) ,X\left( t\right) \right\rangle \\
&&-2\left\langle BX\left( s\right) ,X\left( t\right) \right\rangle
-2\int_{s}^{t}\left\langle Y\left( u\right) ,X\left( t\right) \right\rangle
du
\end{eqnarray*}%
\begin{equation*}
=\left\langle BX\left( t\right) ,X\left( t\right) \right\rangle
-\left\langle BX\left( s\right) ,X\left( s\right) \right\rangle
-2\int_{s}^{t}\left\langle Y\left( u\right) ,X\left( t\right) \right\rangle
du
\end{equation*}%
Therefore,
\begin{eqnarray*}
&&\left\langle BX\left( t\right) ,X\left( t\right) \right\rangle
-\left\langle BX\left( s\right) ,X\left( s\right) \right\rangle \\
&=&2\int_{s}^{t}\left\langle Y\left( u\right) ,X\left( t\right)
\right\rangle du+\left\langle B\left( M\left( t\right) -M\left( s\right)
\right) ,M\left( t\right) -M\left( s\right) \right\rangle
\end{eqnarray*}%
\begin{equation*}
-\left\langle B\left( X\left( t\right) -X\left( s\right) -\left( M\left(
t\right) -M\left( s\right) \right) \right) ,X\left( t\right) -X\left(
s\right) -\left( M\left( t\right) -M\left( s\right) \right) \right\rangle
\end{equation*}%
\begin{equation*}
+2\left\langle BX\left( s\right) ,M\left( t\right) -M\left( s\right)
\right\rangle \
\end{equation*}%
The case with $X_{0}$ is similar. $\blacksquare $

The following lemma is what will be used. It says that you can replace $%
BX\left( t\right) $ in Lemma \ref{z24janl3g} with $P\left( t\right) $ for
all $\omega $ off a set of measure zero. This just involves substituting $%
P\left( t\right) $ for $BX\left( t\right) $ in the above formula since these
are equal.

\begin{lemma}
\label{20julyl1}For given $t,s,s<t,\ s,t\notin S,$ the following holds for
a.e. $\omega $%
\begin{equation*}
\left\langle P\left( t\right) ,X\left( t\right) \right\rangle =\left\langle
P\left( s\right) ,X\left( s\right) \right\rangle +
\end{equation*}%
\begin{equation*}
+2\int_{s}^{t}\left\langle Y\left( u\right) ,X\left( t\right) \right\rangle
du+\left\langle B\left( M\left( t\right) -M\left( s\right) \right) ,M\left(
t\right) -M\left( s\right) \right\rangle
\end{equation*}%
\begin{equation*}
-\left\langle \left( P\left( t\right) -P\left( s\right) -B\left( M\left(
t\right) -M\left( s\right) \right) \right) ,X\left( t\right) -X\left(
s\right) -\left( M\left( t\right) -M\left( s\right) \right) \right\rangle
\end{equation*}%
\begin{equation}
+2\left\langle P\left( s\right) ,M\left( t\right) -M\left( s\right)
\right\rangle  \label{19july1}
\end{equation}%
and for $y\in W,$%
\begin{equation}
\left\vert \left\langle P\left( t\right) ,y\right\rangle \right\vert \leq
\left\langle P\left( t\right) ,X\left( t\right) \right\rangle
^{1/2}\left\langle By,y\right\rangle ^{1/2}  \label{19july3}
\end{equation}
\end{lemma}

Let $\left\{ \mathcal{P}_{k}\right\} $ denote a sequence of nested
partitions of $\left[ 0,T\right] $ which satisfy the conditions needed in
Lemma \ref{z23decl1g} and also for $X_{k}^{l},X_{k}^{r}$ described there,
\begin{equation*}
X_{k}^{l},X_{k}^{r}\rightarrow X\text{ in }K\equiv L^{p}\left( \left[ 0,T%
\right] \times \Omega ;V\right)
\end{equation*}%
Each $\mathcal{P}_{k}$ contains no points of $S$. In what follows $N$ will
be a set of measure zero which includes $N_{0}$. Each $t\in \mathcal{P}_{k}$
has the property that $BX\left( t\right) =P\left( t\right) $ a.e., that is
for $\omega \in N_{t}$ a set of measure zero. Let $N$ also include the set
of measure zero
\begin{equation*}
\cup _{k}\cup \left\{ N_{t}:t\in \mathcal{P}_{k}\right\}
\end{equation*}%
Hence for each $t\in \mathcal{P}_{k},$ $BX\left( t\right) =P\left( t\right) $
for all $\omega \notin N$. Let
\begin{equation*}
\mathcal{D}\equiv \cup _{k}\mathcal{P}_{k}.
\end{equation*}%
Thus $\mathcal{D}$ is dense in $\left[ 0,T\right] $. Since $BX\left(
t\right) =P\left( t\right) $ a.e. $\omega $ for $t\notin S$, it follows that
as $k\rightarrow \infty ,$
\begin{equation*}
P_{k}^{l}\equiv BX_{k}^{l},\ P_{k}^{r}\equiv BX_{k}^{r}\text{ both converge
to }P\text{ in }L^{2}\left( \left[ 0,T\right] \times \Omega ,W^{\prime
}\right)
\end{equation*}%
For convenience, we only consider those points of $\mathcal{P}_{k}$ which
are less than $T$. These are the ones which are used in the left step
functions.

\section{The Main Estimate\label{mainestimate}}

The following estimate holds and it is this estimate which is the main idea
in proving the Ito formula. The last assertion about continuity is like the
well known result that if $y\in L^{p}\left( 0,T;V\right) $ and $y^{\prime
}\in L^{p^{\prime }}\left( 0,T;V^{\prime }\right) ,$ then $y$ is actually
continuous with values in $H$, for $V,H,V^{\prime }$ a Gelfand triple.

In all that follows $\left\{ e_{i}\right\} $ will be the vectors of Lemma %
\ref{7mayl1h}.
\begin{equation*}
\sum_{i=1}^{\infty }\left\langle BX\left( t\right) ,e_{i}\right\rangle
^{2}=\left\langle BX\left( t\right) ,X\left( t\right) \right\rangle
\end{equation*}

\begin{lemma}
\label{z24janl1g}\label{27julyc1} In the situation of Section \ref%
{integralequation}
\begin{eqnarray}
&&E\left( \sup_{t\in \mathcal{D}}\left\langle P\left( t\right) ,X\left(
t\right) \right\rangle \right)  \notag \\
&<&C\left( \left\vert \left\vert Y\right\vert \right\vert _{K^{\prime
}},\left\vert \left\vert X\right\vert \right\vert _{K},\left\vert \left\vert
Z\right\vert \right\vert _{J},\left\Vert \left\langle
BX_{0},X_{0}\right\rangle \right\Vert _{L^{1}\left( \Omega \right) }\right)
<\infty .  \label{9aug2h}
\end{eqnarray}%
where $\mathcal{D}$ is a dense subset of $\left[ 0,T\right] ,$ such that for
$t\in \mathcal{D}$, $BX\left( t\right) =P\left( t\right) $ a.e. $\omega $
where
\begin{eqnarray*}
J &=&L^{2}\left( \left[ 0,T\right] \times \Omega ;\mathcal{L}_{2}\left(
Q^{1/2}U;W\right) \right) ,K\equiv L^{p}\left( \left[ 0,T\right] \times
\Omega ;V\right) , \\
K^{\prime } &\equiv &L^{p^{\prime }}\left( \left[ 0,T\right] \times \Omega
;V^{\prime }\right) ,
\end{eqnarray*}%
the $\sigma $ algebra being the progressively measurable sets. $C$ is a
continuous function of its arguments and $C\left( 0,0,0,0\right) =0$. Also,
for a.e. $\omega ,$%
\begin{equation}
\sup_{t\in \left[ 0,T\right] }\sum_{k}\left\langle P\left( t\right)
,e_{k}\right\rangle ^{2}\leq C\left( \omega \right) <\infty \text{ }
\label{9aug1h}
\end{equation}%
For $t\in \mathcal{D}$, then for all $\omega \notin N,$ a set of measure
zero,%
\begin{equation}
\left\langle P\left( t\right) ,X\left( t\right) \right\rangle =\left\langle
BX\left( t\right) ,X\left( t\right) \right\rangle \leq C\left( \omega
\right) <\infty ,  \label{4aug1h}
\end{equation}%
\begin{equation*}
\int_{\Omega }C\left( \omega \right) dP<\infty .
\end{equation*}%
When each of $\left\vert \left\vert Y\right\vert \right\vert _{K^{\prime
}},\left\vert \left\vert X\right\vert \right\vert _{K},\left\vert \left\vert
Z\right\vert \right\vert _{J},\left\Vert \left\langle
BX_{0},X_{0}\right\rangle \right\Vert _{L^{1}\left( \Omega \right) }$ equal
zero, $C\left( \omega \right) $ can be taken to be 0 also. For a.e. $\omega
, $ $t\rightarrow P\left( t\right) $ is weakly continuous into $W^{\prime }.$
In addition to this, $P$ is progressively measurable into $W^{\prime }$.
\end{lemma}

\textbf{Proof:\ } For $t_{j}>0,X\left( t_{j}\right) $ is just the value of $%
X $ at $t_{j}$ but when $t=0,$ the definition of $X\left( 0\right) $ in this
step function is $X\left( 0\right) \equiv 0$. Consider the formula in Lemma %
\ref{20julyl1}. This is applied to $\mathcal{P}_{k}$ to obtain
\begin{equation*}
\left\langle P\left( t_{m}\right) ,X\left( t_{m}\right) \right\rangle
-\left\langle BX_{0},X_{0}\right\rangle
=2\sum_{j=1}^{m-1}\int_{t_{j}}^{t_{j+1}}\left\langle Y\left( u\right)
,X_{k}^{r}\left( u\right) \right\rangle du+
\end{equation*}%
\begin{eqnarray*}
&&+2\sum_{j=1}^{m-1}\left\langle P\left( t_{j}\right)
,\int_{t_{j}}^{t_{j+1}}Z\left( u\right) dW\right\rangle \\
&&+\sum_{j=1}^{m-1}\left\langle B\left( M\left( t_{j+1}\right) -M\left(
t_{j}\right) \right) ,M\left( t_{j+1}\right) -M\left( t_{j}\right)
\right\rangle
\end{eqnarray*}%
\begin{eqnarray*}
&&-\sum_{j=1}^{m-1}\left\langle P\left( t_{j+1}\right) -P\left( t_{j}\right)
-B\left( M\left( t_{j+1}\right) -M\left( t_{j}\right) \right) \right. , \\
&&\left. X\left( t_{j+1}\right) -X\left( t_{j}\right) -\left( M\left(
t_{j+1}\right) -M\left( t_{j}\right) \right) \right\rangle
\end{eqnarray*}%
\begin{equation*}
+2\int_{0}^{t_{1}}\left\langle Y\left( u\right) ,X\left( t_{1}\right)
\right\rangle du+2\left\langle BX_{0},\int_{0}^{t_{1}}Z\left( u\right)
dW\right\rangle +\left\langle BM\left( t_{1}\right) ,M\left( t_{1}\right)
\right\rangle
\end{equation*}%
\begin{equation}
-\left\langle P\left( t_{1}\right) -BX_{0}-BM\left( t_{1}\right) ,X\left(
t_{1}\right) -X_{0}-M\left( t_{1}\right) \right\rangle  \label{z26jane2g}
\end{equation}%
First consider the terms near to the end of the above expression,
\begin{equation*}
2\int_{0}^{t_{1}}\left\langle Y\left( u\right) ,X\left( t_{1}\right)
\right\rangle du+2\left\langle BX_{0},\int_{0}^{t_{1}}Z\left( u\right)
dW\right\rangle +\left\langle BM\left( t_{1}\right) ,M\left( t_{1}\right)
\right\rangle
\end{equation*}%
Each term of the above converges to 0 for a.e. $\omega $ as $k\rightarrow
\infty $ and in $L^{1}\left( \Omega \right) $ if a suitable subsequence is
used. This follows right away for the second two terms from the Ito isometry
and continuity properties of the stochastic integral. Consider the first
term. This term is dominated by
\begin{eqnarray*}
&&\left( \int_{0}^{t_{1}}\left\Vert Y\left( u\right) \right\Vert ^{p^{\prime
}}du\right) ^{1/p^{\prime }}\left( \int_{0}^{T}\left\Vert X_{k}^{r}\left(
u\right) \right\Vert ^{p}du\right) ^{1/p} \\
&\leq &C\left( \omega \right) \left( \int_{0}^{t_{1}}\left\Vert Y\left(
u\right) \right\Vert ^{p^{\prime }}du\right) ^{1/p^{\prime }},\ C\left(
\omega \right) <\infty
\end{eqnarray*}%
By assumption, and Holder's inequality, the top expression converges to 0 in
$L^{1}\left( \Omega \right) $. Hence there is a further subsequence for
which it converges pointwise.

At this time, not much is known about the last term in \ref{z26jane2g}, but
it is negative and is about to be neglected anyway. The reason it is
negative is that it equals
\begin{equation*}
-\left\langle B\left( X\left( t_{1}\right) -X_{0}-M\left( t_{1}\right)
\right) ,X\left( t_{1}\right) -X_{0}-M\left( t_{1}\right) \right\rangle
\end{equation*}

The term involving the stochastic integral equals
\begin{equation*}
2\sum_{j=1}^{m-1}\left\langle P\left( t_{j}\right)
,\int_{t_{j}}^{t_{j+1}}Z\left( u\right) dW\right\rangle
\end{equation*}%
By Theorem \ref{z19mayt1h}, this equals
\begin{equation*}
2\int_{t_{1}}^{t_{m}}\left( Z\circ J^{-1}\right) ^{\ast }P_{k}^{l}\circ JdW
\end{equation*}%
Also note that since $\left\langle BM\left( t_{1}\right) ,M\left(
t_{1}\right) \right\rangle $ converges to 0 in $L^{1}\left( \Omega \right) $
and for a.e. $\omega ,$ the sum involving%
\begin{equation*}
\left\langle B\left( M\left( t_{j+1}\right) -M\left( t_{j}\right) \right)
,M\left( t_{j+1}\right) -M\left( t_{j}\right) \right\rangle
\end{equation*}%
can be started at 0 rather than 1 at the expense of adding in a term which
converges to 0 a.e. and in $L^{1}\left( \Omega \right) $. Thus \ref%
{z26jane2g} is of the form
\begin{equation*}
\left\langle P\left( t_{m}\right) ,X\left( t_{m}\right) \right\rangle
-\left\langle BX_{0},X_{0}\right\rangle =e\left( k\right)
+2\int_{0}^{t_{m}}\left\langle Y\left( u\right) ,X_{k}^{r}\left( u\right)
\right\rangle du+
\end{equation*}%
\begin{eqnarray*}
&&+2\int_{0}^{t_{m}}\left( Z\circ J^{-1}\right) ^{\ast }P_{k}^{l}\circ JdW \\
&&+\sum_{j=0}^{m-1}\left\langle B\left( M\left( t_{j+1}\right) -M\left(
t_{j}\right) \right) ,M\left( t_{j+1}\right) -M\left( t_{j}\right)
\right\rangle
\end{eqnarray*}%
\begin{eqnarray*}
&&-\sum_{j=1}^{m-1}\left\langle B\left( X\left( t_{j+1}\right) -X\left(
t_{j}\right) -\left( M\left( t_{j+1}\right) -M\left( t_{j}\right) \right)
\right) \right. , \\
&&\left. X\left( t_{j+1}\right) -X\left( t_{j}\right) -\left( M\left(
t_{j+1}\right) -M\left( t_{j}\right) \right) \right\rangle
\end{eqnarray*}%
\begin{equation}
-\left\langle P\left( t_{1}\right) -BX_{0}-BM\left( t_{1}\right) ,X\left(
t_{1}\right) -X_{0}-M\left( t_{1}\right) \right\rangle  \label{z24jane12g}
\end{equation}%
where $e\left( k\right) \rightarrow 0$ for a.e. $\omega $ and also in $%
L^{1}\left( \Omega \right) $.

By definition, $M\left( t_{j+1}\right) -M\left( t_{j}\right)
=\int_{t_{j}}^{t_{j+1}}ZdW.$ Now it follows, on discarding the negative
terms,%
\begin{equation*}
\left\langle P\left( t_{m}\right) ,X\left( t_{m}\right) \right\rangle
-\left\langle BX_{0},X_{0}\right\rangle \leq e\left( k\right)
+2\int_{0}^{t_{m}}\left\langle Y\left( u\right) ,X_{k}^{r}\left( u\right)
\right\rangle du+
\end{equation*}%
\begin{equation*}
+2\int_{0}^{t_{m}}\left( Z\circ J^{-1}\right) ^{\ast }P_{k}^{l}\circ
JdW+\sum_{j=0}^{m-1}\left\langle
B\int_{t_{j}}^{t_{j+1}}ZdW,\int_{t_{j}}^{t_{j+1}}ZdW\right\rangle
\end{equation*}%
Thus also%
\begin{equation*}
\left\langle P\left( t_{m}\right) ,X\left( t_{m}\right) \right\rangle
-\left\langle BX_{0},X_{0}\right\rangle \leq e\left( k\right)
+2\int_{0}^{T}\left\vert \left\langle Y\left( u\right) ,X_{k}^{r}\left(
u\right) \right\rangle \right\vert du+
\end{equation*}%
\begin{equation*}
+2\int_{0}^{t_{m}}\left( Z\circ J^{-1}\right) ^{\ast }P_{k}^{l}\circ
JdW+\sum_{j=0}^{m-1}\left\langle
B\int_{t_{j}}^{t_{j+1}}ZdW,\int_{t_{j}}^{t_{j+1}}ZdW\right\rangle
\end{equation*}

The next task is to somehow take the expectation of both sides. The
difficulty in doing this is that the stochastic integral is only a local
martingale. Let
\begin{equation*}
\tau _{p}=\inf \left\{ t:\left\langle P_{k}^{l}\left( t\right)
,X_{k}^{l}\left( t\right) \right\rangle >p\right\}
\end{equation*}%
By right continuity of $P_{k}^{l}$ and $X_{k}^{l},$ this is a well defined
stopping time. Then you obtain the above inequality stopped with $\tau _{p}$%
. Take the expectation and use the Ito isometry to obtain
\begin{eqnarray*}
&&\int_{\Omega }\left( \sup_{t_{m}\in \mathcal{P}_{k}}\left\langle P\left(
t_{m}\wedge \tau _{p}\right) ,X\left( t_{m}\wedge \tau _{p}\right)
\right\rangle \right) dP \\
&\leq &E\left( \left\langle BX_{0},X_{0}\right\rangle \right) +2\left\vert
\left\vert Y\right\vert \right\vert _{K^{\prime }}\left\vert \left\vert
X_{k}^{r}\right\vert \right\vert _{K}
\end{eqnarray*}%
\begin{equation*}
+\left\Vert B\right\Vert
\sum_{j=0}^{m_{k}-1}\int_{t_{j}}^{t_{j+1}}\int_{\Omega }\left\vert
\left\vert Z\left( u\right) \right\vert \right\vert ^{2}dPdu
\end{equation*}%
\begin{equation*}
+2\int_{\Omega }\left( \sup_{t\in \left[ 0,T\right] }\left\vert
\int_{0}^{t\wedge \tau _{p}}\left( Z\circ J^{-1}\right) ^{\ast }\left(
P_{k}^{l}\right) ^{\tau _{p}}\circ JdW\right\vert \right) dP+E\left(
\left\vert e\left( k\right) \right\vert \right)
\end{equation*}%
\begin{equation}
\leq C+E\left( \left\vert e\left( k\right) \right\vert \right)
+2\int_{\Omega }\left( \sup_{t\in \left[ 0,T\right] }\left\vert
\int_{0}^{t}\left( Z\circ J^{-1}\right) ^{\ast }\left( P_{k}^{l}\right)
^{\tau _{p}}\circ JdW\right\vert \right) dP  \label{z25maye2h}
\end{equation}%
where the result of Lemma \ref{z23decl1g} that $X_{k}^{r}$ converges to $X$
in $K$ shows the term $2\left\vert \left\vert Y\right\vert \right\vert
_{K^{\prime }}\left\vert \left\vert X_{k}^{r}\right\vert \right\vert _{K}$
is bounded. Note that the constant $C$ can be assumed to be a continuous
function of
\begin{equation*}
\left\vert \left\vert Y\right\vert \right\vert _{K^{\prime }},\left\vert
\left\vert X\right\vert \right\vert _{K},\left\vert \left\vert Z\right\vert
\right\vert _{J},\left\Vert \left\langle BX_{0},X_{0}\right\rangle
\right\Vert _{L^{1}\left( \Omega \right) }
\end{equation*}%
which equals zero when all are equal to zero. (We can assume that $%
\left\vert \left\vert X_{k}^{r}\right\vert \right\vert _{K}\leq 2\left\vert
\left\vert X\right\vert \right\vert _{K}$ by taking a suitable subsequence
of the $\mathcal{P}_{K}$ if necessary.) The term involving the stochastic
integral is next.

Let $\mathcal{M}\left( t\right) =\int_{0}^{t}\left( Z\circ J^{-1}\right)
^{\ast }\left( P_{k}^{l}\right) ^{\tau _{p}}\circ JdW.$ Then from the
description of the quadratic variation,
\begin{equation*}
d\left[ \mathcal{M}\right] =\left\Vert \left( Z\circ J^{-1}\right) ^{\ast
}\left( P_{k}^{l}\right) ^{\tau _{p}}\circ J\right\Vert ^{2}ds
\end{equation*}

Applying the Burkholder Davis Gundy inequality, for $F\left( r\right) =r$ in
that stochastic integral,
\begin{equation*}
2\int_{\Omega }\left( \sup_{t\in \left[ 0,T\right] }\left\vert
\int_{0}^{t}\left( Z\circ J^{-1}\right) ^{\ast }\left( P_{k}^{l}\right)
^{\tau _{p}}\circ J\right\vert \right) dP
\end{equation*}%
\begin{equation}
\leq C\int_{\Omega }\left( \int_{0}^{T}\left\Vert \left( Z\circ
J^{-1}\right) ^{\ast }\left( P_{k}^{l}\right) ^{\tau _{p}}\circ J\right\Vert
_{\mathcal{L}_{2}\left( Q^{1/2}U,\mathbb{R}\right) }^{2}ds\right) ^{1/2}dP
\label{z25maye1h}
\end{equation}%
So let $\left\{ g_{i}\right\} $ be an orthonormal basis for $Q^{1/2}U$ and
consider the integrand in the above.
\begin{equation*}
\sum_{i=1}^{\infty }\left( \left( \left( Z\circ J^{-1}\right) ^{\ast }\left(
P_{k}^{l}\right) ^{\tau _{p}}\right) \left( J\left( g_{i}\right) \right)
\right) ^{2}=\sum_{i=1}^{\infty }\left\langle \left( P_{k}^{l}\right) ^{\tau
_{p}},Z\left( g_{i}\right) \right\rangle ^{2}
\end{equation*}%
From \ref{19july3},%
\begin{equation*}
\leq \sum_{i=1}^{\infty }\left\langle \left( P_{k}^{l}\right) ^{\tau
_{p}},\left( X_{k}^{l}\right) ^{\tau _{p}}\right\rangle \left\langle
BZ\left( g_{i}\right) ,Z\left( g_{i}\right) \right\rangle
\end{equation*}%
\begin{equation*}
\leq \left( \sup_{t_{m}\in \mathcal{P}_{k}}\left\langle \left(
P_{k}^{l}\right) ^{\tau _{p}}\left( t_{m}\right) ,\left( X_{k}^{l}\right)
^{\tau _{p}}\left( t_{m}\right) \right\rangle \right) \left\Vert
B\right\Vert \left\Vert Z\right\Vert _{\mathcal{L}_{2}}^{2}
\end{equation*}%
It follows that the integral in \ref{z25maye1h} is dominated by
\begin{equation*}
C\int_{\Omega }\sup_{t_{m}\in \mathcal{P}_{k}}\left\langle \left(
P_{k}^{l}\right) ^{\tau _{p}}\left( t_{m}\right) ,\left( X_{k}^{l}\right)
^{\tau _{p}}\left( t_{m}\right) \right\rangle ^{1/2}\left\Vert B\right\Vert
^{1/2}\left( \int_{0}^{T}\left\Vert Z\right\Vert _{\mathcal{L}%
_{2}}^{2}ds\right) ^{1/2}dP
\end{equation*}%
Now return to \ref{z25maye2h}. From what was just shown,%
\begin{equation*}
E\left( \sup_{t_{m}\in \mathcal{P}_{k}}\left\langle \left( P_{k}^{l}\right)
^{\tau _{p}}\left( t_{m}\right) ,\left( X_{k}^{l}\right) ^{\tau _{p}}\left(
t_{m}\right) \right\rangle \right)
\end{equation*}%
\begin{equation*}
\leq C+E\left( \left\vert e\left( k\right) \right\vert \right)
+2\int_{\Omega }\left( \sup_{t\in \left[ 0,T\right] }\left\vert
\int_{0}^{t}\left( Z\circ J^{-1}\right) ^{\ast }\left( P_{k}^{l}\right)
^{\tau _{p}}\circ JdW\right\vert \right) dP
\end{equation*}%
\begin{equation*}
\leq C+C\int_{\Omega }\sup_{t_{m}\in \mathcal{P}_{k}}\left\langle \left(
P_{k}^{l}\right) ^{\tau _{p}}\left( t_{m}\right) ,\left( X_{k}^{l}\right)
^{\tau _{p}}\left( t_{m}\right) \right\rangle ^{1/2}\left\Vert B\right\Vert
^{1/2}\left( \int_{0}^{T}\left\Vert Z\right\Vert _{\mathcal{L}%
_{2}}^{2}ds\right) ^{1/2}dP
\end{equation*}%
\begin{equation*}
+E\left( \left\vert e\left( k\right) \right\vert \right)
\end{equation*}%
\begin{eqnarray*}
&\leq &C+\frac{1}{2}E\left( \sup_{t_{m}\in \mathcal{P}_{k}}\left\langle
\left( P_{k}^{l}\right) ^{\tau _{p}}\left( t_{m}\right) ,\left(
X_{k}^{l}\right) ^{\tau _{p}}\left( t_{m}\right) \right\rangle \right) \\
&&+C\left\Vert Z\right\Vert _{\mathcal{L}^{2}\left( \left[ 0,T\right] \times
\Omega ,\mathcal{L}_{2}\right) }^{2}+E\left( \left\vert e\left( k\right)
\right\vert \right) .
\end{eqnarray*}%
It follows that
\begin{equation*}
\frac{1}{2}E\left( \sup_{t_{m}\in \mathcal{P}_{k}}\left\langle \left(
P_{k}^{l}\right) ^{\tau _{p}}\left( t_{m}\right) ,\left( X_{k}^{l}\right)
^{\tau _{p}}\left( t_{m}\right) \right\rangle \right) \leq C+E\left(
\left\vert e\left( k\right) \right\vert \right)
\end{equation*}%
Now let $p\rightarrow \infty $ and use the monotone convergence theorem to
obtain%
\begin{eqnarray*}
E\left( \sup_{t_{m}\in \mathcal{P}_{k}}\left\langle P_{k}^{l}\left(
t_{m}\right) ,X_{k}^{l}\left( t_{m}\right) \right\rangle \right) &=&E\left(
\sup_{t_{m}\in \mathcal{P}_{k}}\left\langle P\left( t_{m}\right) ,X\left(
t_{m}\right) \right\rangle \right) \\
&\leq &C+E\left( \left\vert e\left( k\right) \right\vert \right)
\end{eqnarray*}%
The monotone convergence theorem applies because $\tau _{p}$ merely
restricts the values $t_{m}\in \mathcal{P}_{k}$ which can be considered in
the above supremum. As mentioned above, this constant $C$ is a continuous
function of
\begin{equation*}
\left\vert \left\vert Y\right\vert \right\vert _{K^{\prime }},\left\vert
\left\vert X\right\vert \right\vert _{K},\left\vert \left\vert Z\right\vert
\right\vert _{J},\left\Vert \left\langle BX_{0},X_{0}\right\rangle
\right\Vert _{L^{1}\left( \Omega ,H\right) }
\end{equation*}%
and equals zero when all of these quantities equal 0. Also, for each $%
\varepsilon >0,$%
\begin{equation*}
E\left( \sup_{t_{m}\in \mathcal{P}_{k}}\left\langle P\left( t_{m}\right)
,X\left( t_{m}\right) \right\rangle \right) <C+\varepsilon
\end{equation*}%
whenever $k$ is large enough.

Let $\mathcal{D}$ denote the union of all the $\mathcal{P}_{k}.$ Thus $%
\mathcal{D}$ is a dense subset of $\left[ 0,T\right] $ and by the monotone
convergence theorem, it has just been shown, since the $\mathcal{P}_{k}$ are
nested, that for a constant $C$ dependent only on the above quantities,%
\begin{equation*}
E\left( \sup_{t\in \mathcal{D}}\left\langle P\left( t\right) ,X\left(
t\right) \right\rangle \right) \leq C+\varepsilon .
\end{equation*}%
Since $\varepsilon >0$ is arbitrary,
\begin{equation*}
E\left( \sup_{t\in \mathcal{D}}\left\langle P\left( t\right) ,X\left(
t\right) \right\rangle \right) =E\left( \sup_{t\in \mathcal{D}}\left\langle
BX\left( t\right) ,X\left( t\right) \right\rangle \right) \leq C
\end{equation*}%
This establishes \ref{9aug2h}. Now it follows right away that
\begin{equation}
\sup_{t\in \mathcal{D}}\left\langle P\left( t\right) ,X\left( t\right)
\right\rangle =\sup_{t\in \mathcal{D}}\left\langle BX\left( t\right)
,X\left( t\right) \right\rangle \leq C\left( \omega \right) <\infty \text{
a.e. }\omega  \label{9aug3h}
\end{equation}%
where $C=\int_{\Omega }C\left( \omega \right) dP$.

The function $t\rightarrow \sum_{k=1}^{\infty }\left\langle P\left( t\right)
,e_{k}\right\rangle ^{2}$ is obviously lower semicontinuous on $\left[ 0,T%
\right] $. This is because $t\rightarrow P\left( t\right) $ is continuous
into $V^{\prime }$. Now also for $t\in \mathcal{D}$ and $\omega \notin N,$ a
fixed set of measure zero which is defined in terms of $\mathcal{D}$,
\begin{equation*}
\left\langle P\left( t\right) ,X\left( t\right) \right\rangle =\left\langle
BX\left( t\right) ,X\left( t\right) \right\rangle =\sum_{k=1}^{\infty
}\left\langle BX\left( t\right) ,e_{k}\right\rangle ^{2}=\sum_{k=1}^{\infty
}\left\langle P\left( t\right) ,e_{k}\right\rangle ^{2}
\end{equation*}%
and so, for $t\in \mathcal{D}$ this infinite sum equals $\left\langle
P\left( t\right) ,X\left( t\right) \right\rangle ,$ and it was just shown
that $\sup_{t\in \mathcal{D}}\left\langle P\left( t\right) ,X\left( t\right)
\right\rangle \leq C\left( \omega \right) $. Hence, if $t$ is arbitrary and $%
t_{n}\rightarrow t$ for $t_{n}\in \mathcal{D}$, it follows from Fatou's
lemma that
\begin{equation*}
\sum_{k=1}^{\infty }\left\langle P\left( t\right) ,e_{k}\right\rangle
^{2}\leq \lim \inf_{n\rightarrow \infty }\sum_{k=1}^{\infty }\left\langle
P\left( t_{n}\right) ,e_{k}\right\rangle ^{2}=\lim \inf_{n\rightarrow \infty
}\left\langle P\left( t_{n}\right) ,X\left( t_{n}\right) \right\rangle \leq
C\left( \omega \right)
\end{equation*}%
and so%
\begin{equation*}
\sup_{t\in \mathcal{D}}\sum_{k=1}^{\infty }\left\langle P\left( t\right)
,e_{k}\right\rangle ^{2}=\sup_{t\in \left[ 0,T\right] }\sum_{k=1}^{\infty
}\left\langle P\left( t\right) ,e_{k}\right\rangle ^{2}\leq C\left( \omega
\right)
\end{equation*}%
This establishes \ref{9aug1h}.

Finally, consider the claim about weak continuity of $P\left( t\right) .$
From the above estimate, \ref{9aug3h},
\begin{equation*}
\sup_{t\in \mathcal{D}}\left\Vert BX\left( t\right) \right\Vert _{W^{\prime
}}<C\left( \omega \right) .
\end{equation*}%
Now let $t\in \left[ 0,T\right] .$ Then there exists a sequence $%
t_{n}\rightarrow t$ where $t_{n}\in \mathcal{D}$. Then for $\omega \notin N,$
$BX\left( t_{n}\right) $ is bounded in $W^{\prime }.$ But also $BX\left(
t_{n}\right) $ equals $P\left( t_{n}\right) $ for all $\omega $ off a single
set of measure zero, the one which came from $\mathcal{D}$, and so from the
definition of $P\left( t\right) ,$%
\begin{equation*}
P\left( t_{n}\right) \rightarrow P\left( t\right) \text{ in }V^{\prime }
\end{equation*}%
and also $P\left( t_{n}\right) $ is bounded in $W^{\prime }.$ Therefore,
there is a subsequence such that $P\left( t_{n}\right) \rightarrow \zeta $
weakly in $W^{\prime }.$ It follows $P\left( t\right) =\zeta $ and $%
\left\Vert P\left( t\right) \right\Vert _{W^{\prime }}\leq C\left( \omega
\right) $. Thus $t\rightarrow P\left( t\right) $ is continuous into $%
V^{\prime }$ and bounded in $W^{\prime }.$ If $t_{n}\rightarrow t,$ then if $%
P\left( t_{n}\right) $ fails to converge weakly to $P\left( t\right) $ in $%
W^{\prime },$ there would exist a subsequence still called $t_{n}$ such that
$P\left( t_{n}\right) $ converges weakly in $W^{\prime }$ to $\zeta \neq
P\left( t\right) $. However, $P\left( t_{n}\right) \rightarrow P\left(
t\right) $ in $V^{\prime }$ and so $\zeta =P\left( t\right) $ after all.

For $\omega \notin N$ the set of measure zero off which the above
computations were considered which came from the points of $\mathcal{D}$, it
was just shown that $P\left( t,\omega \right) \in W^{\prime }.$ Also, the
formula for $P\left( t,\omega \right) $ implies that this function is
progressively measurable into $V^{\prime }$. Therefore,
\begin{equation*}
\left( t,\omega \right) \rightarrow \left\langle P\left( t,\omega \right)
,v\right\rangle
\end{equation*}%
is progressively measurable if $v\in V$. Thus also, if $\omega \notin N,$
and $v_{n}\rightarrow w\in W,$%
\begin{equation*}
\left\langle P\left( t,\omega \right) ,v_{n}\right\rangle \rightarrow
\left\langle P\left( t,\omega \right) ,w\right\rangle .
\end{equation*}%
Since each $\mathcal{F}_{t}$ is assumed to be complete, this shows from the
Pettis theorem that $\left( t,\omega \right) \rightarrow P\left( t,\omega
\right) $ is progressively measurable. $\blacksquare $

Recall that for a.e. $t,$ $P\left( t,\omega \right) =BX\left( t,\omega
\right) $ for a.e. $\omega $. Also $BX$ is obviously progressively
measurable because $X$ is. However, it is not clear that $t\rightarrow
BX\left( t\right) $ is continuous into $W^{\prime }$ for all $t\in \left[ 0,T%
\right] $ off a set of measure zero. In a sense, $P\left( t,\omega \right) $
is filling in the missing values of $t$ retaining both progressive
measurability and continuity in $t$.

Consider the case where $B=I$ and $W=H$ so that the situation is that of a
Gelfand triple, $V\subseteq H=H^{\prime }\subseteq V^{\prime }$. Then in
this case, $P\left( t\right) =X\left( t\right) $ and the vectors $\left\{
e_{k}\right\} $ reduce to an orthonormal basis for $H$ such that each $%
e_{k}\in V$. Then the above sum
\begin{equation*}
\sum_{k}\left\langle P\left( t\right) ,e_{k}\right\rangle
^{2}=\sum_{k}\left\langle X\left( t\right) ,e_{k}\right\rangle
^{2}=\left\Vert X\left( t\right) \right\Vert _{H}^{2}=\left\langle P\left(
t\right) ,X\left( t\right) \right\rangle .
\end{equation*}%
Thus the main estimate in the above lemma would imply
\begin{equation*}
E\left( \sup_{t\in \left[ 0,T\right] }\left\Vert X\left( t\right)
\right\Vert _{H}^{2}\right) <C.
\end{equation*}%
This is the estimate in this special case which is found in \cite{pre07}. In
this special case, this is also the thing which will be of use in the study
of variational formulations for stochastic equations. However, in the case
considered here in which there is a possibly degenerate operator $B,$ it is
not clear that $\sum_{k}\left\langle P\left( t\right) ,e_{k}\right\rangle
^{2}=\left\langle P\left( t\right) ,X\left( t\right) \right\rangle $ for all
$t$. The vectors $\left\{ e_{k}\right\} $ are not necessarily an orthonormal
basis for $W$.

If there is interest in a more general conclusion which avoids the explicit
reference to $\mathcal{D}$, one can do the following. Let $\left\{
f_{i}\right\} $ be an orthonormal basis for $\left( BW\right) ^{\bot }$ in $%
W^{\prime }$. Then consider the function
\begin{equation*}
F\left( t\right) \equiv \sum_{k=1}^{\infty }\left\langle P\left( t\right)
,e_{k}\right\rangle ^{2}+\sum_{i}\left( P\left( t\right) ,f_{i}\right) ^{2}
\end{equation*}%
From the above lemma, $t\rightarrow P\left( t\right) $ is weakly continuous
into $W^{\prime }.$ Therefore, $F\left( t\right) $ is lower semicontinuous
for each $\omega $ off a single set of measure zero. It follows that for
such $\omega ,$
\begin{equation*}
\sup_{t\in \left[ 0,T\right] }F\left( t\right) =\sup_{t\in \mathcal{D}%
}F\left( t\right)
\end{equation*}%
For $t\in \mathcal{D}$ recall that $BX\left( t\right) =P\left( t\right) $
and so for such $t,$%
\begin{equation*}
F\left( t\right) =\sum_{k=1}^{\infty }\left\langle P\left( t\right)
,e_{k}\right\rangle ^{2}=\left\langle BX\left( t\right) ,X\left( t\right)
\right\rangle =\left\langle P\left( t\right) ,X\left( t\right) \right\rangle
\end{equation*}%
and this holds for all $\omega $ off a single set of measure zero, depending
on $\mathcal{D}$. Therefore, for such $\omega $ and the above lemma,%
\begin{equation*}
\sup_{t\in \left[ 0,T\right] }F\left( t\right) =\sup_{t\in \mathcal{D}%
}F\left( t\right) =\sup_{t\in \mathcal{D}}\left\langle P\left( t\right)
,X\left( t\right) \right\rangle \leq C\left( \omega \right) ,\ \int_{\Omega
}C\left( \omega \right) dP<\infty .
\end{equation*}%
It follows that
\begin{equation*}
E\left( \sup_{t\in \left[ 0,T\right] }F\left( t\right) \right) \leq C<\infty
\end{equation*}%
For any particular $t\notin S,$ we know that $BX\left( t\right) =P\left(
t\right) $ a.e. $\omega .$ Hence%
\begin{equation*}
F\left( t\right) =\sum_{k=1}^{\infty }\left\langle P\left( t\right)
,e_{k}\right\rangle ^{2}=\left\langle BX\left( t\right) ,X\left( t\right)
\right\rangle =\left\langle P\left( t\right) ,X\left( t\right) \right\rangle
\text{ a.e. }\omega
\end{equation*}%
Therefore, for that $t\notin S,$
\begin{equation*}
E\left( \left\langle P\left( t\right) ,X\left( t\right) \right\rangle
\right) =E\left( F\left( t\right) \right) \leq E\left( \sup_{t\in \left[ 0,T%
\right] }F\left( t\right) \right) \leq C<\infty
\end{equation*}

\begin{corollary}
For $C$ in the above lemma, and for any $t\notin S,$%
\begin{equation*}
E\left( \left\langle P\left( t\right) ,X\left( t\right) \right\rangle
\right) \leq C.
\end{equation*}
\end{corollary}

\section{\protect\vspace{1pt}A Simplification Of The Formula\label%
{technicalsimplification}}

This lemma also provides a way to simplify \ref{z24jane12g}. First suppose $%
X_{0}\in L^{p}\left( \Omega ,V\right) $ so that $X-X_{0}\in L^{p}\left( %
\left[ 0,T\right] \times \Omega ,V\right) $. Refer to \ref{z24jane12g}. One
term there is%
\begin{equation*}
\left\langle P\left( t_{1}\right) -BX_{0}-BM\left( t_{1}\right) ,X\left(
t_{1}\right) -X_{0}-M\left( t_{1}\right) \right\rangle
\end{equation*}%
It equals
\begin{equation*}
\left\langle B\left( X\left( t_{1}\right) -X_{0}-M\left( t_{1}\right)
\right) ,X\left( t_{1}\right) -X_{0}-M\left( t_{1}\right) \right\rangle
\end{equation*}%
\begin{equation*}
\leq 2\left\langle B\left( X\left( t_{1}\right) -X_{0}\right) ,X\left(
t_{1}\right) -X_{0}\right\rangle +2\left\langle BM\left( t_{1}\right)
,M\left( t_{1}\right) \right\rangle
\end{equation*}%
\begin{equation*}
=2\left\langle P\left( t_{1}\right) -BX_{0},X\left( t_{1}\right)
-X_{0}\right\rangle +2\left\langle BM\left( t_{1}\right) ,M\left(
t_{1}\right) \right\rangle
\end{equation*}%
It was observed above that $2\left\langle BM\left( t_{1}\right) ,M\left(
t_{1}\right) \right\rangle \rightarrow 0$ a.e. and also in $L^{1}\left(
\Omega \right) $ as $k\rightarrow \infty $. Apply the above lemma for $X$
replaced with $\bar{X}\left( t\right) \equiv X\left( t\right) -X_{0}$ and
denote the resulting $P\left( t\right) $ by $\bar{P}\left( t\right) .$ The
new $X_{0}$ equals 0. Also use $\left[ 0,t_{1}\right] $ instead of $\left[
0,T\right] .$ Thus the above reduces with this new $X$ to
\begin{equation*}
2\left\langle \bar{P}\left( t_{1}\right) ,\bar{X}\left( t_{1}\right)
\right\rangle +2\left\langle BM\left( t_{1}\right) ,M\left( t_{1}\right)
\right\rangle
\end{equation*}%
From the above lemma,
\begin{equation*}
E\left( \left\langle \bar{P}\left( t_{1}\right) ,\bar{X}\left( t_{1}\right)
\right\rangle \right) <C\left( \left\vert \left\vert Y\right\vert
\right\vert _{K_{t_{1}}^{\prime }},\left\vert \left\vert \bar{X}\right\vert
\right\vert _{K_{t_{1}}},\left\vert \left\vert Z\right\vert \right\vert
_{J_{t_{1}}}\right)
\end{equation*}%
where in the definitions of $K,K^{\prime },J,$ replace $\left[ 0,T\right] $
with $\left[ 0,t\right] $ and let the resulting spaces be denoted by $%
K_{t},K_{t}^{\prime },J_{t}$. Therefore, this term converges to 0 in $%
L^{1}\left( \Omega \right) $ as $k\rightarrow \infty $. In addition to this,
the term converges to 0 pointwise for a.e. $\omega $ after passing to a
suitable subsequence. Thus we can enlarge $e\left( k\right) $ and neglect
the last term of \ref{z24jane12g}.

Then, it would follow from \ref{z24jane12g},%
\begin{equation*}
\left\langle P\left( t_{m}\right) ,X\left( t_{m}\right) \right\rangle
-\left\langle BX_{0},X_{0}\right\rangle =e\left( k\right)
+2\int_{0}^{t_{m}}\left\langle Y\left( u\right) ,X_{k}^{r}\left( u\right)
\right\rangle du+
\end{equation*}%
\begin{equation*}
+2\int_{0}^{t_{m}}\left( Z\circ J^{-1}\right) ^{\ast }P_{k}^{l}\circ JdW
\end{equation*}%
\begin{equation*}
+\sum_{j=0}^{m-1}\left\langle B\left( M\left( t_{j+1}\right) -M\left(
t_{j}\right) \right) ,M\left( t_{j+1}\right) -M\left( t_{j}\right)
\right\rangle
\end{equation*}%
\begin{equation}
-\sum_{j=1}^{m-1}\left\langle \Delta P\left( t_{j}\right) -B\Delta M\left(
t_{j}\right) ,\Delta X\left( t_{j}\right) -\Delta M\left( t_{j}\right)
\right\rangle  \label{z3febe1h}
\end{equation}%
where $e\left( k\right) \rightarrow 0$ in $L^{1}\left( \Omega \right) $ and
a.e. $\omega $ and
\begin{equation*}
\Delta X\left( t_{j}\right) \equiv X\left( t_{j+1}\right) -X\left(
t_{j}\right) ,
\end{equation*}%
$\Delta M\left( t_{j}\right) $ being defined similarly.

Can you obtain this equation even in case $X_{0}$ is not assumed to be in $%
L^{p}\left( \Omega ,V\right) ?$ Let $X_{0k}\in L^{p}\left( \Omega ,V\right)
\cap L^{2}\left( \Omega ,W\right) ,X_{0k}\rightarrow X_{0}$ in $L^{2}\left(
\Omega ,W\right) .$
\begin{eqnarray*}
\left\langle P\left( t_{1}\right) -BX_{0},X\left( t_{1}\right)
-X_{0}\right\rangle ^{1/2} &\leq &\left\langle P\left( t_{1}\right)
-BX_{0k},X\left( t_{1}\right) -X_{0k}\right\rangle ^{1/2} \\
&&+\left\langle B\left( X_{0k}-X_{0}\right) ,X_{0k}-X_{0}\right\rangle ^{1/2}
\end{eqnarray*}%
Also, restoring the superscript to identify the partition,
\begin{equation*}
P\left( t_{1}^{k}\right) -BX_{0k}=B\left( X_{0}-X_{0k}\right)
+\int_{0}^{t_{1}^{k}}Y\left( s\right) ds+B\int_{0}^{t_{1}^{k}}Z\left(
s\right) dW.
\end{equation*}%
Of course $\left\Vert X-X_{0k}\right\Vert _{K}$ is not bounded, but for each
$k$ it is finite. Let $n_{k}$ denote a subsequence of $\left\{ k\right\} $
such that
\begin{equation*}
\left\Vert X-X_{0k}\right\Vert _{K_{t_{1}^{n_{k}}}}<1/k.
\end{equation*}%
Then from the above Lemma \ref{z24janl1g},%
\begin{equation*}
E\left( \left\langle P\left( t_{1}^{n_{k}}\right) -BX_{0k},X\left(
t_{1}^{n_{k}}\right) -X_{0k}\right\rangle \right)
\end{equation*}%
\begin{equation*}
\leq C\left( \left\vert \left\vert Y\right\vert \right\vert
_{K_{t_{1}^{n_{k}}}^{\prime }},\left\Vert X-X_{0k}\right\Vert
_{K_{t_{1}^{n_{k}}}},\left\vert \left\vert Z\right\vert \right\vert
_{J_{t_{1}^{n_{k}}}},\left\langle B\left( X_{0}-X_{0k}\right)
,X_{0}-X_{0k}\right\rangle _{L^{1}\left( \Omega \right) }\right)
\end{equation*}%
\begin{equation*}
\leq C\left( \left\vert \left\vert Y\right\vert \right\vert
_{K_{t_{1}^{n_{k}}}^{\prime }},\frac{1}{k},\left\vert \left\vert
Z\right\vert \right\vert _{J_{t_{1}^{n_{k}}}},\left\langle B\left(
X_{0}-X_{0k}\right) ,X_{0}-X_{0k}\right\rangle _{L^{1}\left( \Omega \right)
}\right)
\end{equation*}%
Hence
\begin{equation*}
E\left( \left\langle P\left( t_{1}^{n_{k}}\right) -BX_{0},X\left(
t_{1}^{n_{k}}\right) -X_{0}\right\rangle \right)
\end{equation*}%
\begin{equation*}
\leq 2E\left( \left\langle P\left( t_{1}^{n_{k}}\right) -BX_{0k},X\left(
t_{1}^{n_{k}}\right) -X_{0k}\right\rangle \right) +2E\left( \left\langle
B\left( X_{0k}-X_{0}\right) ,X_{0k}-X_{0}\right\rangle \right)
\end{equation*}%
\begin{eqnarray*}
&\leq &2C\left( \left\vert \left\vert Y\right\vert \right\vert
_{K_{t_{1}^{n_{k}}}^{\prime }},\frac{1}{k},\left\vert \left\vert
Z\right\vert \right\vert _{J_{t_{1}^{n_{k}}}},\left\langle B\left(
X_{0}-X_{0k}\right) ,X_{0}-X_{0k}\right\rangle _{L^{1}\left( \Omega \right)
}\right) \\
&&+2\left\Vert B\right\Vert \left\Vert X_{0k}-X_{0}\right\Vert _{L^{2}\left(
\Omega ,W\right) }^{2}
\end{eqnarray*}%
which converges to 0 as $k\rightarrow \infty $. It follows that there exists
a suitable subsequence such that \ref{z3febe1h} holds even in the case that $%
X_{0}$ is only known to be in $L^{2}\left( \Omega ,W\right) $. From now on,
assume this subsequence for the partitions $\mathcal{P}_{k}$. Thus $k $ will
really be $n_{k}$ and it suffices to consider the limit as $k\rightarrow
\infty $ of the equation of \ref{z3febe1h}. To emphasize this point again,
the reason for the above observations is to neglect
\begin{equation*}
\left\langle P\left( t_{1}\right) -BX_{0}-BM\left( t_{1}\right) ,X\left(
t_{1}\right) -X_{0}-M\left( t_{1}\right) \right\rangle
\end{equation*}%
in passing to the limit as $k\rightarrow \infty $ provided a suitable
subsequence is used.

In order to eventually obtain the Ito formula \ref{z27jane2g}, there is a
technical result which will be needed. It was mostly proved in Lemma \ref%
{z24janl9gaa}.

\begin{lemma}
\label{z24janl9g}Let $P\left( s\right) -P_{k}^{l}\left( s\right) \equiv
\Delta _{k}\left( s\right) .$ Then the following limit occurs.%
\begin{equation*}
\lim_{k\rightarrow \infty }P\left( \left[ \sup_{t\in \left[ 0,T\right]
}\left\vert \int_{0}^{t}\left( Z\circ J^{-1}\right) ^{\ast }\Delta _{k}\circ
JdW\right\vert \geq \varepsilon \right] \right) =0
\end{equation*}%
The stochastic integral%
\begin{equation*}
\int_{0}^{t}\left( Z\circ J^{-1}\right) ^{\ast }P\circ JdW
\end{equation*}%
makes sense because $BX=P$ is $W^{\prime }$ progressively measurable. Also,
there exists a further subsequence, still denoted as $k$ such that%
\begin{equation*}
\int_{0}^{t}\left( Z\circ J^{-1}\right) ^{\ast }P_{k}^{l}\circ
JdW\rightarrow \int_{0}^{t}\left( Z\circ J^{-1}\right) ^{\ast }P\circ JdW
\end{equation*}%
uniformly on $\left[ 0,T\right] $ for $a.e.$ $\omega $.
\end{lemma}

\textbf{Proof: }This follows from Lemma \ref{z24janl9gaa}. The last
conclusion follows from the usual use of the Borel Cantelli lemma, Ito
formula, and the maximal inequalities for submartingales. It was shown in
this lemma that
\begin{equation*}
\lim_{k\rightarrow \infty }P\left( \sup_{t\in \left[ 0,T\right] }\left\vert
\int_{0}^{t}\left( Z\circ J^{-1}\right) ^{\ast }\Delta _{k}\circ
JdW\right\vert >\varepsilon \right) =0
\end{equation*}%
Thus, one can obtain the existence of a subsequence, still denoted as $k$
such that
\begin{equation*}
P\left( \sup_{t\in \left[ 0,T\right] }\left\vert \int_{0}^{t}\left( Z\circ
J^{-1}\right) ^{\ast }\Delta _{k}\circ JdW\right\vert >2^{-k}\right) <2^{-k}
\end{equation*}%
and then uniform convergence is obtained for this subsequence off a set of
meaure zero. $\blacksquare $

From now on, the sequence will either be this subsequence or a further
subsequence. Also $N$ will be enlarged so that for $\omega \notin N,$ the
above uniform convergence of the stochastic integrals takes place in
addition to the other items above.

\section{The Ito Formula\label{itoformula}}

The next lemma is the Ito formula for $t\in \mathcal{D},$ the dense subset
consisting of all the mesh points of all partitions $\mathcal{P}_{k}$.

\begin{proposition}
\label{z29janl1g}Let $X_{0}\in L^{2}\left( \Omega ;W\right) $ and be $%
\mathcal{F}_{0}$ measurable. There exists a dense subset of $\left[ 0,T%
\right] $ denoted as $\mathcal{D}$ such that \ for every $t\in \mathcal{D}$,
\begin{equation*}
\left\langle P\left( t\right) ,X\left( t\right) \right\rangle =\left\langle
BX_{0},X_{0}\right\rangle +\int_{0}^{t}\left( 2\left\langle Y\left( s\right)
,X\left( s\right) \right\rangle +\left\langle BZ,Z\right\rangle _{\mathcal{L}%
_{2}}\right) ds
\end{equation*}%
\begin{equation}
+2\int_{0}^{t}\left( Z\circ J^{-1}\right) ^{\ast }P\circ JdW
\label{z29jane1g}
\end{equation}%
where in the above formula,
\begin{equation*}
\left\langle BZ,Z\right\rangle _{\mathcal{L}_{2}}\equiv \left(
R^{-1}BZ,Z\right) _{\mathcal{L}_{2}\left( Q^{1/2}U,W\right) }
\end{equation*}%
for $R$ the Riesz map from $W$ to $W^{\prime }$. In addition to this, for
all such $t\in \mathcal{D},$ $E\left( \left\langle BX\left( t\right)
,X\left( t\right) \right\rangle \right) =$%
\begin{equation}
E\left( \left\langle P\left( t\right) ,X\left( t\right) \right\rangle
\right) =E\left( \left\langle BX_{0},X_{0}\right\rangle \right) +E\left(
\int_{0}^{t}2\left\langle Y\left( s\right) ,X\left( s\right) \right\rangle
+\left\langle BZ,Z\right\rangle _{\mathcal{L}_{2}}ds\right)  \label{26july2h}
\end{equation}%
In addition to this,
\begin{eqnarray}
&&E\left( \sup_{t\in \mathcal{D}}\left\langle P\left( t\right) ,X\left(
t\right) \right\rangle \right)  \notag \\
&<&C\left( \left\vert \left\vert Y\right\vert \right\vert _{K^{\prime
}},\left\vert \left\vert X\right\vert \right\vert _{K},\left\vert \left\vert
Z\right\vert \right\vert _{J},\left\Vert \left\langle
BX_{0},X_{0}\right\rangle \right\Vert _{L^{1}\left( \Omega \right) }\right)
<\infty .  \label{6aug3h}
\end{eqnarray}
\end{proposition}

Note first that for $\left\{ g_{i}\right\} $ an orthonormal basis for $%
Q^{1/2}\left( U\right) ,$%
\begin{equation*}
\left( R^{-1}BZ,Z\right) _{\mathcal{L}_{2}}\equiv \sum_{i}\left(
R^{-1}BZ\left( g_{i}\right) ,Z\left( g_{i}\right) \right)
_{W}=\sum_{i}\left\langle BZ\left( g_{i}\right) ,Z\left( g_{i}\right)
\right\rangle _{W^{\prime }W}\geq 0
\end{equation*}

\textbf{Proof: }Inequality \ref{6aug3h} follows from the earlier lemma. In
the situation of \ref{z24jane1g}, let $\mathcal{D}$ be the union of the $%
\mathcal{P}_{k}$ described above as the union of all positive mesh points
less than $T$ for all the $\mathcal{P}_{k}.$Then, since these $\mathcal{P}%
_{k}$ are nested, if $t\in \mathcal{D}$, then $t\in \mathcal{P}_{k}$ for all
$k$ large enough. Consider \ref{z3febe1h},
\begin{equation*}
\left\langle P\left( t\right) ,X\left( t\right) \right\rangle -\left\langle
BX_{0},X_{0}\right\rangle =e\left( k\right) +2\int_{0}^{t}\left\langle
Y\left( u\right) ,X_{k}^{r}\left( u\right) \right\rangle du
\end{equation*}%
\begin{equation*}
+2\int_{0}^{t}\left( Z\circ J^{-1}\right) ^{\ast }P_{k}^{l}\circ
JdW+\sum_{j=0}^{q_{k}-1}\left\langle B\left( M\left( t_{j+1}\right) -M\left(
t_{j}\right) \right) ,M\left( t_{j+1}\right) -M\left( t_{j}\right)
\right\rangle
\end{equation*}%
\begin{equation}
-\sum_{j=1}^{q_{k}-1}\left\langle \Delta P\left( t_{j}\right) -\Delta
BM\left( t_{j}\right) ,\Delta X\left( t_{j}\right) -\Delta M\left(
t_{j}\right) \right\rangle  \label{z27jane7g}
\end{equation}%
where $t_{q_{k}}=t,\ \Delta X\left( t_{j}\right) =X\left( t_{j+1}\right)
-X\left( t_{j}\right) $ and $e\left( k\right) \rightarrow 0$ in probability.
By Lemma \ref{z24janl9g} the stochastic integral on the right converges
uniformly to
\begin{equation*}
2\int_{0}^{t}\left( Z\circ J^{-1}\right) ^{\ast }P\circ JdW
\end{equation*}%
off a set of measure zero. The deterministic integral on the right converges
to
\begin{equation*}
2\int_{0}^{t}\left\langle Y\left( u\right) ,X\left( u\right) \right\rangle du
\end{equation*}%
in $L^{1}\left( \Omega \right) $ because $X_{k}^{r}\rightarrow X$ in $K$.
Then
\begin{eqnarray*}
&&P\left( \left[ \sup_{t\in \left[ 0,T\right] }\left\vert
\int_{0}^{t}\left\langle Y\left( u\right) ,X\left( u\right) \right\rangle
du-\int_{0}^{t}\left\langle Y\left( u\right) ,X_{k}^{r}\left( u\right)
\right\rangle du\right\vert >\lambda \right] \right) \\
&\leq &P\left( \left[ \sup_{t\in \left[ 0,T\right] }\int_{0}^{t}\left\vert
\left\langle Y\left( u\right) ,X\left( u\right) -X_{k}^{r}\left( u\right)
\right\rangle \right\vert du>\lambda \right] \right) \\
&\leq &P\left( \left[ \int_{0}^{T}\left\vert \left\langle Y\left( u\right)
,X\left( u\right) -X_{k}^{r}\left( u\right) \right\rangle \right\vert
du>\lambda \right] \right) \\
&\leq &\frac{1}{\lambda }\int_{\Omega }\int_{0}^{T}\left\Vert Y\left(
u\right) \right\Vert \left\Vert X\left( u\right) -X_{k}^{r}\left( u\right)
\right\Vert dudP\leq \frac{1}{\lambda }\left\Vert Y\right\Vert _{K^{\prime
}}\left\Vert X-X_{k}^{r}\right\Vert _{K}
\end{eqnarray*}%
Since $\left\Vert X-X_{k}^{r}\right\Vert _{K}\rightarrow 0$ as $k\rightarrow
\infty ,$ it follows that there is a further subsequence, $n_{k}$ such that
\begin{equation*}
P\left( \left[ \sup_{t\in \left[ 0,T\right] }\left\vert
\int_{0}^{t}\left\langle Y\left( u\right) ,X\left( u\right) \right\rangle
du-\int_{0}^{t}\left\langle Y\left( u\right) ,X_{n_{k}}^{r}\left( u\right)
\right\rangle du\right\vert >2^{-k}\right] \right) \leq 2^{-k}.
\end{equation*}%
To save notation, refer to the subsequence as $k$. Thus, for a suitable
subsequence,
\begin{equation*}
\lim_{k\rightarrow \infty }\int_{0}^{t}\left\langle Y\left( u\right)
,X_{k}^{r}\left( u\right) \right\rangle du=\int_{0}^{t}\left\langle Y\left(
u\right) ,X\left( u\right) \right\rangle du
\end{equation*}%
uniformly off some set of measure zero. Consider the fourth term. It equals
\begin{equation}
\sum_{j=0}^{q_{k}-1}\left( R^{-1}B\left( M\left( t_{j+1}\right) -M\left(
t_{j}\right) \right) ,M\left( t_{j+1}\right) -M\left( t_{j}\right) \right)
_{W}  \label{z1june1h}
\end{equation}%
where $R$ is the Riesz map from $W$ to $W^{\prime }$. This equals
\begin{eqnarray*}
&&\frac{1}{4}\left( \sum_{j=0}^{q_{k}-1}\left\Vert R^{-1}BM\left(
t_{j+1}\right) +M\left( t_{j+1}\right) -\left( R^{-1}BM\left( t_{j}\right)
+M\left( t_{j}\right) \right) \right\Vert ^{2}\right. \\
&&-\left. \sum_{j=0}^{q_{k}-1}\left\Vert R^{-1}BM\left( t_{j+1}\right)
-M\left( t_{j+1}\right) -\left( R^{-1}BM\left( t_{j}\right) -M\left(
t_{j}\right) \right) \right\Vert ^{2}\right)
\end{eqnarray*}%
From Theorem \ref{31mayt1h}, as $k\rightarrow \infty ,$ the above converges
in probability to ($t_{q_{k}}=t$)
\begin{equation*}
\frac{1}{4}\left( \left[ R^{-1}BM+M\right] \left( t\right) -\left[ R^{-1}BM-M%
\right] \left( t\right) \right)
\end{equation*}%
However, from the well known description of the quadratic variation of a
martingale, the above equals
\begin{equation*}
\frac{1}{4}\left( \int_{0}^{t}\left\Vert R^{-1}BZ+Z\right\Vert _{\mathcal{L}%
_{2}}^{2}ds-\int_{0}^{t}\left\Vert R^{-1}BZ-Z\right\Vert _{\mathcal{L}%
_{2}}^{2}ds\right)
\end{equation*}%
which equals
\begin{equation*}
\int_{0}^{t}\left( R^{-1}BZ,Z\right) _{\mathcal{L}_{2}}ds\equiv
\int_{0}^{t}\left\langle BZ,Z\right\rangle _{\mathcal{L}_{2}}ds
\end{equation*}%
This is what was desired.

Note that in the case of a Gelfand triple, when $W=H=H^{\prime },$ the term $%
\left\langle BZ,Z\right\rangle _{\mathcal{L}_{2}}$ will end up reducing to
nothing more than $\left\Vert Z\right\Vert _{\mathcal{L}_{2}}^{2}$.

Thus all the terms in \ref{z27jane7g} converge in probability except for the
last term which also must converge in probability because it equals the sum
of terms which do. It remains to find what this last term converges to. Thus
\begin{equation*}
\left\langle P\left( t\right) ,X\left( t\right) \right\rangle -\left\langle
BX_{0},X_{0}\right\rangle =2\int_{0}^{t}\left\langle Y\left( u\right)
,X\left( u\right) \right\rangle du
\end{equation*}%
\begin{equation*}
+2\int_{0}^{t}\left( Z\circ J^{-1}\right) ^{\ast }P\circ
JdW+\int_{0}^{t}\left\langle BZ,Z\right\rangle _{\mathcal{L}_{2}}ds-a
\end{equation*}%
where $a$ is the limit in probability of the term
\begin{equation}
\sum_{j=1}^{q_{k}-1}\left\langle \Delta P\left( t_{j}\right) -\Delta
BM\left( t_{j}\right) ,\Delta X\left( t_{j}\right) -\Delta M\left(
t_{j}\right) \right\rangle  \label{z1june2h}
\end{equation}%
Let $\pi _{n}$ be the projection onto $\mbox{span}\left( e_{1},\cdots
,e_{n}\right) $ where $\left\{ e_{k}\right\} $ is a complete orthonormal
basis for $W$ with each $e_{k}\in V$. Then using
\begin{equation*}
P\left( t_{j+1}\right) -P\left( t_{j}\right) -\left( BM\left( t_{j+1}\right)
-BM\left( t_{j}\right) \right) =\int_{t_{j}}^{t_{j+1}}Y\left( s\right) ds,
\end{equation*}%
the troublesome term of \ref{z1june2h} above is of the form%
\begin{equation*}
\sum_{j=1}^{q_{k}-1}\int_{t_{j}}^{t_{j+1}}\left\langle Y\left( s\right)
,\Delta X\left( t_{j}\right) -\Delta M\left( t_{j}\right) \right\rangle ds
\end{equation*}%
\begin{eqnarray*}
&=&\sum_{j=1}^{q_{k}-1}\int_{t_{j}}^{t_{j+1}}\left\langle Y\left( s\right)
,\Delta X\left( t_{j}\right) -\pi _{n}\Delta M\left( t_{j}\right)
\right\rangle ds \\
&&+\sum_{j=1}^{q_{k}-1}\int_{t_{j}}^{t_{j+1}}\left\langle Y\left( s\right)
,-\left( I-\pi _{n}\right) \Delta M\left( t_{j}\right) \right\rangle ds
\end{eqnarray*}%
which equals
\begin{equation}
\sum_{j=1}^{q_{k}-1}\int_{t_{j}}^{t_{j+1}}\left\langle Y\left( s\right)
,X\left( t_{j+1}\right) -X\left( t_{j}\right) -\pi _{n}\left( M\left(
t_{j+1}\right) -M\left( t_{j}\right) \right) \right\rangle ds
\label{z24jane15g}
\end{equation}%
\begin{equation}
+\sum_{j=1}^{q_{k}-1}\left\langle \Delta P\left( t_{j}\right) -B\Delta
M\left( t_{j}\right) ,-\left( I-\pi _{n}\right) \left( M\left(
t_{j+1}\right) -M\left( t_{j}\right) \right) \right\rangle
\label{z24jane16g}
\end{equation}%
Since $P\left( t\right) =BX\left( t\right) $ for the $t$ of interest in the
above, the Cauchy Schwarz inequality implies the term of \ref{z24jane16g} is
dominated by
\begin{equation*}
\left( \sum_{j=1}^{q_{k}-1}\left\langle \Delta P\left( t_{j}\right) -\Delta
BM\left( t_{j}\right) ,\left( \Delta X\left( t_{j}\right) -\Delta M\left(
t_{j}\right) \right) \right\rangle \right) ^{1/2}\cdot
\end{equation*}%
\begin{equation}
\left( \sum_{j=1}^{q_{k}-1}\left\vert \left\langle B\left( \ I-\pi
_{n}\right) \Delta M\left( t_{j}\right) ,\left( \ I-\pi _{n}\right) \Delta
M\left( t_{j}\right) \right\rangle \right\vert ^{2}\right) ^{1/2}
\label{z27jane9g}
\end{equation}%
Now it is known that $\sum_{j=1}^{q_{k}-1}\left\langle \Delta P\left(
t_{j}\right) -\Delta BM\left( t_{j}\right) ,\left( \Delta X\left(
t_{j}\right) -\Delta M\left( t_{j}\right) \right) \right\rangle $ converges
in probability to $a\geq 0.$ If you take the expectation of the square of
the other factor above, it is no larger than
\begin{equation*}
\left\Vert B\right\Vert E\left( \sum_{j=1}^{q_{k}-1}\left\Vert \left( \
I-\pi _{n}\right) \Delta M\left( t_{j}\right) \right\Vert _{W}^{2}\right)
\end{equation*}%
\begin{equation*}
=\left\Vert B\right\Vert E\left( \sum_{j=1}^{q_{k}-1}\left\Vert \left( \
I-\pi _{n}\right) \int_{t_{j}}^{t_{j+1}}Z\left( s\right) dW\left( s\right)
\right\Vert _{W}^{2}\right)
\end{equation*}%
\begin{equation*}
=\left\Vert B\right\Vert \sum_{j=1}^{q_{k}-1}E\left( \left\Vert
\int_{t_{j}}^{t_{j+1}}\left( \ I-\pi _{n}\right) Z\left( s\right) dW\left(
s\right) \right\Vert ^{2}\right)
\end{equation*}%
\begin{eqnarray*}
&=&\left\Vert B\right\Vert \sum_{j=1}^{q_{k}-1}E\left(
\int_{t_{j}}^{t_{j+1}}\left\vert \left\vert \left( \ I-\pi _{n}\right)
Z\left( s\right) \right\vert \right\vert _{\mathcal{L}_{2}\left(
Q^{1/2}U,W\right) }^{2}ds\right) \\
&\leq &\left\Vert B\right\Vert E\left( \int_{0}^{T}\left\vert \left\vert
\left( \ I-\pi _{n}\right) Z\left( s\right) \right\vert \right\vert _{%
\mathcal{L}_{2}\left( Q^{1/2}U,W\right) }^{2}ds\right)
\end{eqnarray*}%
Letting $\left\{ g_{i}\right\} $ be an orthonormal basis for $Q^{1/2}U,$%
\begin{equation}
=\left\Vert B\right\Vert \int_{\Omega }\int_{0}^{T}\sum_{i=1}^{\infty
}\left\Vert \left( \ I-\pi _{n}\right) Z\left( s\right) \left( g_{i}\right)
\right\Vert _{W}^{2}dsdP  \label{z3june1h}
\end{equation}%
The integrand $\sum_{i=1}^{\infty }\left\Vert \left( \ I-\pi _{n}\right)
Z\left( s\right) \left( g_{i}\right) \right\Vert _{W}^{2}$ converges to 0.
Also, it is dominated by%
\begin{equation*}
\sum_{i=1}^{\infty }\left\Vert Z\left( s\right) \left( g_{i}\right)
\right\Vert _{W}^{2}\equiv \left\Vert Z\right\Vert _{\mathcal{L}_{2}\left(
Q^{1/2}U,W\right) }^{2}
\end{equation*}%
which is given to be in $L^{1}\left( \left[ 0,T\right] \times \Omega \right)
.$ Therefore, from the dominated convergence theorem, the expression in \ref%
{z3june1h} converges to 0 as $n\rightarrow \infty $ independent of $k$.

Thus the expression in \ref{z27jane9g} is of the form $f_{k}g_{nk}$ where $%
f_{k}$ converges in probability to $a^{1/2}$ as $k\rightarrow \infty $ and $%
g_{nk}$ converges in probability to 0 as $n\rightarrow \infty $ independent
of $k.$ Now this implies $f_{k}g_{nk}$ converges in probability to 0. Here
is why.
\begin{eqnarray*}
P\left( \left[ \left\vert f_{k}g_{nk}\right\vert >\varepsilon \right]
\right) &\leq &P\left( 2\delta \left\vert f_{k}\right\vert >\varepsilon
\right) +P\left( 2C_{\delta }\left\vert g_{nk}\right\vert >\varepsilon
\right) \\
&\leq &P\left( 2\delta \left\vert f_{k}-a^{1/2}\right\vert +2\delta
\left\vert a^{1/2}\right\vert >\varepsilon \right) +P\left( 2C_{\delta
}\left\vert g_{nk}\right\vert >\varepsilon \right)
\end{eqnarray*}%
where $\delta \left\vert f_{k}\right\vert +C_{\delta }\left\vert
g_{kn}\right\vert >\left\vert f_{k}g_{nk}\right\vert $ and $\lim_{\delta
\rightarrow 0}C_{\delta }=\infty .$ Pick $\delta $ small enough that $%
\varepsilon -2\delta a^{1/2}>\varepsilon /2.$ Then this is dominated by
\begin{equation*}
\leq P\left( 2\delta \left\vert f_{k}-a^{1/2}\right\vert >\varepsilon
/2\right) +P\left( 2C_{\delta }\left\vert g_{nk}\right\vert >\varepsilon
\right)
\end{equation*}%
Fix $n$ large enough that the second term is less than $\eta $ for all $k.$
Now taking $k$ large enough, the above is less than $\eta $. It follows the
expression in \ref{z27jane9g} and co nsequently in \ref{z24jane16g}
converges to 0 in probability.

Now consider the other term \ref{z24jane15g} using the $n$ just determined.
This term is of the form
\begin{eqnarray*}
&&\sum_{j=1}^{q_{k}-1}\int_{t_{j}}^{t_{j+1}}\left\langle Y\left( s\right)
,X_{k}^{r}\left( s\right) -X_{k}^{l}\left( s\right) -\pi _{n}\left(
M_{k}^{r}\left( s\right) -M_{k}^{l}\left( s\right) \right) \right\rangle ds
\\
&=&\int_{t_{1}}^{t}\left\langle Y\left( s\right) ,X_{k}^{r}\left( s\right)
-X_{k}^{l}\left( s\right) -\pi _{n}\left( M_{k}^{r}\left( s\right)
-M_{k}^{l}\left( s\right) \right) \right\rangle ds
\end{eqnarray*}%
where $M_{k}^{r}$ denotes the step function
\begin{equation*}
M_{k}^{r}\left( t\right) =\sum_{i=0}^{m_{k}-1}M\left( t_{i+1}\right)
\mathcal{X}_{(t_{i},t_{i+1}]}\left( t\right)
\end{equation*}%
with $M_{k}^{l}$ defined similarly as a step function featuring the value of
$M$ at the left end of each interval. The term
\begin{equation*}
\int_{t_{1}}^{t}\left\langle Y\left( s\right) ,\pi _{n}\left(
M_{k}^{r}\left( s\right) -M_{k}^{l}\left( s\right) \right) \right\rangle ds
\end{equation*}%
converges to 0 for $a.e.$ $\omega $ as $k\rightarrow \infty $ thanks to
continuity of $t\rightarrow M\left( t\right) $. However, more is needed than
this. Define the stopping time
\begin{equation*}
\tau _{p}=\inf \left\{ t>0:\left\Vert M\left( t\right) \right\Vert
_{W}>p\right\} .
\end{equation*}%
Then $\tau _{p}\rightarrow \infty $ $\ a.e.$ $\omega .$ Let
\begin{equation*}
A_{k}=\left[ \left\vert \int_{t_{1}}^{t}\left\langle Y\left( s\right) ,\pi
_{n}\left( M_{k}^{r}\left( s\right) -M_{k}^{l}\left( s\right) \right)
\right\rangle ds\right\vert >\varepsilon \right]
\end{equation*}%
\begin{equation}
P\left( A_{k}\right) =\sum_{p=1}^{\infty }P\left( A_{k}\cap \left( \left[
\tau _{p}=\infty \right] \setminus \left[ \tau _{p-1}<\infty \right] \right)
\right)  \label{7june4h}
\end{equation}%
Now
\begin{equation*}
P\left( A_{k}\cap \left( \left[ \tau _{p}=\infty \right] \setminus \left[
\tau _{p-1}<\infty \right] \right) \right) \leq P\left( A_{k}\cap \left( %
\left[ \tau _{p}=\infty \right] \right) \right)
\end{equation*}%
\begin{equation*}
\leq P\left( \left[ \left\vert \int_{t_{1}}^{t}\left\langle Y\left( s\right)
,\pi _{n}\left( \left( M^{\tau _{p}}\right) _{k}^{r}\left( s\right) -\left(
M^{\tau _{p}}\right) _{k}^{l}\left( s\right) \right) \right\rangle
ds\right\vert >\varepsilon \right] \right)
\end{equation*}%
This is no larger than an expression of the form
\begin{equation}
\frac{C_{n}}{\varepsilon }\int_{\Omega }\int_{0}^{T}\left\Vert Y\left(
s\right) \right\Vert _{V^{\prime }}\left\Vert \left( M^{\tau _{p}}\right)
_{k}^{r}\left( s\right) -\left( M^{\tau _{p}}\right) _{k}^{l}\left( s\right)
\right\Vert _{W}dsdP  \label{7june3h}
\end{equation}%
The inside integral converges to 0 by continuity of $M$. Also, thanks to the
stopping time, the inside integral is dominated by an expression of the form
\begin{equation*}
\int_{0}^{T}\left\Vert Y\left( s\right) \right\Vert _{V^{\prime }}2pds
\end{equation*}%
and this is a function in $L^{1}\left( \Omega \right) $ by assumption on $Y$%
. It follows that the integral in \ref{7june3h} converges to 0 as $%
k\rightarrow \infty $ by the dominated convergence theorem. Hence
\begin{equation*}
\lim_{k\rightarrow \infty }P\left( A_{k}\cap \left( \left[ \tau _{p}=\infty %
\right] \right) \right) =0.
\end{equation*}%
Since the sets $\left[ \tau _{p}=\infty \right] \setminus \left[ \tau
_{p-1}<\infty \right] $ are disjoint, the sum of their probabilities is
finite. Hence by the dominated convergence theorem applied to the sum,
\begin{equation*}
\lim_{k\rightarrow \infty }P\left( A_{k}\right) =\sum_{p=0}^{\infty
}\lim_{k\rightarrow \infty }P\left( A_{k}\cap \left( \left[ \tau _{p}=\infty %
\right] \setminus \left[ \tau _{p-1}<\infty \right] \right) \right) =0
\end{equation*}%
Thus $\int_{t_{1}}^{t}\left\langle Y\left( s\right) ,\pi _{n}\left(
M_{k}^{r}\left( s\right) -M_{k}^{l}\left( s\right) \right) \right\rangle ds$
converges to 0 in probability as $k\rightarrow \infty $.

Now consider the other part of this expression,
\begin{equation*}
\int_{t_{1}}^{t}\left\langle Y\left( s\right) ,X_{k}^{r}\left( s\right)
-X_{k}^{l}\left( s\right) \right\rangle ds.
\end{equation*}%
This converges to 0 in $L^{1}\left( \Omega \right) $ because it is of the
form
\begin{equation*}
\int_{t_{1}}^{t}\left\langle Y\left( s\right) ,X_{k}^{r}\left( s\right)
\right\rangle ds-\int_{t_{1}}^{t}\left\langle Y\left( s\right)
,X_{k}^{l}\left( s\right) \right\rangle ds
\end{equation*}%
and both $X_{k}^{l}$ and $X_{k}^{r}$ converge to $X$ in $K$. Therefore, the
expression%
\begin{equation*}
\sum_{j=1}^{q_{k}-1}\left\langle \Delta P\left( t_{j}\right) -\Delta
BM\left( t_{j}\right) ,\Delta X\left( t_{j}\right) -\Delta M\left(
t_{j}\right) \right\rangle
\end{equation*}%
converges to 0 in probability. This establishes the desired formula for $%
t\in \mathcal{D}$.

To verify the last formula, let $t\in \mathcal{D}$. Then $t\in \mathcal{P}%
_{k}$ for some $k.$ Define
\begin{equation*}
\tau _{p}=\inf \left\{ t\in \mathcal{P}_{k}:\left\Vert P\left( t\right)
\right\Vert _{W^{\prime }}>p\right\}
\end{equation*}%
This is just the first hitting time of an adapted process so this is a well
defined stopping time. Then stop both sides of \ref{z29jane1g}. Thus
\begin{equation*}
\left\langle P\left( t\wedge \tau _{p}\right) ,X\left( t\wedge \tau
_{p}\right) \right\rangle =\left\langle BX_{0},X_{0}\right\rangle
+\int_{0}^{t\wedge \tau _{p}}\left( 2\left\langle Y\left( s\right) ,X\left(
s\right) \right\rangle +\left\langle BZ,Z\right\rangle _{\mathcal{L}%
_{2}}\right) ds
\end{equation*}%
\begin{equation*}
+2\int_{0}^{t}\mathcal{X}_{\left[ 0,\tau _{p}\right] }\left( Z\circ
J^{-1}\right) ^{\ast }P^{\tau _{p}}\circ JdW
\end{equation*}%
Now the last term is a martingale and you can take expectations of both
sides. Then
\begin{equation*}
E\left\langle P\left( t\wedge \tau _{p}\right) ,X\left( t\wedge \tau
_{p}\right) \right\rangle =E\left\langle BX_{0},X_{0}\right\rangle
+E\int_{0}^{t\wedge \tau _{p}}\left( 2\left\langle Y\left( s\right) ,X\left(
s\right) \right\rangle +\left\langle BZ,Z\right\rangle _{\mathcal{L}%
_{2}}\right) ds
\end{equation*}%
Then the integrands
\begin{equation*}
\left\langle P\left( t\wedge \tau _{p}\right) ,X\left( t\wedge \tau
_{p}\right) \right\rangle
\end{equation*}%
are uniformly integrable because
\begin{equation*}
\left\langle P\left( t\wedge \tau _{p}\right) ,X\left( t\wedge \tau
_{p}\right) \right\rangle \leq \sup_{t\in \mathcal{D}}\left\langle P\left(
t\right) ,X\left( t\right) \right\rangle
\end{equation*}%
which was shown to be in $L^{1}\left( \Omega \right) ,$ \ref{6aug3h}. Then
apply the Vitali convergence theorem to the left and the dominated
convergence theorem on the right to obtain the formula
\begin{equation*}
E\left\langle P\left( t\right) ,X\left( t\right) \right\rangle
=E\left\langle BX_{0},X_{0}\right\rangle +E\int_{0}^{t}\left( 2\left\langle
Y\left( s\right) ,X\left( s\right) \right\rangle +\left\langle
BZ,Z\right\rangle _{\mathcal{L}_{2}}\right) ds\ \blacksquare
\end{equation*}

Also we have the following improved version of Lemma \ref{27julyc1} in the
case that the integral equation holds for all $t$ off a set of measure zero.
See \cite{pre07} for a similar special case involving a Gelfand triple and $%
B=I$. That is, for $\omega $ off a set of measure zero,%
\begin{equation*}
BX\left( t\right) =BX_{0}+\int_{0}^{t}Y\left( s\right)
ds+B\int_{0}^{t}Z\left( s\right) dW\left( s\right)
\end{equation*}%
for all $t\in \left[ 0,T\right] $.

\begin{lemma}
\label{z24janl1gh} In the above situation where, off a set of measure zero, %
\ref{25oct1h}, the above integral equation holds for all $t\in \left[ 0,T%
\right] $, and $X$ is progressively measurable into $V$,
\begin{eqnarray*}
&&E\left( \sup_{t\in \left[ 0,T\right] }\left\langle BX\left( t\right)
,X\left( t\right) \right\rangle \right) \\
&<&C\left( \left\vert \left\vert Y\right\vert \right\vert _{K^{\prime
}},\left\vert \left\vert X\right\vert \right\vert _{K},\left\vert \left\vert
Z\right\vert \right\vert _{J},\left\Vert \left\langle
BX_{0},X_{0}\right\rangle \right\Vert _{L^{1}\left( \Omega \right) }\right)
<\infty .
\end{eqnarray*}%
where
\begin{eqnarray*}
J &=&L^{2}\left( \left[ 0,T\right] \times \Omega ;\mathcal{L}_{2}\left(
Q^{1/2}U;W\right) \right) ,K\equiv L^{p}\left( \left[ 0,T\right] \times
\Omega ;V\right) , \\
K^{\prime } &\equiv &L^{p^{\prime }}\left( \left[ 0,T\right] \times \Omega
;V^{\prime }\right) .
\end{eqnarray*}%
Also, $C$ is a continuous function of its arguments and $C\left(
0,0,0,0\right) =0$. Thus for a.e. $\omega ,$%
\begin{equation*}
\sup_{t\in \left[ 0,T\right] }\left\langle BX\left( t,\omega \right)
,X\left( t,\omega \right) \right\rangle \leq C\left( \omega \right) <\infty .
\end{equation*}%
For a.e. $\omega ,t\rightarrow BX\left( t,\omega \right) $ is weakly
continuous with values in $W^{\prime }$. \ Also $t\rightarrow \left\langle
BX\left( t\right) ,X\left( t\right) \right\rangle $ is lower semicontinuous.
\end{lemma}

\textbf{Proof of Lemma }\ref{z24janl1gh}: In the situation of this lemma, $%
P\left( t\right) =BX\left( t\right) $ for all $t$ provided $\omega $ is off
a single set of measure zero. Thus, there is a countable dense set $\mathcal{%
D}$ such that%
\begin{equation*}
E\left( \sup_{t\in \mathcal{D}}\left\langle BX\left( t\right) ,X\left(
t\right) \right\rangle \right) =E\left( \sup_{t\in \mathcal{D}%
}\sum_{k=1}^{\infty }\left\langle BX\left( t\right) ,e_{k}\right\rangle
^{2}\right)
\end{equation*}%
\begin{eqnarray}
&=&E\left( \sup_{t\in \mathcal{D}}\sum_{k=1}^{\infty }\left\langle P\left(
t\right) ,e_{k}\right\rangle ^{2}\right)  \notag \\
&<&C\left( \left\vert \left\vert Y\right\vert \right\vert _{K^{\prime
}},\left\vert \left\vert X\right\vert \right\vert _{K},\left\vert \left\vert
Z\right\vert \right\vert _{J},\left\Vert \left\langle
BX_{0},X_{0}\right\rangle \right\Vert _{L^{1}\left( \Omega \right) }\right)
<\infty .  \label{19septe1h}
\end{eqnarray}%
Now the function $t\rightarrow \sum_{k=1}^{\infty }\left\langle P\left(
t\right) ,e_{k}\right\rangle ^{2}$ is clearly lower semicontinuous. This is
because the partial sums are all continuous. Therefore, off the exceptional
set,
\begin{equation*}
\sup_{t\in \mathcal{D}}\sum_{k=1}^{\infty }\left\langle P\left( t\right)
,e_{k}\right\rangle ^{2}=\sup_{t\in \left[ 0,T\right] }\sum_{k=1}^{\infty
}\left\langle P\left( t\right) ,e_{k}\right\rangle ^{2}=\sup_{t\in \left[ 0,T%
\right] }\left\langle BX\left( t\right) ,X\left( t\right) \right\rangle \
\end{equation*}%
It follows that the desired estimates of Lemma \ref{z24janl1gh} are valid. $%
\blacksquare $

Proposition \ref{z29janl1g} along with the fundamental estimate of Lemma \ref%
{z24janl1gh} can be used to prove the following version of the Ito formula.
In proving this, we are considering the context of the integral equation \ref%
{25oct1h} holding for all $t$ provided $\omega $ is off a single set of
measure zero. The proof of this theorem follows the same methods used for a
similar result in \cite{pre07}.

\begin{theorem}
\label{z12febt1g}Suppose that off a set of measure zero, \ref{25oct1h} holds
for all $t$ so that $BX\left( t\right) =P\left( t\right) $. Then off a set
of measure zero, for every $t\in \left[ 0,T\right] ,$%
\begin{equation*}
\left\langle BX\left( t\right) ,X\left( t\right) \right\rangle =\left\langle
BX_{0},X_{0}\right\rangle +\int_{0}^{t}\left( 2\left\langle Y\left( s\right)
,X\left( s\right) \right\rangle +\left\langle BZ,Z\right\rangle _{\mathcal{L}%
_{2}}\right) ds
\end{equation*}%
\begin{equation}
+2\int_{0}^{t}\left( Z\circ J^{-1}\right) ^{\ast }BX\circ JdW
\label{z4febe1g}
\end{equation}%
Also%
\begin{equation*}
E\left( \left\langle BX\left( t\right) ,X\left( t\right) \right\rangle
\right) =
\end{equation*}%
\begin{equation}
E\left( \left\langle BX_{0},X_{0}\right\rangle \right) +E\left(
\int_{0}^{t}\left( 2\left\langle Y\left( s\right) ,X\left( s\right)
\right\rangle +\left\langle BZ,Z\right\rangle _{\mathcal{L}_{2}}\right)
ds\right)  \label{z12febe1g}
\end{equation}%
The quadratic variation of the stochastic integral is dominated by
\begin{equation}
C\int_{0}^{t}\left\Vert Z\right\Vert _{\mathcal{L}_{2}}^{2}\left\Vert
BX\right\Vert _{W^{\prime }}^{2}ds
\end{equation}%
for a suitable constant $C$. Also $t\rightarrow BX\left( t\right) $ is
continuous into $W^{\prime }$.
\end{theorem}

\textbf{Proof: }Let $t\notin \mathcal{D}.$ For $t>0,$ let $t\left( k\right) $
denote the largest point of $\mathcal{P}_{k}$ which is less than $t.$
Suppose $t\left( m\right) <t\left( k\right) $. Hence $m\leq k.$ Then
\begin{equation*}
P\left( t\left( m\right) \right) =BX_{0}+\int_{0}^{t\left( m\right) }Y\left(
s\right) ds+B\int_{0}^{t\left( m\right) }Z\left( s\right) dW\left( s\right) ,
\end{equation*}%
Thus for $t>t\left( m\right) ,$%
\begin{equation*}
P\left( t\right) -P\left( t\left( m\right) \right) =\int_{t\left( m\right)
}^{t}Y\left( s\right) ds+B\int_{t\left( m\right) }^{t}Z\left( s\right)
dW\left( s\right)
\end{equation*}%
which is the same sort of thing studied so far except that it starts at $%
t\left( m\right) $ rather than at $0$ and $BX_{0}=0.$ Therefore, from
Proposition \ref{z29janl1g} it follows
\begin{eqnarray*}
&&\left\langle P\left( t\left( k\right) \right) -P\left( t\left( m\right)
\right) ,X\left( t\left( k\right) \right) -X\left( t\left( m\right) \right)
\right\rangle \\
&=&\int_{t\left( m\right) }^{t\left( k\right) }\left( 2\left\langle Y\left(
s\right) ,X\left( s\right) -X\left( t\left( m\right) \right) \right\rangle
+\left\langle BZ,Z\right\rangle _{\mathcal{L}_{2}}\right) ds
\end{eqnarray*}%
\begin{equation}
+2\int_{t\left( m\right) }^{t\left( k\right) }\left( Z\circ J^{-1}\right)
^{\ast }\left( P\left( s\right) -P\left( t\left( m\right) \right) \right)
\circ JdW  \label{z31jane1g}
\end{equation}%
Consider that last term. It equals
\begin{equation}
2\int_{t\left( m\right) }^{t\left( k\right) }\left( Z\circ J^{-1}\right)
^{\ast }\left( P\left( s\right) -P_{m}^{l}\left( s\right) \right) \circ JdW
\label{z29jane2g}
\end{equation}%
This is dominated by
\begin{equation*}
2\left\vert \int_{0}^{t\left( k\right) }\left( Z\circ J^{-1}\right) ^{\ast
}\left( P\left( s\right) -P_{m}^{l}\left( s\right) \right) \circ JdW\right.
\end{equation*}%
\begin{equation*}
-\left. \int_{0}^{t\left( m\right) }\left( Z\circ J^{-1}\right) ^{\ast
}\left( P\left( s\right) -P_{m}^{l}\left( s\right) \right) \circ
JdW\right\vert
\end{equation*}%
\begin{equation*}
\leq 4\sup_{t\in \left[ 0,T\right] }\left\vert \int_{0}^{t}\left( Z\circ
J^{-1}\right) ^{\ast }\left( P\left( s\right) -P_{m}^{l}\left( s\right)
\right) \circ JdW\right\vert
\end{equation*}%
In Lemma \ref{z24janl9g} the above expression was shown to converge to 0 in
probability. Therefore, by the usual appeal to the Borel Canteli lemma,
there is a subsequence still referred to as $\left\{ m\right\} ,$ such that
the above expression converges to 0 pointwise in $\omega $ for all $\omega $
off some set of measure 0 as $m\rightarrow \infty $. It follows there is a
set of measure 0 such that for $\omega $ not in that set, \ref{z29jane2g}
converges to 0 in $\mathbb{R}$. Similar reasoning shows the first term in
the non stochastic integral of \ref{z31jane1g} is dominated by an expression
of the form
\begin{equation*}
4\int_{0}^{T}\left\vert \left\langle Y\left( s\right) ,X\left( s\right)
-X_{m}^{l}\left( s\right) \right\rangle \right\vert ds
\end{equation*}%
which clearly has a subsequence which converges to 0 for $\omega $ not in
some set of measure zero because $X_{m}^{l}$ converges in $K$ to $X$.
Finally, it is obvious that
\begin{equation*}
\lim_{m\rightarrow \infty }\int_{t\left( m\right) }^{t\left( k\right)
}\left\langle BZ,Z\right\rangle _{\mathcal{L}_{2}}ds=0\text{ for a.e. }\omega
\end{equation*}%
due to the assumptions on $Z$. For $\left\{ g_{i}\right\} $ an orthonormal
basis of $Q^{1/2}\left( U\right) ,$
\begin{eqnarray*}
\left\langle BZ,Z\right\rangle _{\mathcal{L}_{2}} &\equiv &\sum_{i}\left(
R^{-1}BZ\left( g_{i}\right) ,Z\left( g_{i}\right) \right)
=\sum_{i}\left\langle BZ\left( g_{i}\right) ,Z\left( g_{i}\right)
\right\rangle \\
&\leq &\left\Vert B\right\Vert \sum_{i}\left\Vert Z\left( g_{i}\right)
\right\Vert _{W}^{2}\in L^{1}\left( 0,T\right) \text{ a.e.}
\end{eqnarray*}

This shows that for $\omega $ off a set of measure 0%
\begin{equation}
\lim_{m,k\rightarrow \infty }\left\langle P\left( t\left( k\right) \right)
-P\left( t\left( m\right) \right) ,X\left( t\left( k\right) \right) -X\left(
t\left( m\right) \right) \right\rangle =0  \label{1aug1h}
\end{equation}%
Then for $x\in W,$%
\begin{eqnarray*}
&&\left\vert \left\langle P\left( t\left( k\right) \right) -P\left( t\left(
m\right) \right) ,x\right\rangle \right\vert \\
&\leq &\left\langle P\left( t\left( k\right) \right) -P\left( t\left(
m\right) \right) ,X\left( t\left( k\right) \right) -X\left( t\left( m\right)
\right) \right\rangle ^{1/2}\left\langle Bx,x\right\rangle ^{1/2} \\
&\leq &\left\langle P\left( t\left( k\right) \right) -P\left( t\left(
m\right) \right) ,X\left( t\left( k\right) \right) -X\left( t\left( m\right)
\right) \right\rangle ^{1/2}\left\Vert B\right\Vert ^{1/2}\left\Vert
x\right\Vert _{W}
\end{eqnarray*}%
and so
\begin{equation*}
\lim_{m,k\rightarrow \infty }\left\Vert P\left( t\left( k\right) \right)
-P\left( t\left( m\right) \right) \right\Vert _{W^{\prime }}=0
\end{equation*}%
Recall $t$ was arbitrary and $\left\{ t\left( k\right) \right\} $ is a
sequence converging to $t$. Then the above has shown that $\left\{ P\left(
t\left( k\right) \right) \right\} _{k=1}^{\infty }$ is a convergent sequence
in $W^{\prime }$. Does it converge to $P\left( t\right) ?$ Let $\xi \left(
t\right) \in W^{\prime }$ be what it converges to. Letting $v\in V$ then,
since $t\rightarrow P\left( t\right) $ is continuous into $V^{\prime },$
\begin{equation*}
\left\langle \xi \left( t\right) ,v\right\rangle =\lim_{k\rightarrow \infty
}\left\langle P\left( t\left( k\right) \right) ,v\right\rangle =\left\langle
P\left( t\right) ,v\right\rangle ,
\end{equation*}%
and now, since $V$ is dense in $W,$ this implies $\xi \left( t\right)
=P\left( t\right) $. Thus $P\left( t\right) =\lim_{k\rightarrow \infty
}P\left( t\left( k\right) \right) .$ Next consider the product $\left\langle
P\left( t\right) ,X\left( t\right) \right\rangle $.

For every $t\in \mathcal{D},$%
\begin{equation*}
\left\langle P\left( t\right) ,X\left( t\right) \right\rangle =\left\langle
BX_{0},X_{0}\right\rangle +\int_{0}^{t}\left( 2\left\langle Y\left( s\right)
,X\left( s\right) \right\rangle +\left\langle BZ,Z\right\rangle _{\mathcal{L}%
_{2}}ds\right) ds
\end{equation*}%
\begin{equation}
+2\int_{0}^{t}\left( Z\circ J^{-1}\right) ^{\ast }P\circ JdW
\label{z2june2h}
\end{equation}%
Does this formula hold for all $t\in \left[ 0,T\right] $?
\begin{equation*}
\left\vert \left\langle P\left( t\left( k\right) \right) ,X\left( t\left(
k\right) \right) \right\rangle -\left\langle P\left( t\right) ,X\left(
t\right) \right\rangle \right\vert
\end{equation*}%
\begin{eqnarray*}
&\leq &\left\vert \left\langle P\left( t\left( k\right) \right) ,X\left(
t\left( k\right) \right) \right\rangle -\left\langle P\left( t\right)
,X\left( t\left( k\right) \right) \right\rangle \right\vert \\
&&+\left\vert \left\langle P\left( t\right) ,X\left( t\left( k\right)
\right) \right\rangle -\left\langle P\left( t\right) ,X\left( t\right)
\right\rangle \right\vert
\end{eqnarray*}%
\begin{equation*}
=\left\vert \left\langle P\left( t\left( k\right) \right) -P\left( t\right)
,X\left( t\left( k\right) \right) \right\rangle \right\vert +\left\vert
\left\langle P\left( t\right) ,X\left( t\left( k\right) \right) -X\left(
t\right) \right\rangle \right\vert
\end{equation*}%
Since $BX\left( t\right) =P\left( t\right) ,$ the Cauchy Schwarz inequality
implies that the above expression is dominated by
\begin{eqnarray*}
&\leq &\left\langle P\left( t\left( k\right) \right) -P\left( t\right)
,X\left( t\left( k\right) \right) -X\left( t\right) \right\rangle ^{1/2} \\
&&\cdot \left( \left\langle BX\left( t\left( k\right) \right) ,X\left(
t\left( k\right) \right) \right\rangle ^{1/2}+\left\langle BX\left( t\right)
,X\left( t\right) \right\rangle ^{1/2}\right)
\end{eqnarray*}%
Also$,$
\begin{equation*}
\left\langle P\left( t\left( k\right) \right) -P\left( t\right) ,X\left(
t\left( k\right) \right) -X\left( t\right) \right\rangle =\left\langle
BX\left( t\right) -P\left( t\left( k\right) \right) ,X\left( t\right)
-X\left( t\left( k\right) \right) \right\rangle
\end{equation*}%
\begin{equation*}
=\left\langle BX\left( t\right) ,X\left( t\right) \right\rangle
-\left\langle P\left( t\right) ,X\left( t\left( k\right) \right)
\right\rangle -\left\langle P\left( t\left( k\right) \right) ,X\left(
t\right) \right\rangle +\left\langle P\left( t\left( k\right) \right)
,X\left( t\left( k\right) \right) \right\rangle
\end{equation*}%
The expression above simplifies to
\begin{equation*}
\left\langle BX\left( t\right) ,X\left( t\right) \right\rangle
-2\left\langle P\left( t\right) ,X\left( t\left( k\right) \right)
\right\rangle +\left\langle P\left( t\left( k\right) \right) ,X\left(
t\left( k\right) \right) \right\rangle
\end{equation*}%
which is clearly lower semicontinuous in $t$ due to the continuity of $%
P\left( t\right) $ into $V^{\prime }$ and the equation
\begin{equation*}
\left\langle BX\left( t\right) ,X\left( t\right) \right\rangle
=\sum_{k=1}^{\infty }\left\langle BX\left( t\right) ,e_{k}\right\rangle
^{2}=\sum_{k=1}^{\infty }\left\langle P\left( t\right) ,e_{k}\right\rangle
^{2}
\end{equation*}%
Summarizing the above, this has shown that
\begin{equation*}
\left\vert \left\langle P\left( t\left( k\right) \right) ,X\left( t\left(
k\right) \right) \right\rangle -\left\langle P\left( t\right) ,X\left(
t\right) \right\rangle \right\vert \leq \left\langle P\left( t\left(
k\right) \right) -P\left( t\right) ,X\left( t\left( k\right) \right)
-X\left( t\right) \right\rangle ^{1/2}\cdot
\end{equation*}%
\begin{equation}
\left( \left\langle BX\left( t\left( k\right) \right) ,X\left( t\left(
k\right) \right) \right\rangle ^{1/2}+\left\langle BX\left( t\right)
,X\left( t\right) \right\rangle ^{1/2}\right)  \label{23july1}
\end{equation}%
and also that $t\rightarrow \left\langle P\left( t\left( k\right) \right)
-P\left( t\right) ,X\left( t\left( k\right) \right) -X\left( t\right)
\right\rangle ^{1/2}$ is lower semicontinuous.

Consider the right side of the above.
\begin{equation*}
t\rightarrow \left\langle P\left( t\left( k\right) \right) -BX\left(
t\right) ,X\left( t\left( k\right) \right) -X\left( t\right) \right\rangle
\end{equation*}%
Since $\left\langle BX\left( t\right) ,X\left( t\right) \right\rangle =$ $%
\left\langle P\left( t\right) ,X\left( t\right) \right\rangle $ is bounded, it follows that
\begin{eqnarray*}
&&\left\vert \left\langle P\left( t\left( k\right) \right) ,X\left( t\left(
k\right) \right) \right\rangle -\left\langle P\left( t\right) ,X\left(
t\right) \right\rangle \right\vert \\
&\leq &C\left\langle P\left( t\left( k\right) \right) -P\left( t\right)
,X\left( t\left( k\right) \right) -X\left( t\right) \right\rangle ^{1/2}
\end{eqnarray*}%
From the above, the right side equals a lower semicontinuous function.
Therefore, passing to a limit and using the lower semicontinuity,
\begin{equation*}
\left\vert \left\langle P\left( t\left( k\right) \right) ,X\left( t\left(
k\right) \right) \right\rangle -\left\langle P\left( t\right) ,X\left(
t\right) \right\rangle \right\vert \leq C\left\langle P\left( t\left(
k\right) \right) -P\left( t\right) ,X\left( t\left( k\right) \right)
-X\left( t\right) \right\rangle ^{1/2}
\end{equation*}%
\begin{equation*}
\leq C\lim \inf_{m\rightarrow \infty }\left\langle P\left( t\left( k\right)
\right) -P\left( t\left( m\right) \right) ,X\left( t\left( k\right) \right)
-X\left( t\left( m\right) \right) \right\rangle ^{1/2}<\varepsilon
\end{equation*}%
provided $k$ is sufficiently large (by \ref{1aug1h}). Since $\varepsilon $
is arbitrary,
\begin{equation*}
\lim_{k\rightarrow \infty }\left\langle P\left( t\left( k\right) \right)
,X\left( t\left( k\right) \right) \right\rangle =\left\langle BX\left(
t\right) ,X\left( t\right) \right\rangle .
\end{equation*}%
It follows that for $\omega $ off the set of measure zero $N,$ the formula %
\ref{z2june2h} is valid for all $t.$ Now this formula shows that off a set
of measure zero, $t\rightarrow \left\langle P\left( t\right) ,X\left(
t\right) \right\rangle $ is continuous.

This implies that $t\rightarrow P\left( t\right) =BX\left( t\right) $ is
continuous with values in $W^{\prime }$. Here is why. The fact that the
formula \ref{z2june2h} holds for all $t$ implies that $t\rightarrow
\left\langle BX\left( t\right) ,X\left( t\right) \right\rangle $ is
continuous. Then for $x\in W,$
\begin{equation}
\left\vert \left\langle BX\left( t\right) -BX\left( s\right) ,x\right\rangle
\right\vert \leq \left\langle B\left( X\left( t\right) -X\left( s\right)
\right) ,X\left( t\right) -X\left( s\right) \right\rangle ^{1/2}\left\Vert
B\right\Vert ^{1/2}\left\Vert x\right\Vert _{W}.  \label{z2june3h}
\end{equation}%
Also
\begin{eqnarray*}
&&\left\langle B\left( X\left( t\right) -X\left( s\right) \right) ,X\left(
t\right) -X\left( s\right) \right\rangle \\
&=&\left\langle BX\left( t\right) ,X\left( t\right) \right\rangle
+\left\langle BX\left( s\right) ,X\left( s\right) \right\rangle
-2\left\langle BX\left( t\right) ,X\left( s\right) \right\rangle
\end{eqnarray*}%
By weak continuity of $t\rightarrow BX\left( t\right) $ shown earlier,
\begin{equation*}
\lim_{t\rightarrow s}\left\langle BX\left( t\right) ,X\left( s\right)
\right\rangle =\left\langle BX\left( s\right) ,X\left( s\right)
\right\rangle .
\end{equation*}%
Therefore,
\begin{equation*}
\lim_{t\rightarrow s}\left\langle B\left( X\left( t\right) -X\left( s\right)
\right) ,X\left( t\right) -X\left( s\right) \right\rangle =0
\end{equation*}%
and so the inequality \ref{z2june3h} implies the continuity of $t\rightarrow
BX\left( t\right) $ into $W^{\prime }$.

Now consider the claim about the expectation. Since the stochastic integral%
\begin{equation*}
2\int_{0}^{t}\left( Z\circ J^{-1}\right) ^{\ast }P\circ JdW
\end{equation*}%
is only a local martingale, it is necessary to employ a stopping time. Since
$t\rightarrow \left\langle BX\left( t\right) ,X\left( t\right) \right\rangle
$ is continuous, one can define a stopping time
\begin{equation*}
\tau _{p}\equiv \inf \left\{ t>0:\left\langle BX\left( t\right) ,X\left(
t\right) \right\rangle >p\right\}
\end{equation*}%
Then use the stopping time in both sides of \ref{z4febe1g} and take the
expectation. The stopped local martingale has expectation equal to 0. Thus
\begin{equation*}
E\left( \left\langle BX^{\tau _{p}}\left( t\right) ,X^{\tau _{p}}\left(
t\right) \right\rangle \right) =E\left( \left\langle
BX_{0},X_{0}\right\rangle \right)
\end{equation*}%
\begin{equation*}
+E\left( \int_{0}^{t}\mathcal{X}_{\left[ 0,\tau _{p}\right] }\left(
2\left\langle Y\left( s\right) ,X\left( s\right) \right\rangle +\left\langle
BZ,Z\right\rangle _{\mathcal{L}_{2}}\right) ds\right)
\end{equation*}%
Next use the dominated convergence theorem on the right and the monotone
convergence theorem on the left to let $p\rightarrow \infty $ and obtain the
desired result. The claim about the quadratic variation follows from the
description of the quadratic variation for a stochastic integral.

Another interesting observation is that $t\rightarrow BX\left( t\right) $ is
continuous into $W^{\prime }.$%
\begin{equation*}
\left\langle BX\left( t\right) -BX\left( s\right) ,X\left( t\right) -X\left(
s\right) \right\rangle =\left\langle BX\left( t\right) ,X\left( t\right)
\right\rangle +\left\langle BX\left( s\right) ,X\left( s\right)
\right\rangle -2\left\langle BX\left( t\right) ,X\left( s\right)
\right\rangle
\end{equation*}%
From the above formula, it is known that $t\rightarrow \left\langle BX\left(
t\right) ,X\left( t\right) \right\rangle $ is continuous. It was also shown
above that $t\rightarrow BX\left( t\right) $ is weakly continuous into $%
W^{\prime }$. Therefore, you could let $t\rightarrow s$ and conclude that
\begin{equation*}
\lim_{t\rightarrow s}\left\langle BX\left( t\right) -BX\left( s\right)
,X\left( t\right) -X\left( s\right) \right\rangle =0
\end{equation*}%
It follows that for $w\in W,$
\begin{eqnarray*}
\left\langle BX\left( t\right) -BX\left( s\right) ,w\right\rangle  &\leq
&\left\langle BX\left( t\right) -BX\left( s\right) ,X\left( t\right)
-X\left( s\right) \right\rangle ^{1/2}\left\langle Bw,w\right\rangle ^{1/2}
\\
\, &\leq &\left\langle BX\left( t\right) -BX\left( s\right) ,X\left(
t\right) -X\left( s\right) \right\rangle ^{1/2}\left\Vert B\right\Vert
^{1/2}\left\Vert w\right\Vert
\end{eqnarray*}%
and so
\begin{equation*}
\left\Vert BX\left( t\right) -BX\left( s\right) \right\Vert _{W^{\prime
}}\leq \left\langle BX\left( t\right) -BX\left( s\right) ,X\left( t\right)
-X\left( s\right) \right\rangle ^{1/2}\left\Vert B\right\Vert ^{1/2}
\end{equation*}%
which converges to 0 as $t\rightarrow s$. $\blacksquare $

\section{An application to evolution equations}

\label{evolutionequation}

First we consider the case of a stochastic equation in a single Hilbert
space. Here we give an example of how the Ito formula can be used to obtain
theorems of existence and uniqueness. This begins with an introductory
result on evolution equations in a single Hilbert space which is included
for the sake of completeness. In what follows, $H$ is a separable Hilbert
space. It will be assumed that for each $t,\omega ,$%
\begin{equation*}
u\rightarrow A\left( t,u,\omega \right)
\end{equation*}%
is a mapping from $H$ to $H$. Assume also that
\begin{equation*}
\left( t,u,\omega \right) \rightarrow A\left( t,u,\omega \right)
\end{equation*}%
is progressively measurable.

It is possible to assume only that $u\rightarrow A\left( t,u,\omega \right) $
is continuous and base the theory on this. It is more troublesome because
you end up having to consider finite dimensional subspaces and it would
distract attention from the issue of interest in this paper. Therefore, it
will be assumed here that for each $B\left( \mathbf{0,}r\right) $ the
restriction of $A\left( t,\cdot ,\omega \right) $ to $B\left( \mathbf{0,}%
r\right) $ is Lipschitz continuous. Thus
\begin{equation}
\left\vert A\left( t,u,\omega \right) -A\left( t,v,\omega \right)
\right\vert \leq K_{r}\left\vert u-v\right\vert  \label{20jane5h}
\end{equation}%
whenever $u,v\in \overline{B\left( \mathbf{0,}r\right) }.$ Also assume in
addition to the above Lipschitz condition, the estimate,
\begin{equation}
\ \left\langle A\left( t,u\right) ,u\right\rangle \geq -k\left\vert
u\right\vert ^{2}-C.  \label{20jane6h}
\end{equation}%
where $C\in L^{1}\left( \left[ 0,T\right] \times \Omega \right) ,C\geq 0$.
It is also assumed that
\begin{equation*}
\Phi \in L^{2}\left( \left[ 0,T\right] \times \Omega ,\mathcal{L}_{2}\left(
Q^{1/2}U,H\right) \right)
\end{equation*}%
is given. It is routine to generalize this to the case where $\Phi $ depends
on the unkown function $u$.

Then under these conditions, we can prove the following theorem.

\begin{theorem}
\label{10jant1h} Let $u\rightarrow A\left( t,u,\omega \right) $ be locally
Lipschitz in the sense that for each $B\left( \mathbf{0},r\right) ,$ $A$
restricted to $B\left( \mathbf{0},r\right) $ is Lipschitz. \ Also suppose $%
\left( t,u,\omega \right) \rightarrow A\left( t,u,\omega \right) $ is
progressively measurable. Then if $u_{0}\in L^{2}\left( \Omega ,H\right) $
with $u_{0}$ measurable in $\mathcal{F}_{0},$
\begin{equation*}
\left\langle A\left( t,u\right) ,u\right\rangle \geq -k\left\vert
u\right\vert ^{2}-C,\text{ }C\geq 0,
\end{equation*}%
where $C\in L^{1}\left( \left[ 0,T\right] \times \Omega \right) ,$ it
follows that there exists a progressively measurable function $u$ and a set
of measure zero $N$, such that for $\omega \notin N,$%
\begin{equation}
u\left( t\right) -u_{0}+\int_{0}^{t}A\left( u\right) ds=\int_{0}^{t}\Phi dW.
\label{24jane1h}
\end{equation}
\end{theorem}

\vspace{1pt}\textbf{Proof: } Let $P_{n}$ denote the projection onto $%
\overline{B\left( \mathbf{0,}9^{n}\right) }$ and let $u_{n}$ be the solution
to
\begin{equation*}
u_{n}\left( t\right) -u_{0}+\int_{0}^{t}A\left( P_{n}u_{n}\right)
ds=\int_{0}^{t}\Phi dW
\end{equation*}%
That a unique progressively measurable solution exists follows readily from
showing that a high enough power of an operator is a contraction map, just
as in the deterministic case. The solution holds for all $t\in \left[ 0,T%
\right] $ for $\omega $ off a set of measure zero.

Next let
\begin{equation*}
\tau _{n}\equiv \inf \left\{ t>0:\left\vert u_{n}\left( t\right) \right\vert
^{2}>2^{n}\right\}
\end{equation*}%
Thus from localization as described in \cite{pre07} and \cite{dap92},
\begin{equation*}
u_{n}^{\tau _{n}}\left( t\right) -u_{0}+\int_{0}^{t}\mathcal{X}_{\left[
0,\tau _{n}\right] }A\left( u_{n}^{\tau _{n}}\right) ds=\int_{0}^{t}\mathcal{%
X}_{\left[ 0,\tau _{n}\right] }\Phi dW
\end{equation*}%
It is important to get an estimate now. From the standard Ito formula or
Theorem \ref{z12febt1g}, letting $F\left( u\right) =\left\vert u\right\vert
_{H}^{2},$ then using the boundedness of $u_{n}^{\tau _{n}},$
\begin{equation*}
\frac{1}{2}\left\vert u_{n}^{\tau _{n}}\left( t\right) \right\vert ^{2}-%
\frac{1}{2}\left\vert u_{0}\right\vert ^{2}+\int_{0}^{t}\mathcal{X}_{\left[
0,\tau _{n}\right] }\left\langle A\left( u_{n}^{\tau _{n}}\right)
,u_{n}^{\tau _{n}}\right\rangle ds=\int_{0}^{t}\mathcal{X}_{\left[ 0,\tau
_{n}\right] }\left\Vert \Phi \right\Vert ^{2}ds+M\left( t\right)
\end{equation*}%
where $M\left( t\right) $ is a local martingale with
\begin{equation*}
\left[ M\right] \left( T\right) \leq C\int_{0}^{T}\left\Vert \Phi
\right\Vert _{\mathcal{L}_{2}}^{2}\left\vert u_{n}^{\tau _{n}}\right\vert
_{H}^{2}dt
\end{equation*}%
Then from maximal estimates and Burkholder Davis Gundy inequality,
\begin{equation*}
P\left( \left[ \sup_{t\in \left[ 0,T\right] }\left\vert
\begin{array}{c}
\frac{1}{2}\left\vert u_{n}^{\tau _{n}}\left( t\right) \right\vert ^{2}-%
\frac{1}{2}\left\vert u_{0}\right\vert ^{2}+ \\
\int_{0}^{t}\mathcal{X}_{\left[ 0,\tau _{n}\right] }\left\langle A\left(
u_{n}^{\tau _{n}}\right) ,u_{n}^{\tau _{n}}\right\rangle ds-\int_{0}^{t}%
\mathcal{X}_{\left[ 0,\tau _{n}\right] }\left\Vert \Phi \right\Vert ^{2}ds%
\end{array}%
\right\vert >\lambda \right] \right)
\end{equation*}%
\begin{equation*}
\leq \frac{1}{\lambda }\int_{\Omega }\sup \left\{ \left\vert M\left(
t\right) \right\vert ,t\in \left[ 0,T\right] \right\} dP\leq C\frac{1}{%
\lambda }\int_{\Omega }\left[ M\right] \left( T\right) ^{1/2}dP
\end{equation*}%
Now from the description of the quadratic variation for stochastic integrals
and using the stopping time,%
\begin{equation*}
\leq \frac{C}{\lambda }\int_{\Omega }\left( \int_{0}^{T}\left\Vert \Phi
\right\Vert ^{2}2^{n}\right) ^{1/2}dP=\frac{C\left( \Phi \right) 2^{n/2}}{%
\lambda }
\end{equation*}%
The above holds for each $n.$ Let $\lambda =\left( \frac{3}{2}\right) ^{n}$.
Then the above implies
\begin{equation*}
P\left( \left[ \sup_{t\in \left[ 0,T\right] }\left\vert
\begin{array}{c}
\frac{1}{2}\left\vert u_{n}^{\tau _{n}}\left( t\right) \right\vert ^{2}-%
\frac{1}{2}\left\vert u_{0}\right\vert ^{2}-k\int_{0}^{t}\left\vert
u_{n}^{\tau _{n}}\left( s\right) \right\vert ^{2} \\
-\int_{0}^{t}\left\Vert \Phi \right\Vert ^{2}ds-\int_{0}^{t}Cds%
\end{array}%
\right\vert >\left( \frac{3}{2}\right) ^{n}\right] \right)
\end{equation*}%
\begin{equation*}
\leq C\left( \Phi \right) \frac{2^{n/2}}{\left( 3/2\right) ^{n}}\leq C\left(
\Phi \right) \left( .96\right) ^{n}
\end{equation*}

By the Borel Cantelli lemma, it follows that there exists a set of measure
zero $N$ such that for $\omega \notin N,$ all $n$ large enough, say $n\geq
M\left( \omega \right) $ and $t,$
\begin{equation*}
\frac{1}{2}\left\vert u_{n}^{\tau _{n}}\left( t\right) \right\vert ^{2}-%
\frac{1}{2}\left\vert u_{0}\right\vert ^{2}-k\int_{0}^{t}\left\vert
u_{n}^{\tau _{n}}\left( s\right) \right\vert ^{2}ds-\int_{0}^{t}\left\Vert
\Phi \right\Vert ^{2}ds-\int_{0}^{t}Cds\leq \left( \frac{3}{2}\right) ^{n}.
\end{equation*}%
for such $\omega $ and $n,$
\begin{equation*}
\left\vert u_{n}^{\tau _{n}}\left( t\right) \right\vert ^{2}\leq 2\left(
\left( \frac{3}{2}\right) ^{n}+\left\vert u_{0}\right\vert
^{2}+2\int_{0}^{T}\left\Vert \Phi \right\Vert ^{2}ds+2\int_{0}^{T}Cds\right)
+2k\int_{0}^{t}\left\vert u_{n}^{\tau _{n}}\left( s\right) \right\vert ^{2}ds
\end{equation*}%
Apply Gronwall's inequality to conclude that for all $t\in \left[ 0,T\right]
,$%
\begin{equation*}
\left\vert u_{n}^{\tau _{n}}\left( t\right) \right\vert ^{2}\leq 2\left(
\left( \frac{3}{2}\right) ^{n}+\left\vert u_{0}\right\vert
^{2}+2\int_{0}^{T}\left\Vert \Phi \right\Vert ^{2}ds+2\int_{0}^{T}Cds\right)
e^{2kT}
\end{equation*}%
If $n$ is sufficiently large, the right side is smaller than $2^{n}$ for all
$t\in \left[ 0,T\right] $. Thus for such $\omega $ and $n,$%
\begin{equation*}
\left\vert u_{n}\left( t\wedge \tau _{n}\right) \right\vert ^{2}<2^{n}
\end{equation*}%
It follows that $t<\tau _{n}$ for all $t\in \left[ 0,T\right] $ because if
not, then you would have
\begin{equation*}
2^{n}=\left\vert u_{n}\left( \tau _{n}\right) \right\vert ^{2}<2^{n}.
\end{equation*}%
Hence for all $n$ large enough, $\tau _{n}=\infty .$ For $n$ and $\omega $
as just described,
\begin{equation}
u_{n}\left( t\right) -u_{0}+\int_{0}^{t}A\left( P_{n}u_{n}\right)
ds=\int_{0}^{t}\Phi dW.  \label{10jane1h}
\end{equation}

Of course the problem here is that $n$ depends on $\omega $ and we\ need a
single function $u$. Suppose then that $\omega \notin N$ and both $m,n$ are
so large that $\tau _{m}\left( \omega \right) =\tau _{n}\left( \omega
\right) =\infty $. Say $n>m$. Then
\begin{eqnarray*}
u_{n}\left( t\right) -u_{0}+\int_{0}^{t}A\left( P_{n}u_{n}\right) ds
&=&\int_{0}^{t}\Phi dW \\
u_{m}\left( t\right) -u_{0}+\int_{0}^{t}A\left( P_{m}u_{m}\right) ds
&=&\int_{0}^{t}\Phi dW
\end{eqnarray*}%
However, $\left\vert u_{n}\left( t\right) \right\vert ^{2}<3^{n},\left\vert
u_{m}\left( t\right) \right\vert ^{2}<3^{m}$ and so the $P_{n}u_{n}$ and $%
P_{m}u_{m}$ equal $u_{n}$ and $u_{m}$ respectively. Hence
\begin{equation*}
u_{n}\left( t\right) -u_{m}\left( t\right) +\int_{0}^{t}A\left( u_{n}\right)
-A\left( u_{m}\right) ds=0.
\end{equation*}%
Furthermore, all the values of these two functions are in $\overline{B\left(
\mathbf{0,}3^{n}\right) }$. Note that this is a deterministic integral, not
one of those stochastic integrals. Therefore, by the local Lipschitz
assumption, there exists a $K$ such that
\begin{equation*}
\left\vert u_{n}\left( t\right) -u_{m}\left( t\right) \right\vert \leq
\int_{0}^{t}\left\vert A\left( u_{n}\right) -A\left( u_{m}\right)
\right\vert ds\leq K\int_{0}^{t}\left\vert u_{n}-u_{m}\right\vert ds
\end{equation*}%
and so, by Gronwall's inequality, $u_{n}\left( t\right) =u_{m}\left(
t\right) $ for that $\omega $.

Because of this, define for $\omega \notin N,$%
\begin{equation*}
u\left( t\right) \equiv u_{n}\left( t\right) \text{ where }\tau _{n}\left(
\omega \right) =\infty \text{.}
\end{equation*}%
It was just shown that this is well defined. Also from \ref{10jane1h}, it
follows that
\begin{equation*}
u\left( t\right) -u_{0}+\int_{0}^{t}A\left( u\right) ds=\int_{0}^{t}\Phi dW.
\end{equation*}%
This proves existence.

It only remains to verify uniqueness. If $v$ is another such solution, then
taking the union of the two exceptional sets, it follows that for $\omega $
not in this union,
\begin{equation*}
u\left( t\right) -v\left( t\right) +\int_{0}^{t}\left( Au-Av\right) ds=0
\end{equation*}%
Thus, since both $u$ and $v$ are bounded, there is a Lipschitz constant $K$
such that
\begin{equation*}
\left\vert u\left( t\right) -v\left( t\right) \right\vert \leq
\int_{0}^{t}\left\vert Au-Av\right\vert ds\leq K\int_{0}^{t}\left\vert
u\left( t\right) -v\left( t\right) \right\vert dt
\end{equation*}%
and so, by Gronwall's inequality, $u\left( t\right) =v\left( t\right) $. $%
\blacksquare $

Note that there is no monotonicity required on $A$ in order to obtain
existence.

\section{Multiple Spaces}

Next we consider the case of variational evolution equations in infinite
dimensional spaces. Consider the case of a reflexive separable Banach space $%
V$ and a Hilbert space $W$ such that $V\subseteq W$ with $V$ dense in $W.$
Thus $W^{\prime }\subseteq V^{\prime }$. Suppose there exists a Hilbert
space $E$ which is dense in $V$. Thus
\begin{equation}
E\subseteq V\subseteq W,\ W^{\prime }\subseteq V^{\prime }\subseteq
E^{\prime }  \label{30octe2h}
\end{equation}%
and let $R:E\rightarrow E^{\prime }$ be the Riesz map. Let $\left(
t,u,\omega \right) \rightarrow A\left( t,u,\omega \right) $ where $A\left(
t,u,\omega \right) \in V^{\prime }.$ Suppose
\begin{equation}
\left( t,u,\omega \right) \rightarrow A\left( t,u,\omega \right)
\label{4may1f}
\end{equation}%
is progressively measurable. Also assume the coercivity condition
\begin{equation}
\left\langle A\left( t,u,\omega \right) ,u\right\rangle \geq k\left\Vert
u\right\Vert _{V}^{p}-C\left( t,\omega \right)  \label{4may2f}
\end{equation}%
where $C\in L^{1}\left( \left[ 0,T\right] \times \Omega \right) $. Let $%
\mathcal{V}\equiv L^{p}\left( \left[ 0,T\right] \times \Omega ,V\right) $
and let $\mathcal{V}^{\prime }$ be its dual space $L^{p^{\prime }}\left( %
\left[ 0,T\right] \times \Omega ,V^{\prime }\right) $. We assume the
operator $A:\mathcal{V\rightarrow V}^{\prime }$ is type $M$ \cite{Lio69}.

This is a more general condition than monotone and hemicontinuous. However,
the question whether there exist meaningful examples which are type $M$ on $%
\mathcal{V}$ which are not also monotone and hemicontinuous is being left
open for now. We have no such examples. However, the type $M$ condition is
convenient to use and so this is why we make this theoretically more general
assumption.

Also let $B:W\rightarrow W^{\prime }$ be nonnegative and self adjoint. In
all of the above, the $\sigma $ algebra will be the product measurable sets $%
\mathcal{B}\left( \left[ 0,T\right] \right) \times \mathcal{F}_{T}$. Then we
need some sort of continuity condition on $u\rightarrow A\left( t,u,\omega
\right) $. In general, is suffices to assume this map is demicontinuous,
possibly even less. However, here we will assume more for the sake of
convenience. Assume
\begin{equation}
u\rightarrow A\left( t,u,\omega \right) \text{ is locally Lipschitz}
\label{4may4.5f}
\end{equation}%
as a map from $E$ to $E^{\prime }$. We note that this condition is often
true in many applications of interest thanks to the Sobolev embedding
theorem. One takes $E$ to be a suitable closed subspace of $H^{k}\left(
U\right) $ for $k$ sufficiently large.

Also let $\Phi \in L^{2}\left( \left( \left[ 0,T\right] \times \Omega
\right) ,\mathcal{L}_{2}\left( JQ^{1/2}U,W\right) \right) $ so we can
consider $\int_{0}^{t}\Phi dW,$ and it has values in the space $W$.

Then the main result to be proved is Theorem \ref{7mayt1f} and its
corollaries stated below. They give existence for a solution to the integral
equation
\begin{equation*}
Bu\left( t\right) -Bu_{0}+\int_{0}^{t}A\left( u\right) ds=\int_{0}^{t}B\Phi
dW
\end{equation*}%
in the sense that for a.e. $\omega ,$ the equation holds for all $t\in \left[
0,T\right] .$

This theorem is proved by using Theorem \ref{10jant1h} to obtain existence
for a regularized problem. The Ito formula is then used to obtain estimates
on these solutions. After this, weakly convergent subsequences are obtained
which are then shown to converge to the desired solution through the use of
the Ito formula presented above, along with the assumption that $A$ is type $%
M$.

\begin{lemma}
\label{6augl2h} Let $u_{0}\in L^{q}\left( \Omega ,E\right) $ where $q=\max
\left( p,2\right) $. Also let $R$ be the Riesz map from $E$ to $E^{\prime }.$
Then there exists a solution to the integral equation
\begin{equation*}
u\left( t\right) -u_{0}+\int_{0}^{t}\left( B+\varepsilon R\right)
^{-1}\left( A\left( u\right) +\varepsilon R\left( u\right) \right)
ds=\int_{0}^{t}\left( B+\varepsilon R\right) ^{-1}B\Phi dW+\left(
B+\varepsilon R\right) ^{-1}\int_{0}^{t}fds
\end{equation*}%
in the sense that off a set of measure zero, the equation holds for all $t$.
This solution satisfies the estimate%
\begin{equation*}
\frac{1}{2}E\left\langle \left( B+\varepsilon R\right) u\left( t\right)
,u\left( t\right) \right\rangle -\frac{1}{2}E\left\langle \left(
B+\varepsilon R\right) u_{0},u_{0}\right\rangle +E\int_{0}^{t}\left\langle
Au,u\right\rangle +\varepsilon \left\langle Ru,u\right\rangle ds
\end{equation*}%
\begin{equation*}
\leq \frac{1}{2}E\int_{0}^{t}\left( \mathcal{R}^{-1}B\Phi ,\Phi \right) _{%
\mathcal{L}_{2}\left( Q^{1/2}U,W\right) }ds+E\int_{0}^{t}\left\langle
f,u\right\rangle ds
\end{equation*}%
where $\mathcal{R}$ is the Riez map from $W$ to $W^{\prime }$.
\end{lemma}

\textbf{Proof: }Let $R$ be the Riesz map from $E$ to $E^{\prime }$. Then
there exists an equivalent Hilbert space norm on $E$ such that for fixed $%
\varepsilon >0,$ the Riesz map is $B+\varepsilon R$. To simplify the
notation, let
\begin{equation*}
A_{\varepsilon }\left( u\right) \equiv A\left( u\right) +\varepsilon R\left(
u\right)
\end{equation*}%
By Theorem \ref{10jant1h}, for $u_{0}\in L^{2}\left( \Omega ,E\right) $ with
$u_{0}$ an $\mathcal{F}_{0}$ measurable function, there exists a unique
progressively measurable function $u$ having values in $E$ such that
\begin{equation}
u\left( t\right) -u_{0}+\int_{0}^{t}\left( B+\varepsilon R\right)
^{-1}A_{\varepsilon }\left( u\right) ds=\int_{0}^{t}\left( B+\varepsilon
R\right) ^{-1}B\Phi dW+\left( B+\varepsilon R\right) ^{-1}\int_{0}^{t}fds
\label{4may5f}
\end{equation}%
This is because the integrand is locally Lipschitz and it satisfies
\begin{equation*}
\left( \left( B+\varepsilon R\right) ^{-1}A_{\varepsilon }\left( t,u,\omega
\right) ,u\right) _{E}=\left\langle A_{\varepsilon }\left( t,u,\omega
\right) ,u\right\rangle _{V^{\prime },V}\geq C_{\varepsilon }\left\Vert
u\right\Vert _{E}^{2}+k\left\Vert u\right\Vert _{V}^{p}-C\left( t,\omega
\right)
\end{equation*}%
Multiplying through by $\left( B+\varepsilon R\right) ,$ this shows that
there exists a unique progressively measurable $u$ which is the solution to
\begin{equation}
\left( B+\varepsilon R\right) u\left( t\right) -\left( B+\varepsilon
R\right) u_{0}+\int_{0}^{t}A_{\varepsilon }\left( u\right)
ds=\int_{0}^{t}B\Phi dW+\int_{0}^{t}fds  \label{4may6f}
\end{equation}%
That stochastic integral on the right equals
\begin{equation*}
\left( B+\varepsilon R\right) \int_{0}^{t}\left( B+\varepsilon R\right)
^{-1}B\Phi dW
\end{equation*}

From now on, we use the usual norm on $E$ and usual Riesz map $R$ mapping $E$
to $E^{\prime }$. At this point, use the above implicit Ito formula \ref%
{z12febt1g} on \ref{4may5f} to obtain
\begin{equation*}
\frac{1}{2}E\left\langle \left( B+\varepsilon R\right) u\left( t\right)
,u\left( t\right) \right\rangle -\frac{1}{2}E\left\langle \left(
B+\varepsilon R\right) u_{0},u_{0}\right\rangle +E\int_{0}^{t}\left\langle
A_{\varepsilon }u,u\right\rangle ds
\end{equation*}%
\begin{eqnarray}
&=&\frac{1}{2}E\int_{0}^{t}\left\langle \left( B+\varepsilon R\right) \left(
B+\varepsilon R\right) ^{-1}B\Phi ,\left( B+\varepsilon R\right) ^{-1}B\Phi
\right\rangle ds+E\int_{0}^{t}\left\langle f,u\right\rangle ds  \notag \\
&=&\frac{1}{2}E\int_{0}^{t}\left( R^{-1}B\Phi ,\left( B+\varepsilon R\right)
^{-1}B\Phi \right) _{\mathcal{L}_{2}}ds+E\int_{0}^{t}\left\langle
f,u\right\rangle ds  \label{4may7f}
\end{eqnarray}%
where the symbol $\mathcal{L}_{2}$ signifies $\mathcal{L}_{2}\left(
Q^{1/2}U,E\right) $. Letting $\left\{ g_{i}\right\} $ be an orthonormal
basis in $Q^{1/2}U,$
\begin{equation*}
\left( R^{-1}B\Phi ,\left( B+\varepsilon R\right) ^{-1}B\Phi \right) _{%
\mathcal{L}_{2}\left( Q^{1/2}U,E\right) }\equiv \sum_{j}\left( R^{-1}B\Phi
g_{j},\left( B+\varepsilon R\right) ^{-1}B\Phi g_{j}\right) _{E}
\end{equation*}%
\begin{equation*}
=\sum_{j}\left\langle B\Phi g_{j},\left( B+\varepsilon R\right) ^{-1}B\Phi
g_{j}\right\rangle _{E^{\prime },E}=\sum_{j}\left\langle B\Phi g_{j},\left(
B+\varepsilon R\right) ^{-1}B\Phi g_{j}\right\rangle _{W^{\prime },W}
\end{equation*}%
Consider
\begin{eqnarray}
&&\sum_{j}\left\langle B\Phi g_{j},\left( B+\varepsilon R\right) ^{-1}B\Phi
g_{j}\right\rangle _{W^{\prime },W}-\sum_{j}\left\langle B\Phi g_{j},\Phi
g_{j}\right\rangle _{W^{\prime },W}  \notag \\
&=&\sum_{j}\left\langle B\Phi g_{j},\left( B+\varepsilon R\right) ^{-1}B\Phi
g_{j}-\Phi g_{j}\right\rangle  \label{4may8f}
\end{eqnarray}%
\begin{eqnarray*}
\left\langle Bh,\left( B+\varepsilon R\right) ^{-1}Bh\right\rangle
&=&\left\langle \left( B+\varepsilon R\right) \left( B+\varepsilon R\right)
^{-1}Bh,\left( B+\varepsilon R\right) ^{-1}Bh\right\rangle \\
&\geq &\left\langle B\left( B+\varepsilon R\right) ^{-1}Bh,\left(
B+\varepsilon R\right) ^{-1}Bh\right\rangle
\end{eqnarray*}%
Hence
\begin{eqnarray*}
&&\left\langle Bh,h\right\rangle ^{1/2}\left\langle B\left( B+\varepsilon
R\right) ^{-1}Bh,\left( B+\varepsilon R\right) ^{-1}Bh\right\rangle ^{1/2} \\
&\geq &\left\langle B\left( B+\varepsilon R\right) ^{-1}Bh,\left(
B+\varepsilon R\right) ^{-1}Bh\right\rangle
\end{eqnarray*}%
and so
\begin{equation*}
\left\langle Bh,h\right\rangle ^{1/2}\geq \left\langle B\left( B+\varepsilon
R\right) ^{-1}Bh,\left( B+\varepsilon R\right) ^{-1}Bh\right\rangle ^{1/2}
\end{equation*}%
Therefore,
\begin{eqnarray*}
\left\langle Bh,\left( B+\varepsilon R\right) ^{-1}Bh\right\rangle &\leq
&\left\langle Bh,h\right\rangle ^{1/2}\left\langle B\left( B+\varepsilon
R\right) ^{-1}Bh,\left( B+\varepsilon R\right) ^{-1}Bh\right\rangle ^{1/2} \\
&\leq &\left\langle Bh,h\right\rangle ^{1/2}\left\langle Bh,h\right\rangle
^{1/2}=\left\langle Bh,h\right\rangle
\end{eqnarray*}%
Therefore, \ref{4may8f} is non positive.

Return to \ref{4may7f}. The above has shown that
\begin{eqnarray*}
&&\frac{1}{2}E\int_{0}^{t}\left( R^{-1}B\Phi ,\left( B+\varepsilon R\right)
^{-1}B\Phi \right) _{\mathcal{L}_{2}}ds \\
&=&\frac{1}{2}E\int_{0}^{t}\sum_{j}\left\langle B\Phi g_{j},\left(
B+\varepsilon R\right) ^{-1}B\Phi g_{j}\right\rangle _{W^{\prime },W}
\end{eqnarray*}%
\begin{equation*}
\leq \frac{1}{2}E\int_{0}^{t}\sum_{j}\left\langle B\Phi g_{j},\Phi
g_{j}\right\rangle _{W^{\prime },W}=\frac{1}{2}E\int_{0}^{t}\left( \mathcal{R%
}^{-1}B\Phi ,\Phi \right) _{\mathcal{L}_{2}\left( Q^{1/2}U,W\right) }ds
\end{equation*}%
where $\mathcal{R}$ is the Riesz map from $W$ to $W^{\prime }$, distinct
from $R$ the Riesz map from $E$ to $E^{\prime }$. Summarizing this, the
following inequality has been established.
\begin{equation*}
\frac{1}{2}E\left\langle \left( B+\varepsilon R\right) u\left( t\right)
,u\left( t\right) \right\rangle -\frac{1}{2}E\left\langle \left(
B+\varepsilon R\right) u_{0},u_{0}\right\rangle +E\int_{0}^{t}\left\langle
A_{\varepsilon }u_{\varepsilon },u_{\varepsilon }\right\rangle ds
\end{equation*}%
\begin{equation}
\leq \frac{1}{2}E\int_{0}^{t}\left( \mathcal{R}^{-1}B\Phi ,\Phi \right) _{%
\mathcal{L}_{2}\left( Q^{1/2}U,W\right) }ds\ +E\int_{0}^{t}\left\langle
f,u\right\rangle ds\ \ \blacksquare  \label{4may9f}
\end{equation}

From now on, we will use a subscript of $\varepsilon $ on $u$ because we are
about to take limits as $\varepsilon \rightarrow 0$. From the coercivity
condition \ref{4may2f}, the following inequality is obtained.
\begin{equation*}
E\left\langle \left( B+\varepsilon R\right) u_{\varepsilon }\left( t\right)
,u_{\varepsilon }\left( t\right) \right\rangle +E\int_{0}^{t}\left\Vert
u_{\varepsilon }\right\Vert _{V}^{p}ds+\varepsilon E\int_{0}^{t}\left\Vert
u_{\varepsilon }\right\Vert _{E}^{2}ds\leq C\left( \Phi ,u_{0},f\right)
+CE\left\langle \left( B+\varepsilon R\right) u_{0},u_{0}\right\rangle
\end{equation*}%
Since $u_{0}$ is in $L^{2}\left( \Omega ,E\right) ,$ the right side is
bounded independent of $t\leq T$ and $\varepsilon $. In particular, for some
constant $C$ independent of $\varepsilon ,$
\begin{equation*}
E\left\langle Bu_{\varepsilon }\left( t\right) ,u_{\varepsilon }\left(
t\right) \right\rangle
\end{equation*}%
is bounded independent of $\varepsilon .$ Therefore, if $v\in L^{2}\left(
\Omega ,W\right) ,$%
\begin{eqnarray}
\left\vert E\left\langle Bu_{\varepsilon }\left( t\right) ,v\right\rangle
\right\vert &\leq &\left( E\left\langle Bu_{\varepsilon }\left( t\right)
,u_{\varepsilon }\left( t\right) \right\rangle \right) ^{1/2}\left(
E\left\langle Bv,v\right\rangle \right) ^{1/2}  \notag \\
&\leq &C\left\Vert B\right\Vert ^{1/2}\left\Vert v\right\Vert _{L^{2}\left(
\Omega ,W\right) }  \label{29octe6h}
\end{eqnarray}

It follows that there exists a subsequence still called $\varepsilon $ such
that
\begin{equation}
\varepsilon Ru_{\varepsilon }\left( t\right) \rightarrow 0\text{ in }%
L^{2}\left( \Omega ,E^{\prime }\right) \text{ uniformly in }t
\label{4may10f}
\end{equation}%
\begin{equation*}
u_{\varepsilon }\rightarrow u\text{ weakly in }\mathcal{V}
\end{equation*}%
\begin{equation*}
\varepsilon Ru_{\varepsilon }\rightarrow 0\text{ in }\mathcal{E}^{\prime }
\end{equation*}%
where $\mathcal{E}^{\prime }\equiv L^{q^{\prime }}\left( \left[ 0,T\right]
\times \Omega ,E^{\prime }\right) $%
\begin{equation*}
Au_{\varepsilon }\rightarrow \xi \text{ weakly in }\mathcal{V}^{\prime }
\end{equation*}%
This last convergence implies that
\begin{equation*}
\int_{0}^{t}Au_{\varepsilon }ds\rightarrow \int_{0}^{t}\xi ds\text{ weakly
in }L^{p^{\prime }}\left( \Omega ,V^{\prime }\right)
\end{equation*}

From the integral equation \ref{4may6f}, and boundedness of $A,$ it also
follows that a further subsequence satisfies
\begin{equation*}
\left( \left( B+\varepsilon R\right) u_{\varepsilon }-\left( B+\varepsilon
R\right) u_{0}-B\int_{0}^{\left( \cdot \right) }\Phi dW\right) ^{\prime
}\rightarrow \zeta \text{ weakly in }\mathcal{E}^{\prime }
\end{equation*}%
where $\mathcal{E}^{\prime }\equiv L^{q^{\prime }}\left( \left[ 0,T\right]
\times \Omega ,E^{\prime }\right) $ for $q=\max \left( p,2\right) $ and as
usual, $1/q^{\prime }+1/q=1$. Thus
\begin{equation*}
\zeta +\xi =f\text{ in }\mathcal{E}^{\prime }
\end{equation*}%
However, both $f$ and $\xi $ are in $\mathcal{V}^{\prime }$ so in fact $%
\zeta \in \mathcal{V}^{\prime }$ also from the fact that $\mathcal{E}$ is
dense in $\mathcal{V}$. Thus the equation actually holds in $\mathcal{V}%
^{\prime }$.

Consider $\zeta $. Let $g\in L^{q}\left( \Omega ,E\right) ,q=\max \left(
2,p\right) ,$ and let $\psi $ be infinitely differentiable and equal to 0
near $T$. Then since $Bu_{\varepsilon }\left( 0\right) =Bu_{0}$ a.e. $\omega
,$
\begin{equation*}
\int_{0}^{T}\int_{\Omega }\left\langle \zeta ,\psi g\right\rangle
dPdt=\lim_{\varepsilon \rightarrow 0}\int_{0}^{T}\int_{\Omega }\left\langle
\left( \left( B+\varepsilon R\right) u_{\varepsilon }-B\int_{0}^{\left(
\cdot \right) }\Phi dW-\left( B+\varepsilon R\right) u_{0}\right) ^{\prime
},\psi g\right\rangle dPdt
\end{equation*}%
Then using \ref{4may10f} and $u_{0}\in L^{2}\left( \Omega ,E\right) $,%
\begin{equation*}
=-\lim_{\varepsilon \rightarrow 0}\int_{0}^{T}\int_{\Omega }\left\langle
\left( \left( B+\varepsilon R\right) u_{\varepsilon }-B\int_{0}^{\left(
\cdot \right) }\Phi dW-\left( B+\varepsilon R\right) u_{0}\right) ,\psi
^{\prime }g\right\rangle dPdt
\end{equation*}%
\begin{equation*}
=-\lim_{\varepsilon \rightarrow 0}\int_{0}^{T}\int_{\Omega }\left\langle
\left( Bu_{\varepsilon }-B\int_{0}^{\left( \cdot \right) }\Phi
dW-Bu_{0}\right) ,\psi ^{\prime }g\right\rangle dPdt
\end{equation*}%
\begin{equation*}
=-\lim_{\varepsilon \rightarrow 0}\int_{0}^{T}\int_{\Omega }\left\langle
\psi ^{\prime }Bg,u_{\varepsilon }-\int_{0}^{\left( \cdot \right) }\Phi
dW-u_{0}\right\rangle dPdt
\end{equation*}%
\begin{eqnarray*}
&=&-\int_{0}^{T}\int_{\Omega }\left\langle \psi ^{\prime
}Bg,u-\int_{0}^{\left( \cdot \right) }\Phi dW-u_{0}\right\rangle dPdt \\
&=&-\int_{0}^{T}\int_{\Omega }\left\langle Bu-B\int_{0}^{\left( \cdot
\right) }\Phi dW-Bu_{0},\psi ^{\prime }g\right\rangle dPdt
\end{eqnarray*}%
Since $g$ is arbitrary, this shows that $\zeta =\left( Bu-B\int_{0}^{\left(
\cdot \right) }\Phi dW-Bu_{0}\right) ^{\prime }$ in $\mathcal{E}^{\prime }$.
Also, it shows that, on integrating by parts,
\begin{equation*}
\int_{0}^{T}\int_{\Omega }\left\langle \zeta ,\psi g\right\rangle
dPdt=\int_{\Omega }\left\langle Bu\left( 0\right) -Bu_{0},\psi \left(
0\right) g\right\rangle dP+
\end{equation*}%
\begin{equation*}
\int_{0}^{T}\int_{\Omega }\left\langle \left( Bu-B\int_{0}^{\left( \cdot
\right) }\Phi dW-Bu_{0}\right) ^{\prime },\psi g\right\rangle dPdt
\end{equation*}%
\begin{equation*}
=\int_{0}^{T}\int_{\Omega }\left\langle \zeta ,\psi g\right\rangle dPdt
\end{equation*}%
and so
\begin{equation*}
\int_{\Omega }\left\langle Bu\left( 0\right) -Bu_{0},\psi \left( 0\right)
g\right\rangle dP=0
\end{equation*}%
which shows that in $L^{q^{\prime }}\left( \Omega ,E^{\prime }\right) ,$ you
have
\begin{equation}
Bu\left( 0\right) =Bu_{0}.  \label{29octe1h}
\end{equation}%
By density considerations, this implies the equation also holds in $%
L^{2}\left( \Omega ,W^{\prime }\right) $. In particular, off a set of
measure zero, for all $t,$
\begin{equation*}
\int_{0}^{t}\zeta ds=Bu\left( t\right) -B\int_{0}^{t}\Phi dW-Bu_{0}
\end{equation*}

Also from \ref{4may10f}, and the above weak convergence in $\mathcal{E}%
^{\prime }$ of%
\begin{equation*}
\left( \left( B+\varepsilon R\right) u_{\varepsilon }-\left( B+\varepsilon
R\right) u_{0}-B\int_{0}^{\left( \cdot \right) }\Phi dW\right) ^{\prime },
\end{equation*}%
\begin{equation*}
Bu_{\varepsilon }\left( t\right) \rightarrow Bu\left( t\right) \text{ weakly
in }L^{q^{\prime }}\left( \Omega ,E^{\prime }\right)
\end{equation*}%
for each $t$. By \ref{29octe6h}, there is a further subsequence such that
\begin{equation*}
Bu_{\varepsilon }\left( T\right) \rightarrow Bu\left( T\right) \text{ weakly
in }L^{2}\left( \Omega ,W^{\prime }\right)
\end{equation*}%
Now if $e\in W,$ consider the functional defined on $L^{2}\left( \Omega
,W^{\prime }\right) $ which is given by
\begin{equation*}
v\rightarrow \int_{\Omega }\left\langle v,e\right\rangle ^{2}dP
\end{equation*}%
This is clearly convex and lower semicontinuous. Therefore, it is also
weakly lower semicontinuous. It follows that
\begin{equation*}
\lim \inf_{\varepsilon \rightarrow 0}\int_{\Omega }\left\langle
Bu_{\varepsilon }\left( T\right) ,e\right\rangle ^{2}dP\geq \int_{\Omega
}\left\langle Bu\left( T\right) ,e\right\rangle ^{2}dP
\end{equation*}%
Letting $\left\{ e_{i}\right\} $ be the vectors of Lemma \ref{7mayl1h}, it
follows from the above observation and Fatou's lemma,%
\begin{eqnarray}
\lim \inf_{\varepsilon \rightarrow 0}E\left\langle Bu_{\varepsilon }\left(
T\right) ,u_{\varepsilon }\left( T\right) \right\rangle &=&\lim
\inf_{\varepsilon \rightarrow 0}\sum_{i=1}^{\infty }E\left\langle
Bu_{\varepsilon }\left( T\right) ,e_{i}\right\rangle ^{2}  \notag \\
&\geq &\sum_{i=1}^{\infty }\lim \inf_{\varepsilon \rightarrow
0}E\left\langle Bu_{\varepsilon }\left( T\right) ,e_{i}\right\rangle ^{2}
\notag \\
&\geq &\sum_{i=1}^{\infty }E\left\langle Bu\left( T\right)
,e_{i}\right\rangle ^{2}=E\left\langle Bu\left( T\right) ,u\left( T\right)
\right\rangle  \label{30octe1h}
\end{eqnarray}

As explained above,
\begin{equation}
\zeta +\xi =f\text{ in }\mathcal{V}^{\prime }  \label{29octe2h}
\end{equation}%
It follows that there is a set of measure zero $N$ such that for $\omega
\notin N,$
\begin{equation*}
\zeta \left( t\right) +\xi \left( t\right) =f\left( t\right) \text{ a.e. }t,
\end{equation*}%
the equation holding in $V^{\prime }$. In addition to this, we have also
obtained
\begin{equation}
\zeta =\left( Bu-B\int_{0}^{\left( \cdot \right) }\Phi dW-Bu_{0}\right)
^{\prime }  \label{29octe3h}
\end{equation}%
Enlarging $N$ if necessary, it follows that for $\omega \notin N,$%
\begin{equation}
Bu\left( t\right) -Bu_{0}+\int_{0}^{t}\xi \left( s\right)
ds=\int_{0}^{t}fds+B\int_{0}^{t}\Phi dW  \label{30octe3h}
\end{equation}

It follows, since $A$ is progressively measurable, each $A_{\varepsilon
}u_{\varepsilon }$ is progressively measurable and so an application of the
Pettis theorem implies $\xi $ is also progressively measurable. Therefore,
we can apply the implicit Ito formula to the above integral equation of \ref%
{30octe3h} and obtain
\begin{equation*}
\frac{1}{2}E\left\langle Bu\left( t\right) ,u\left( t\right) \right\rangle -%
\frac{1}{2}E\left\langle Bu_{0},u_{0}\right\rangle
+E\int_{0}^{t}\left\langle \xi ,u\right\rangle -\frac{1}{2}\left\langle
B\Phi ,\Phi \right\rangle ds=\int_{0}^{t}\left\langle f,u\right\rangle ds
\end{equation*}%
Thus, letting $t=T,$
\begin{equation*}
E\int_{0}^{T}\left\langle \xi ,u\right\rangle ds=\int_{0}^{T}\frac{1}{2}%
\left\langle B\Phi ,\Phi \right\rangle +\left\langle f,u\right\rangle ds+%
\frac{1}{2}\left\langle Bu_{0},u_{0}\right\rangle -\frac{1}{2}E\left\langle
Bu\left( T\right) ,u\left( T\right) \right\rangle
\end{equation*}%
Recall the integral equation for the approximate solution,
\begin{equation*}
Bu_{\varepsilon }\left( t\right) -Bu_{0}+\int_{0}^{t}A_{\varepsilon }\left(
u_{\varepsilon }\right) ds=\int_{0}^{t}B\Phi dW+\int_{0}^{t}fds-\varepsilon
Ru_{\varepsilon }\left( t\right) +\varepsilon Ru_{0}
\end{equation*}%
and Lemma \ref{6augl2h}. Using the result of this lemma, when $t=T,$%
\begin{equation*}
\frac{1}{2}E\left\langle \left( B+\varepsilon R\right) u_{\varepsilon
}\left( T\right) ,u_{\varepsilon }\left( T\right) \right\rangle -\frac{1}{2}%
E\left\langle \left( B+\varepsilon R\right) u_{0},u_{0}\right\rangle
\end{equation*}%
\begin{equation*}
+E\int_{0}^{T}\left\langle A_{\varepsilon }u_{\varepsilon },u_{\varepsilon
}\right\rangle -\frac{1}{2}\left\langle B\Phi ,\Phi \right\rangle
ds=\int_{0}^{T}\left\langle f,u_{\varepsilon }\right\rangle ds
\end{equation*}%
Then, dropping the term $\left\langle \varepsilon Ru_{\varepsilon
},u_{\varepsilon }\right\rangle $ from $A_{\varepsilon },$
\begin{equation*}
E\int_{0}^{T}\left\langle Au_{\varepsilon },u_{\varepsilon }\right\rangle
ds\leq \frac{1}{2}E\int_{0}^{T}\left\langle B\Phi ,\Phi \right\rangle
+\left\langle f,u_{\varepsilon }\right\rangle ds+\frac{1}{2}E\left\langle
\left( B+\varepsilon R\right) u_{0},u_{0}\right\rangle
\end{equation*}%
\begin{equation*}
-\frac{1}{2}E\left\langle Bu_{\varepsilon }\left( T\right) ,u_{\varepsilon
}\left( T\right) \right\rangle
\end{equation*}%
Now take $\lim \sup $ of both sides and use \ref{30octe1h} to write
\begin{eqnarray*}
\lim \sup_{\varepsilon \rightarrow 0}\left( -\frac{1}{2}E\left\langle
Bu_{\varepsilon }\left( T\right) ,u_{\varepsilon }\left( T\right)
\right\rangle \right) &=&-\lim \inf_{\varepsilon \rightarrow 0}\frac{1}{2}%
E\left\langle Bu_{\varepsilon }\left( T\right) ,u_{\varepsilon }\left(
T\right) \right\rangle \\
&\leq &-\frac{1}{2}E\left\langle Bu\left( T\right) ,u\left( T\right)
\right\rangle
\end{eqnarray*}%
Thus
\begin{eqnarray*}
\lim \sup_{\varepsilon \rightarrow 0}E\int_{0}^{T}\left\langle
Au_{\varepsilon },u_{\varepsilon }\right\rangle ds &\leq &\frac{1}{2}%
E\int_{0}^{T}\left\langle B\Phi ,\Phi \right\rangle +\left\langle
f,u\right\rangle ds \\
&&+\frac{1}{2}E\left\langle Bu_{0},u_{0}\right\rangle -\frac{1}{2}%
E\left\langle Bu\left( T\right) ,u\left( T\right) \right\rangle
\end{eqnarray*}%
\begin{equation*}
=E\int_{0}^{T}\left\langle \xi ,u\right\rangle ds
\end{equation*}%
Since $A$ is assumed type $M$ on $\mathcal{V}$, it follows that $\xi =Au$.
Then referring to \ref{30octe3h}, this has proved the following theorem
which is the main result.

\begin{theorem}
\label{7mayt1f}\label{25augt1h}Let the spaces $E,V,W$ be as described in \ref%
{30octe2h} and suppose
\begin{equation*}
u\rightarrow A\left( t,u,\omega \right)
\end{equation*}%
is locally Lipschitz as a map from $E$ to $E^{\prime },$%
\begin{equation*}
\left( t,u,\omega \right) \rightarrow A\left( t,u,\omega \right)
\end{equation*}%
is progressively measurable. Also suppose that the map $A:\mathcal{%
V\rightarrow V}^{\prime }$ is type $M$ where
\begin{equation*}
\mathcal{V}\equiv L^{p}\left( \left[ 0,T\right] \times \Omega ;V\right)
\end{equation*}%
with the $\sigma $ algebra equal to $\mathcal{B}\left( \left[ 0,T\right]
\right) \times \mathcal{F}_{T}$ and there is a coercivity condition
\begin{equation*}
\left\langle A\left( t,u,\omega \right) ,u\right\rangle \geq k\left\Vert
u\right\Vert _{V}^{p}-C\left( t,\omega \right)
\end{equation*}%
where $C\in L^{1}\left( \left[ 0,T\right] \times \Omega \right) .$ Also let $%
u_{0}\in L^{q}\left( \Omega ,E\right) ,u_{0}$ being $\mathcal{F}_{0}$
measurable, where $q=\max \left( p,2\right) $ and let $f\in \mathcal{V}%
^{\prime }$. Then there exists a progressively measurable function $u\in
\mathcal{V}$ which is a solution to the integral equation
\begin{equation*}
Bu\left( t\right) -Bu_{0}+\int_{0}^{t}Au\left( s\right)
ds=\int_{0}^{t}fds+B\int_{0}^{t}\Phi dW,\ t\in \left[ 0,T\right]
\end{equation*}%
in the sense that the equation holds in $V^{\prime }$ for all $\omega \notin
N$ where $N$ is a set of measure zero. In terms of the weak derivative, this
solution is of the form
\begin{eqnarray*}
\left( Bu-B\int_{0}^{\left( \cdot \right) }\Phi dW\right) ^{\prime }+Au &=&f%
\text{ in }\mathcal{V}^{\prime } \\
Bu\left( 0\right) &=&Bu_{0}\text{ in }L^{2}\left( \Omega ,W\right)
\end{eqnarray*}
\end{theorem}

It is easy to generalize to assume only that $u_{0}\in L^{2}\left( \Omega
,W\right) .$

\begin{corollary}
\label{30octc1h}Let the spaces $E,V,W$ be as described in \ref{30octe2h} and
suppose
\begin{equation*}
u\rightarrow A\left( t,u,\omega \right)
\end{equation*}%
is locally Lipschitz as a map from $E$ to $E^{\prime },$%
\begin{equation*}
\left( t,u,\omega \right) \rightarrow A\left( t,u,\omega \right)
\end{equation*}%
is progressively measurable. Also suppose that the map $A:\mathcal{%
V\rightarrow V}^{\prime }$ is type $M$ where
\begin{equation*}
\mathcal{V}\equiv L^{p}\left( \left[ 0,T\right] \times \Omega ;V\right)
\end{equation*}%
with the $\sigma $ algebra equal to $\mathcal{B}\left( \left[ 0,T\right]
\right) \times \mathcal{F}_{T}$ and there is a coercivity condition
\begin{equation*}
\left\langle A\left( t,u,\omega \right) ,u\right\rangle \geq k\left\Vert
u\right\Vert _{V}^{p}-C\left( t,\omega \right)
\end{equation*}%
where $C\in L^{1}\left( \left[ 0,T\right] \times \Omega \right) .$ Also let $%
u_{0}\in L^{2}\left( \Omega ,W\right) ,$ $u_{0}$ being $\mathcal{F}_{0}$
measurable, and let $f\in \mathcal{V}^{\prime }$. Then there exists a
progressively measurable function $u\in \mathcal{V}$ which is a solution to
the integral equation
\begin{equation*}
Bu\left( t\right) -Bu_{0}+\int_{0}^{t}Au\left( s\right)
ds=\int_{0}^{t}fds+B\int_{0}^{t}\Phi dW,\ t\in \left[ 0,T\right]
\end{equation*}%
in the sense that the equation holds in $V^{\prime }$ for all $\omega \notin
N$ where $N$ is a set of measure zero. In terms of the weak derivative, this
solution is of the form
\begin{eqnarray*}
\left( Bu-B\int_{0}^{\left( \cdot \right) }\Phi dW\right) ^{\prime }+Au &=&f%
\text{ in }\mathcal{V}^{\prime } \\
Bu\left( 0\right) &=&Bu_{0}\text{ in }L^{2}\left( \Omega ,W\right)
\end{eqnarray*}
\end{corollary}

\textbf{Proof:\ }Let $u_{n}$ be the solution to the above theorem satisfying
the integral equation
\begin{equation*}
Bu_{n}\left( t\right)
-Bu_{0n}+\int_{0}^{t}Au_{n}ds=\int_{0}^{t}fds+B\int_{0}^{t}\Phi dW
\end{equation*}%
as described there, where $u_{0n}\in L^{q}\left( \Omega ,E\right) $ and
\begin{equation*}
u_{0n}\rightarrow u_{0}\text{ in }L^{2}\left( \Omega ,W\right) .
\end{equation*}%
Then by the implicit Ito formula,
\begin{equation*}
\frac{1}{2}E\left\langle Bu_{n}\left( t\right) ,u_{n}\left( t\right)
\right\rangle -\frac{1}{2}E\left\langle Bu_{0n},u_{0n}\right\rangle
+\int_{0}^{t}\left\langle Au_{n},u_{n}\right\rangle
\end{equation*}%
\begin{equation}
-\frac{1}{2}\left\langle B\Phi ,\Phi \right\rangle
ds=\int_{0}^{t}\left\langle f,u_{n}\right\rangle ds  \label{30octe7h}
\end{equation}%
Then, as in the above argument, there is a subsequence, still denoted by $n$
such that
\begin{equation*}
u_{n}\rightarrow u\text{ weakly in }\mathcal{V}
\end{equation*}%
\begin{equation*}
Au_{n}\rightarrow \xi \text{ weakly in }\mathcal{V}^{\prime }
\end{equation*}%
\begin{equation*}
\left( Bu_{n}-Bu_{0n}-B\int_{0}^{\left( \cdot \right) }\Phi dW\right)
^{\prime }\rightarrow \left( Bu-Bu_{0}-B\int_{0}^{\left( \cdot \right) }\Phi
dW\right) ^{\prime }\text{ in }\mathcal{V}^{\prime }
\end{equation*}%
\begin{equation*}
Bu\left( 0\right) =Bu_{0}\text{ in }L^{2}\left( \Omega ,W\right)
\end{equation*}%
As before, the integral equation implies%
\begin{equation*}
Bu_{\varepsilon }\left( t\right) \rightarrow Bu\left( t\right) \text{ weakly
in }L^{q^{\prime }}\left( \Omega ,E^{\prime }\right)
\end{equation*}%
and there is a subsequence such that also
\begin{equation}
Bu_{n}\left( T\right) \rightarrow Bu\left( T\right) \text{ weakly in }%
L^{2}\left( \Omega ,W^{\prime }\right)  \label{9nove1h}
\end{equation}%
Then passing to a limit,
\begin{equation}
Bu\left( t\right) -Bu_{0}+\int_{0}^{t}\xi
ds=\int_{0}^{t}fds+B\int_{0}^{t}\Phi dW  \label{30octe8h}
\end{equation}%
By the implicit Ito formula,
\begin{equation*}
E\int_{0}^{t}\left\langle \xi ,u\right\rangle ds=E\int_{0}^{t}\left\langle
f,u\right\rangle ds+E\int_{0}^{t}\frac{1}{2}\left\langle B\Phi ,\Phi
\right\rangle ds+\frac{1}{2}E\left\langle Bu_{0},u_{0}\right\rangle -\frac{1%
}{2}E\left\langle Bu\left( t\right) ,u\left( t\right) \right\rangle
\end{equation*}%
Then from \ref{30octe7h} and \ref{9nove1h} along with similar arguments
given in the above theorem,
\begin{equation*}
\lim \sup_{n\rightarrow \infty }\int_{0}^{T}\left\langle
Au_{n},u_{n}\right\rangle \leq \frac{1}{2}E\left\langle
Bu_{0},u_{0}\right\rangle -\frac{1}{2}E\left\langle Bu\left( T\right)
,u\left( T\right) \right\rangle
\end{equation*}%
\begin{equation*}
+\int_{0}^{T}\frac{1}{2}\left\langle B\Phi ,\Phi \right\rangle
ds+\int_{0}^{T}\left\langle f,u\right\rangle ds=\int_{0}^{T}\left\langle \xi
,u\right\rangle ds
\end{equation*}%
and so $Au=\xi $. With \ref{30octe8h}, this proves the corollary. $%
\blacksquare $

Note that there is no conclusion of uniqueness in the above theorem and
corollary.

One can make the assumptions on $\lambda B+A$ rather than $A$ and get the
same conclusions. Also, there is a uniqueness result available under an
assumption of weak monotonicity.

\begin{corollary}
\label{5septc1h}\label{6augc1h}Suppose the situation of the above corollary
but replace the coercivity, and type $M$ conditions, with the following
weaker conditions
\begin{equation}
\lambda \left\langle Bu,u\right\rangle +\left\langle A\left( t,u,\omega
\right) ,u\right\rangle _{V}\geq \delta \left\Vert u\right\Vert
_{V}^{p}-C\left( t,\omega \right)  \label{10septe4h}
\end{equation}%
for all $\lambda $ large enough where $C\in L^{1}\left( \left[ 0,T\right]
\times \Omega \right) .$ Also
\begin{equation*}
\lambda B+A:\mathcal{V\rightarrow V}^{\prime }\text{ is type }M
\end{equation*}%
Then the conclusion of Theorem \ref{25augt1h} is still valid. There exists a
progressively measurable function $u\in \mathcal{V}$ which is a solution to
the integral equation
\begin{equation}
Bu\left( t\right) -Bu_{0}+\int_{0}^{t}Au\left( s\right)
ds=\int_{0}^{t}fds+B\int_{0}^{t}\Phi dW,\ t\in \left[ 0,T\right]
\label{30octe10h}
\end{equation}%
in the sense that the equation holds in $V^{\prime }$ for all $\omega \notin
N$ where $N$ is a set of measure zero. If the weak monotonicity condition
\begin{equation*}
\left\langle \lambda Bu+A\left( t,u,\omega \right) -\left( \lambda
Bv+A\left( t,v,\omega \right) \right) ,u-v\right\rangle \geq 0
\end{equation*}%
is valid, then if $u,v$ are two solutions to \ref{30octe10h}, it follows
that off a set of measure zero, $Bu\left( t\right) =Bv\left( t\right) .$
\end{corollary}

\textbf{Proof: }Define $A_{\lambda }$ by
\begin{equation*}
\left\langle A_{\lambda }\left( t,w,\omega \right) ,v\right\rangle
_{V^{\prime },V}\equiv \left\langle e^{-\lambda t}A\left( t,e^{\lambda
t}w,\omega \right) ,v\right\rangle _{V^{\prime },V}
\end{equation*}%
Then
\begin{equation*}
\lambda \left\langle Bu,u\right\rangle +\left\langle A_{\lambda }\left(
t,u,\omega \right) ,u\right\rangle _{V}\geq e^{-2\lambda t}\left( \lambda
\left\langle B\left( e^{\lambda t}u\right) ,e^{\lambda t}u\right\rangle
+\left\langle A\left( t,e^{\lambda t}u,\omega \right) ,e^{\lambda
t}u\right\rangle \right)
\end{equation*}%
\begin{equation*}
\geq e^{-2\lambda t}\left( \delta \left\Vert e^{\lambda t}u\right\Vert
_{V}^{p}-C\left( t,\omega \right) \right) \geq e^{-2\lambda t}\left( \delta
\left\Vert e^{\lambda t}u\right\Vert _{V}^{p}-e^{\lambda pt}e^{-\lambda
pt}C\left( t,\omega \right) \right)
\end{equation*}%
\begin{equation*}
\geq e^{-2\lambda t}e^{p\lambda t}\left( \delta \left\Vert u\right\Vert
_{V}^{p}-e^{-\lambda pt}C\left( t,\omega \right) \right) \geq \bar{\delta}%
\left\Vert u\right\Vert _{V}^{p}-\bar{C}\left( t,\omega \right)
\end{equation*}%
which is of the right form. By Corollary \ref{30octc1h}, there exists a
solution in $\mathcal{V}^{\prime }$ to
\begin{equation*}
\left( Bu-Be^{-\lambda \left( \cdot \right) }\int_{0}^{\left( \cdot \right)
}\Phi dW\right) ^{\prime }+\lambda Bu+A_{\lambda }u=e^{-\lambda \left( \cdot
\right) }f+\lambda e^{-\lambda \left( \cdot \right) }B\int_{0}^{\left( \cdot
\right) }\Phi dW\text{ in }\mathcal{V}^{\prime }
\end{equation*}%
\begin{equation*}
Bu\left( 0\right) =Bu_{0}\text{ in }L^{2}\left( \Omega ,W\right)
\end{equation*}%
Now let $e^{\lambda t}u\left( t\right) =w\left( t\right) $. Writing in terms
of $w,$%
\begin{equation*}
\left( Be^{-\lambda \left( \cdot \right) }w-Be^{-\lambda \left( \cdot
\right) }\int_{0}^{\left( \cdot \right) }\Phi dW\right) ^{\prime }+\lambda
Be^{-\lambda \left( \cdot \right) }w+e^{-\lambda \left( \cdot \right) }Aw
\end{equation*}%
\begin{equation*}
=e^{-\lambda \left( \cdot \right) }f+\lambda e^{-\lambda \left( \cdot
\right) }B\int_{0}^{\left( \cdot \right) }\Phi dW\text{ in }\mathcal{V}%
^{\prime }
\end{equation*}%
\begin{equation*}
Bw\left( 0\right) =Bu_{0}\text{ in }L^{2}\left( \Omega ,W\right)
\end{equation*}%
It follows that
\begin{equation*}
e^{-\lambda \left( \cdot \right) }\left( Bw-B\int_{0}^{\left( \cdot \right)
}\Phi dW\right) ^{\prime }-\lambda e^{-\lambda \left( \cdot \right) }\left(
Bw-B\int_{0}^{\left( \cdot \right) }\Phi dW\right)
\end{equation*}%
\begin{equation*}
+\lambda Be^{-\lambda \left( \cdot \right) }w+e^{-\lambda \left( \cdot
\right) }Aw=e^{-\lambda \left( \cdot \right) }f+\lambda e^{-\lambda \left(
\cdot \right) }B\int_{0}^{\left( \cdot \right) }\Phi dW
\end{equation*}%
After cancelling terms and multiplying by $e^{-\lambda \left( \cdot \right)
},$ this yields%
\begin{equation*}
\left( Bw-B\int_{0}^{\left( \cdot \right) }\Phi dW\right) ^{\prime }+Aw=f%
\text{ in }\mathcal{V}^{\prime }
\end{equation*}%
along with the initial condition
\begin{equation*}
Bw\left( 0\right) =Bu_{0}
\end{equation*}%
Then integrating, one obtains that for $\omega $ not in a suitable set of
measure zero,%
\begin{equation*}
Bw\left( t\right) -Bu_{0}-B\int_{0}^{t}\Phi
dW+\int_{0}^{t}Awds=\int_{0}^{t}fds\
\end{equation*}

Next suppose $u,v$ are two solutions to the above integral equation as
described above. Then off the union of the two exceptional sets,
\begin{equation*}
B\left( u\left( t\right) -v\left( t\right) \right) +\int_{0}^{t}Au-Avds=0
\end{equation*}%
and so by the implicit Ito formula,
\begin{equation*}
\frac{1}{2}\left\langle Bu\left( t\right) -v\left( t\right) ,u\left(
t\right) -v\left( t\right) \right\rangle +\int_{0}^{t}\left\langle
Au-Av,u-v\right\rangle ds=0
\end{equation*}%
and so from the monotonicity condition,
\begin{equation*}
\frac{1}{2}\left\langle Bu\left( t\right) -v\left( t\right) ,u\left(
t\right) -v\left( t\right) \right\rangle \leq \lambda
\int_{0}^{t}\left\langle B\left( u-v\right) ,u-v\right\rangle ds
\end{equation*}%
Apply Gronwall's inequality. $\blacksquare $

There is an assumption that $A:\mathcal{V\rightarrow V}^{\prime }$ is type $%
M $. This is made because this is what will be used. Recall that monotone
and hemicontinuous implies type $M$ \cite{sho96}. What are some easy to
check conditions which will imply $A$ is hemicontinuous on $\mathcal{V}$?
Make the following specific assumption involving an estimate.

There exists $c\in \lbrack 0,\infty )$ and $g\in L^{p^{\prime }}\left( \left[
0,T\right] \times \Omega \right) $ which is $\mathcal{B}\left( \left[ 0,T%
\right] \right) \times \mathcal{F}$ measurable such that for all $v\in
V,t\in \left[ 0,T\right] $%
\begin{equation}
\left\vert \left\vert A\left( t,v\right) \right\vert \right\vert _{V^{\prime
}}\leq g\left( t\right) +c\left\vert \left\vert v\right\vert \right\vert
_{V}^{p-1}.  \label{z26april31f}
\end{equation}%
Here $1/p+1/p^{\prime }=1.$

If $\lim_{\lambda \rightarrow 0}\left\langle A\left( t,u+\lambda v,\omega
\right) ,w\right\rangle =\left\langle A\left( t,u,\omega \right)
,w\right\rangle ,$ then it will also follow that $A$ is hemicontinuous on $%
\mathcal{V}$. You need to verify that
\begin{equation*}
\lim_{\lambda \rightarrow 0}\int_{\Omega }\int_{0}^{T}\left\langle A\left(
t,u+\lambda v,\omega \right) ,w\right\rangle dtdP=\int_{\Omega
}\int_{0}^{T}\left\langle A\left( t,u,\omega \right) ,w\right\rangle dtdP
\end{equation*}%
for $w,u,v\in \mathcal{V}$. But this follows from the growth condition
above, Vitali's convergence theorem, and showing that the integrand is
uniformly integrable. Letting $Q\subseteq \Omega \times \left[ 0,T\right] ,$
\begin{equation*}
\int_{Q}\left\vert \left\langle A\left( t,u+\lambda v,\omega \right)
,w\right\rangle \right\vert dtdP\leq \int_{Q}\left( \left\vert g\left(
t,\omega \right) \right\vert +c\left\vert \left\vert u+\lambda v\right\vert
\right\vert _{V}^{p-1}\right) \left\Vert w\right\Vert _{V}dtdP
\end{equation*}%
\begin{eqnarray*}
&\leq &\left( \int_{Q}\left\Vert g\right\Vert ^{p^{\prime }}dtdP\right)
^{1/p^{\prime }}\left( \int_{Q}\left\Vert w\right\Vert _{V}^{p}dtdP\right)
^{1/p} \\
&&+c\left( \int_{Q}\left\Vert u+\lambda v\right\Vert ^{p}dtdP\right)
^{1/p^{\prime }}\left( \int_{Q}\left\Vert w\right\Vert _{V}^{p}dtdP\right)
^{1/p}
\end{eqnarray*}%
which is small independent of $\lambda $ for $\left\vert \lambda \right\vert
<1$ provided $Q$ has sufficiently small measure.

It would be very interesting to get interesting examples where $\lambda B+A$
is type $M$ on $\mathcal{V}$ without being monotone for any $\lambda $. It
would also be interesting to obtain theorems which involve $\lambda B+A$
being type $M$ on $L^{p}\left( 0,T,V\right) $ rather than $\mathcal{V}$ or
on a suitable space of solutions. The difficulty in doing this latter
problem is obtaining the appropriate measurability which seems to require
the need for the solution to the non stochastic evolution problem to have a
unique solution. Of course this is usually achieved by having some
combination of the operators $B,A$ being monotone. However, one can obtain
more general conditions which do include the case of non monotone operators
by using the theory presented here as a starting point and adding various
nonlinear operators as compact perturbations. This will be explored later.

\section{Some examples\label{examples}}

Here we give a few examples. The first is a standard example, the porous
media equation, which is discussed well in \cite{sho96}. For stochastic
versions of this example, see \cite{roc07} and \cite{pre07}. The
generalization to stochastic equations does not require the theory of this
paper. We will show, however, that it can be considered in terms of the
theory of this paper without much difficulty using an approach proposed in
\cite{bre68}.

\begin{example}
The stochastic porous media equation is
\begin{equation*}
u_{t}-\Delta \left( u\left\vert u\right\vert ^{p-2}\right) =f,\ u\left(
0\right) =u_{0},\ u=0\text{ on }\partial U
\end{equation*}%
where here $U$ is a bounded open set in $\mathbb{R}^{n},n\leq 3$ having
Lipschitz boundary. One can consider a stochastic version of this as a
solution to the following integral equation
\begin{equation}
u\left( t\right) -u_{0}+\int_{0}^{t}\left( -\Delta \right) \left(
u\left\vert u\right\vert ^{p-2}\right) ds=\int_{0}^{t}\Phi dW  \label{7may1f}
\end{equation}%
where here $\Phi \in L^{2}\left( \left[ 0,T\right] \times \Omega ,\mathcal{L}%
_{2}\left( Q^{1/2}U,H\right) \right) $, $H=L^{2}\left( U\right) $ and the
equation holds in the manner described above in $H^{-1}\left( U\right) .$
Assume $p\geq 2$.
\end{example}

One can consider this as an implicit integral equation of the form
\begin{equation}
\left( -\Delta \right) ^{-1}u\left( t\right) -\left( -\Delta \right)
^{-1}u_{0}+\int_{0}^{t}u\left\vert u\right\vert ^{p-2}ds=\left( -\Delta
\right) ^{-1}\int_{0}^{t}\Phi dW  \label{7may2f}
\end{equation}%
where $-\Delta $ is the Riesz map of $H_{0}^{1}\left( U\right) $ to $%
H^{-1}\left( U\right) $. Then we can also consider $\left( -\Delta \right)
^{-1}$ as a map from $L^{2}\left( U\right) $ to $L^{2}\left( U\right) $ as
follows.
\begin{equation*}
\left( -\Delta \right) ^{-1}f=u\text{ where }-\Delta u=f,\ u=0\text{ on }%
\partial U.
\end{equation*}%
Thus we let $W=L^{2}\left( U\right) $ and $V=L^{p}\left( U\right) .$ Let $%
B\equiv \left( -\Delta \right) ^{-1}$ on $L^{2}\left( U\right) $ as just
described. Let $A\left( u\right) =u\left\vert u\right\vert ^{p-2}$. Then
clearly there exists a Hilbert space $E$ dense in $V$ which has the property
that $A:E\rightarrow E^{\prime }$ is locally Lipschitz. Just pick $%
E=H^{2}\left( U\right) $ for example. By the Sobolev embedding theorem, $%
H^{2}\left( U\right) $ embeds continuously into $C\left( \overline{U}\right)
$. Furthermore, the function $F\left( x\right) =x\left\vert x\right\vert
^{p-2}$ is differentiable as a map from $\mathbb{R}$ to $\mathbb{R}$. Thus
for $u,v\in E,$
\begin{equation*}
F\left( u+v\right) =F\left( u\right) +F^{\prime }\left( u\right) v+o\left(
v\right)
\end{equation*}%
If $w\in E,$ then if $\left\Vert v\right\Vert _{H^{2}\left( U\right) }$ is
small enough, it follows that $\left\Vert v\right\Vert _{C\left( \bar{U}%
\right) }$ is also small enough that
\begin{equation*}
\left\vert \left\langle o\left( v\right) ,w\right\rangle \right\vert \leq
\varepsilon \left\Vert w\right\Vert _{C\left( \bar{U}\right) }\leq
\varepsilon \left\Vert w\right\Vert _{H^{2}}
\end{equation*}%
which shows that this function is differentiable as a map from $E$ to $%
E^{\prime }$. Similarly, the derivative is continuous. Hence it is locally
Lipschitz as a map from $E$ to $E^{\prime }$. Thus $A:E\rightarrow E^{\prime
}$ is locally Lipschitz on $E.$ It is obvious that the necessary coercivity
condition holds. In addition, there is a monotonicity condition which holds.
In fact, $A$ is monotone and hemicontinuous so it also is type $M$.
Therefore, if $u_{0}\in L^{2}\left( \Omega ,L^{2}\left( U\right) \right) $
and $\mathcal{F}_{0}$ measurable, Corollary \ref{6augc1h} applies and we can
conclude that there exists a unique solution to the integral equation \ref%
{7may2f} in the sense described there. Here $u\in L^{p}\left( \left[ 0,T%
\right] \times \Omega ,L^{p}\left( U\right) \right) $ and is progressively
measurable, the integral equation holding for all $t$ for $a.e.\omega $.
Since $A$ satisfies for some $\delta >0$ an inequality of the form
\begin{equation*}
\left\langle Au-Av,u-v\right\rangle \geq \delta \left\Vert u-v\right\Vert
_{L^{p}\left( U\right) }^{p},
\end{equation*}%
it follows easily from the above methods that the solution is also unique.

Next we give a simple example which is a singular and degenerate equation.

\begin{example}
Suppose $U$ is a bounded open set in $\mathbb{R}^{3}$ and $b\left( \mathbf{x}%
\right) \geq 0,$ $b\in L^{p}\left( U\right) ,$ $p\geq 4$ for simplicity.
Consider the degenerate stochastic initial boundary value problem%
\begin{eqnarray*}
b\left( \mathbf{\cdot }\right) u\left( t,\mathbf{\cdot }\right) -b\left(
\mathbf{\cdot }\right) u_{0}\left( \mathbf{\cdot }\right)
-\int_{0}^{t}\nabla \cdot \left( \left\vert \nabla u\right\vert ^{p-2}\nabla
u\right) &=&b\int_{0}^{t}\Phi dW \\
u &=&0\text{ on }\partial U
\end{eqnarray*}%
where $\Phi \in L^{2}\left( \left[ 0,T\right] \times \Omega ,\mathcal{L}%
_{2}\left( Q^{1/2}U,W\right) \right) $ for $W=H_{0}^{1}\left( U\right) .$
\end{example}

To consider this equation and initial condition, it suffices to let $%
W=H_{0}^{1}\left( U\right) ,V=W_{0}^{1,p}\left( U\right) ,$
\begin{eqnarray*}
A &:&V\rightarrow V^{\prime },\ \left\langle Au,v\right\rangle
=\int_{U}\left\vert \nabla u\right\vert ^{p-2}\nabla u\cdot \nabla vdx, \\
\ B &:&W\rightarrow W^{\prime },\left\langle Bu,v\right\rangle
=\int_{U}b\left( x\right) u\left( x\right) v\left( x\right) dx
\end{eqnarray*}%
Then by the Sobolev embedding theorem, $B$ is obviously self adjoint,
bounded and nonnegative. This follows from a short computation:
\begin{equation*}
\left\vert \int_{U}b\left( x\right) u\left( x\right) v\left( x\right)
dx\right\vert \leq \left\Vert v\right\Vert _{L^{4}\left( U\right) }\left(
\int_{U}\left\vert b\left( x\right) \right\vert ^{4/3}\left\vert u\left(
x\right) \right\vert ^{4/3}dx\right) ^{3/4}
\end{equation*}%
\begin{eqnarray*}
&\leq &\left\Vert v\right\Vert _{H_{0}^{1}\left( U\right) }\left( \left(
\int \left\vert b\left( x\right) \right\vert ^{4}dx\right) ^{1/3}\left( \int
\left( \left\vert u\left( x\right) \right\vert ^{4/3}\right) ^{3/2}\right)
^{2/3}\right) ^{3/4} \\
&=&\left\Vert v\right\Vert _{H_{0}^{1}\left( U\right) }\left\Vert
b\right\Vert _{L^{4}\left( U\right) }\left\Vert u\right\Vert _{L^{2}\left(
U\right) }\leq C\left\Vert b\right\Vert _{L^{4}}\left\Vert u\right\Vert
_{H_{0}^{1}}\left\Vert v\right\Vert _{H_{0}^{1}}
\end{eqnarray*}%
The nonlinear operator is obviously monotone and hemicontinuous, so it is
pseudomonotone. The technical requirement that the operator $A$ be locally
Lipschitz on some Hilbert space $E$ which is dense in $V$ and embeds
continuously into $V$ is easily satisfied by taking $E=H^{m+1}\left(
U\right) $ for $2m>3$ then using the Sobolev embedding theorem, similar to
the above example. As for $u_{0},$ it is only necessary to assume $u_{0}\in
L^{2}\left( \Omega ,W\right) $ and $\mathcal{F}_{0}$ measurable. Then
Corollary \ref{6augc1h} gives the existence of a solution. Note that $b$ can
be unbounded and may also vanish. Thus the equation can degenerate to the
case of a non stochastic nonlinear elliptic equation.

The existence theorems can easily be extended to include the situation where
$\Phi $ is replaced with a function of the unknown function $u$. This is
done by splitting the time interval into small sub intervals of length $h$
and retarding the function in the stochastic integral, like a standard proof
of the Peano existence theorem. Then the Ito formula is applied to obtain
estimates and a limit is taken. However, this will be done later.

\end{document}